\documentclass[12pt,american]{article}
\usepackage[T1]{fontenc}
\usepackage[utf8]{inputenc}
\usepackage{amsmath}
\usepackage{amsthm}
\usepackage{amssymb}
\usepackage{geometry}
\usepackage{setspace}
\usepackage[authoryear]{natbib}
\makeatletter
\theoremstyle{plain}
\newtheorem{thm}{\protect\theoremname}
\theoremstyle{plain}
\newtheorem{lem}[thm]{\protect\lemmaname}
\theoremstyle{remark}
\newtheorem{rem}[thm]{\protect\remarkname}
\theoremstyle{plain}
\newtheorem{cor}[thm]{\protect\corollaryname}
\theoremstyle{plain}
\newtheorem{prop}[thm]{\protect\propositionname}
\theoremstyle{plain}
\newtheorem{assumption}{Assumption}

\@ifundefined{date}{}{\date{}}
\usepackage{babel}
\providecommand{\corollaryname}{Corollary}
\providecommand{\lemmaname}{Lemma}
\providecommand{\theoremname}{Theorem}

\providecommand{\remarkname}{Remark}

\providecommand{\propositionname}{Proposition}

\makeatother

\usepackage{babel}
\providecommand{\corollaryname}{Corollary}
\providecommand{\lemmaname}{Lemma}
\providecommand{\propositionname}{Proposition}
\providecommand{\remarkname}{Remark}
\providecommand{\theoremname}{Theorem}

\begin{document}
\title{Approximation rates for finite mixtures of location-scale models and
fast least-squares estimators}
\author{Hien Duy Nguyen\thanks{School of Computing, Engineering and Mathematical Sciences, La Trobe
University, Bundoora, VIC 3086, Australia; Institute of Mathematics
for Industry, Kyushu University, Nishi Ward, Fukuoka 819-0395, Japan.
Email: h.nguyen5@latrobe.edu.au}\and TrungTin Nguyen\thanks{School of Mathematical Sciences, Faculty of Science, Queensland University
of Technology, Brisbane, QLD 4000, Australia.}\and Jacob Westerhout\thanks{School of Mathematics and Physics, The University of Queensland, St
Lucia, QLD 4072, Australia.}\and Xin Guo\thanks{School of Mathematics and Physics, The University of Queensland, St
Lucia, QLD 4072, Australia.}}
\maketitle
\begin{abstract}
Finite mixture models provide a flexible framework for approximating
and estimating multivariate probability densities. We study mixtures
formed from translated and rescaled copies of a fixed density kernel
and obtain explicit results for both approximation and least-squares
estimation. Our main deterministic result is a quantisation theorem
showing that, after smoothing the target density at a fixed resolution,
the resulting convolution can be compressed into a finite location
mixture with controlled error. Combining this with the smoothing bias
yields approximation rates in $\mathcal{L}_{p}$ over Sobolev classes.
For estimation, we analyse least-squares $\varepsilon$-minimisers
over suitably tuned mixture sieves. Under exponential decay of the
Fourier transform of the kernel, a matching moment condition,
and bounded Sobolev targets, the estimator attains a squared $\mathcal{L}_{2}$
risk bound whose rate matches the Sobolev minimax benchmark up to
a logarithmic factor. If, in addition, the kernel is bandlimited,
then the same theorem recovers the Sobolev rate $n^{-2s/\left(2s+d\right)}$.
We further report a slower convergence rate under weaker VC-type assumptions.
At fixed scale, the Fourier-based approach also gives a nearly parametric
risk bound for the associated location-mixture class, and the same
bandlimited simplification removes the logarithmic correction. In
the Gaussian case, this recovers the known Gaussian location-mixture
rate. We also prove matching lower
bounds on Gaussian convolution submodels, including strict submodels
of the Gaussian location-mixture class, and on the tensor-product odd-degree
Student-$t$ location-mixture family. 
\end{abstract}

\section{Introduction}

\label{sec:Introduction}

Finite mixture models appear ubiquitously throughout the statistics
literature, and their applications and theory are developed extensively
in classic works such as \citet{EverittHand1981FiniteMixture}, \citet{Titterington1985StatisticalAnalysis},
\citet{McLachlanBasford1988MixtureModels}, \citet{Lindsay1995MixtureModels},
\citet{mclachlan2000finite}, and \citet{Fruhwirth-Schnatter2006FiniteMixture},
along with more recent texts such as \citet{McNicholas2015MixtureModel},
\citet{Ng2019MixtureModelling}, \citet{Robert2019HandbookMixture},
\citet{Chen2023StatisticalInference}, and \citet{Yao2024MixtureModels},
among other sources. In application, finite mixtures are typically
employed to perform one of two primary tasks: to provide a statistical
model for a population that is hypothesized or known to have a heterogeneous
subpopulation structure, or to provide a flexible approximation to
some unknown probability distribution. When finite mixtures are employed
for the latter task, the Gaussian mixture model is typically considered
for its well-known ability to provide accurate approximations to the
probability models commonly encountered in applications (cf. \citealp{nguyen2019approximations}
and references therein).

On $\mathbb{R}^{d}$, the fact that linear combinations of Gaussian
kernels can approximate sufficiently regular functions in ${\cal L}_{p}$
spaces goes back to the Tauberian theory of \citet{wiener1932tauberian}
(cf. \citealp[Ch. VI]{Achieser1956TheoryApproximation}). Subsequent
generalizations and developments on such results appear throughout
the approximation theory literature, and include the works of \citet{park1991universal},
\citet{park1993approximation}, \citet{Xu1993Constructive}, \citet{sandberg2001gaussian},
and \citet{nestoridis2011universal}---see also \citet[Ch. 20--24]{CheneyLight2000Approximation}.

When we restrict the approximand $f_{0}:\mathbb{R}^{d}\to\mathbb{R}_{\ge0}$
to be a probability density function, and the approximator class (${\cal A}$,
say) to be convex combinations of probability density kernels, the
approximation problem requires different tools. When $f_{0}$ is within
the closure of ${\cal A}$, and the class ${\cal A}$ is uniformly
bounded in ${\cal L}_{p}$ norm, a Maurey--Barron--Jones type result
of \citet{donahue1997rates} provides an approximation theorem with
known and optimal rate (see, e.g., \citealp{zeevi1997density} and
\citealp{nguyen2019approximations}, for example). However, when ${\cal A}$
is unbounded in ${\cal L}_{p}$ norm, the study of approximation demands
more care.

In \citet{bacharoglou2010approximation}, it is proved that, on $\mathbb{R}$,
one can arbitrarily well approximate $f_{0}$ in any ${\cal L}_{p}$
norm ($p\in\left[1,\infty\right]$), under integrability and tail
conditions, by taking ${\cal A}$ to be the class of convex combinations
of Gaussian densities on $\mathbb{R}$. Then, in \citet{nguyen2020approximation}
and \citet{nguyen2023approximation}, generalizations are obtained,
where it is proved that on $\mathbb{R}^{d}$, $f_{0}$ can be approximated
arbitrarily well when ${\cal A}$ is taken to be location-scale mixtures
of some fixed probability density function, under various regularity
conditions.

We extend the qualitative results of \citet{nguyen2020approximation}
and \citet{nguyen2023approximation} by proving quantitative ${\cal L}_{p}$
approximation rates in the number of mixture components. More precisely,
when approximation is carried out by location-scale mixtures and the
target $f_{0}$ satisfies suitable smoothness assumptions, we obtain
explicit rates in the size of the approximating convex combination.
In particular, for classes of targets $f_{0}$ in the Sobolev space
${\cal W}^{1,p}\left(\mathbb{R}^{d}\right)$ \citep{Leoni2017SobolevSpaces}
or the fractional Sobolev spaces ${\cal W}^{s,p}\left(\mathbb{R}^{d}\right)$
($s\in\left(0,1\right)$; \citealp{Leoni2023FractionalSobolev}) with
bounded norms, $f_{0}$ has smoothness coefficient $\alpha=1$ or
$\alpha=s$, respectively. Then, for every $f_{0}$ whose Sobolev
(or fractional Sobolev) norm is bounded, if ${\cal A}$ and $f_{0}$
satisfy appropriate integrability assumptions, then there exist sequences
of location-scale mixtures $\left(f_{m}\right)_{m}\subset{\cal A}$
such that 
\[
\left\Vert f_{m}-f_{0}\right\Vert _{p}\lesssim\begin{cases}
m^{-\frac{\alpha}{\alpha q+d}} & \text{if }1<p<2\text{,}\\
m^{-\frac{\alpha q}{2\left(\alpha q+d\right)}} & \text{if }2\le p<\infty\text{,}
\end{cases}
\]
where $1/p+1/q=1$.

As an application of our approximation result, we study empirical
least-squares estimation over fixed-scale mixture sieves generated
by kernels whose Fourier transforms decay exponentially. The deterministic
approximation theory yields the bias term $\nu^{d}/m+\nu^{-2s}$,
and we then prove a localized empirical-process bound that exploits
this Fourier decay. Under a supersmoothness condition (cf. \citealt[Sec.~2.4.3]{meister2009deconvolution}), together
with the corresponding moment hypothesis on the kernel,
this yields a fast oracle inequality for least-squares $\epsilon$-minimizers
with a logarithmic correction and, after tuning the sieve, an expected
squared ${\cal L}_{2}$ risk of order 
\[
\left(\frac{\left(1+\log n\right)^{d/\beta}}{n}\right)^{2s/\left(2s+d\right)}
\]
for bounded targets $f_{0}\in{\cal W}^{s,2}\cap{\cal L}_{\infty}$,
where $\beta$ denotes the Fourier-decay exponent. Since the classical
minimax benchmark over Sobolev balls is of order $n^{-2s/\left(2s+d\right)}$
(cf. \citealp[Cor. 11.11]{Klemela2009Smoothing} and \citealp[Exercise 6.3.2]{GineNickl2016}),
the rate appearing here is that benchmark multiplied by a logarithmic correction
factor. If the kernel is additionally bandlimited, then the logarithmic
correction disappears from the localized modulus and the bandlimited
part of Theorem~\ref{thm:FastGaussianLeastSquares} recovers the
classical Sobolev rate $n^{-2s/\left(2s+d\right)}$. The same localized
Fourier argument also yields a fixed-scale theorem of order $\nu^{d}\left(1+\log n\right)^{d/\beta}/n$
over the associated location-mixture class at scale $\nu$, and in
the bandlimited case the logarithmic factor again disappears. The
Gaussian kernel corresponds to $\beta=2$; taking $\nu=1$ then gives
the near-parametric order $\left(1+\log n\right)^{d/2}/n$ on the
Gaussian location-mixture model itself, matching the minimax order
of \citet{kim2022minimax} and providing a fast analogue of \citet[Thm. 2]{klemela2007density}.
For the tensor-product odd-degree Student-$t$ family on $\mathbb{R}^{d}$,
the Fourier transform still has exponential decay with $\beta=1$,
and we prove the matching lower bound of order 
\[
\nu^{d}\frac{\left(\log n\right)^{d}}{n}
\]
at fixed scale. We further show that the Gaussian lower-bound order persists on Gaussian
convolution families of models in which the kernel can be written as
a Gaussian convolved with a probability measure whose Fourier transform
is bounded away from zero, including strict submodels of the Gaussian
location-mixture class.

Before we proceed, we note that the
literature on density estimation and approximation by finite mixture
models is far broader and richer than what has been cited. An incomplete
list of such works is as follows. For instance, the results of \citet{nguyen2020approximation}
and \citet{nguyen2023approximation} have been extended to the mixture
of experts setting of \citet{jacobs1991adaptive} in \citet{nguyen2021approximations}
and to the spherical manifold setting in \citet{ng2022universal}.
Schemes relating to weak approximations by finite mixture models appear,
for example, in \citet{Dalal1983Approximating}, \citet{petrone2010feller},
and \citet{Lee2012Modeling}, while best approximation and estimation
results in various divergence measures are obtained by \citet{Ma2025BestApproximation}
for the case where $f_{0}$ is in the closure of the Gaussian mixture
models. In the Bayesian finite-mixture literature, \citet{norets2012bayesian}
characterize the Kullback--Leibler closure of finite mixtures of
multivariate normals and establish consistency of the induced joint
and conditional predictive densities. Within the stagewise density-estimation
literature initiated by \citet{klemela2007density}, \citet{naito2013density}
extend the procedure to general $U$-divergence losses. Joint approximation
and estimation results in related Kullback--Leibler settings have
also been considered by \citet{Li1999Mixture}, \citet{Rakhlin2005RiskBounds},
and \citet{Dalalyan2018Optimal}, while \citet{ChongRiskBounds} study
a modified Kullback--Leibler divergence.

The rest of the manuscript proceeds as follows. In Section~\ref{sec:Preliminaries},
we outline the notational conventions of the manuscript and provide
some preliminary setup. In Section~\ref{sec:Approximation-rates-for},
we isolate the fixed-scale quantization lemma and derive the deterministic
approximation rates as corollaries. Section~\ref{sec:Finite-Gaussian-mixture-estimation}
then studies least-squares estimation over finite mixture sieves:
it proves a fast least-squares theorem for kernels with exponentially
decaying Fourier transform and reports a slower least-squares benchmark
under weaker assumptions for more general kernel classes. Section~\ref{sec:LocationMixtureTargets}
turns to fixed-scale location-mixture targets: it establishes the
corresponding fixed-scale risk bound over the associated location-mixture
class, with the Gaussian case recovering the minimax order from \citet{kim2022minimax}
on the Gaussian location-mixture model, and proves matching lower bounds on a Gaussian convolution
family of submodels of that Gaussian class and on the tensor-product
odd-degree Student-$t$ location-mixture family. A discussion of the
results appears in Section~\ref{sec:Discussion}. The Appendix then
collects the auxiliary results used throughout the paper and the more
technical proofs deferred from the main text.

\section{Preliminaries}

\label{sec:Preliminaries}

\paragraph{General notation and conventions}

For $m\in\mathbb{N}$, we write $\left[m\right]=\left\{ 1,\dots,m\right\} $.
For nonnegative quantities $a$ and $b$, we write $a\lesssim b$
if $a\le Cb$ for some finite constant $C>0$ independent of the principal
asymptotic variables under discussion, and $a\asymp b$ if both $a\lesssim b$
and $b\lesssim a$. We write $B\left(x,r\right)=\left\{ y\in\mathbb{R}^{d}\text{ | }\left\Vert y-x\right\Vert _{2}\le r\right\} $
for the closed Euclidean ball of radius $r>0$ centered at $x$, and
$\log_{+}\left(t\right)=\max\left\{ \log t,0\right\} $ for $t>0$.
If $p\in\left(1,\infty\right)$, its Hölder conjugate is denoted by
$q=p/\left(p-1\right)$, so $1/p+1/q=1$. Unnamed positive constants
may change from line to line.

\paragraph{${\cal L}_{p}$ spaces}

For each $p\in\left[1,\infty\right)$, say that $f\in{\cal L}_{p}\left(\mathbb{R}^{d}\right)$
if $\left\Vert f\right\Vert _{p}<\infty$, where 
\[
\left\Vert f\right\Vert _{p}=\left(\int_{\mathbb{R}^{d}}\left|f\left(x\right)\right|^{p}\mathrm{d}x\right)^{1/p}\text{.}
\]
When $p=\infty$, we say that $f\in{\cal L}_{\infty}\left(\mathbb{R}^{d}\right)$
if its essential supremum satisfies $\left\Vert f\right\Vert _{\infty}<\infty$.
Here $\lambda$ denotes the Lebesgue measure whenever the ambient
measure needs to be named explicitly. When there is no confusion,
we simplify ${\cal L}_{p}\left(\mathbb{R}^{d}\right)$ to ${\cal L}_{p}$,
along with other such notation. For vectors $x\in\mathbb{R}^{d}$,
we also write $\left\Vert x\right\Vert _{p}$ for the usual $\ell_{p}$
norm unless that would be misleading. Equalities involving elements
of ${\cal L}_{p}$ are understood up to Lebesgue-null sets. Whenever
an ${\cal L}_{p}$ object is evaluated pointwise, for instance in an
empirical criterion or in a pointwise mixture class, a specified measurable
representative is being used.

\paragraph{Sobolev spaces}

We write ${\cal C}_{c}^{\infty}\left(\mathbb{R}^{d}\right)$ for the
space of smooth compactly supported functions on $\mathbb{R}^{d}$,
and ${\cal S}\left(\mathbb{R}^{d}\right)$ for the Schwartz space
of smooth functions whose derivatives decay faster than every polynomial
at infinity (cf. \citealp[Def. 7.1.2]{Hormander2003}).

For each $p\in\left[1,\infty\right]$, we say that $f$ admits first-order
weak derivatives in ${\cal L}_{p}$ if for each $k\in\left[d\right]$,
there exists a $v_{k}\in{\cal L}_{p}$, such that 
\[
\int_{\mathbb{R}^{d}}f\left(x\right)D_{k}u\left(x\right)\mathrm{d}x=-\int_{\mathbb{R}^{d}}v_{k}\left(x\right)u\left(x\right)\mathrm{d}x\text{,}
\]
for all smooth and compactly supported functions $u\in{\cal C}_{c}^{\infty}\left(\mathbb{R}^{d}\right)$,
where $D_{k}u=\partial u/\partial x_{k}$ and we write $v_{k}=D_{k}f$.
The vector of first-order weak derivatives can then be written as
$x\mapsto\nabla f=\left(D_{1}f,\dots,D_{d}f\right)\left(x\right)$
and we say that $f\in{\cal W}^{1,p}\left(\mathbb{R}^{d}\right)$ if
\[
\left\Vert f\right\Vert _{{\cal W}^{1,p}}=\left\Vert f\right\Vert _{p}+\sum_{k=1}^{d}\left\Vert D_{k}f\right\Vert _{p}<\infty\text{.}
\]
Similarly, for each $p\in\left[1,\infty\right)$ and $s\in\left(0,1\right)$,
we say that $f\in{\cal W}^{s,p}\left(\mathbb{R}^{d}\right)$ if 
\[
\left\Vert f\right\Vert _{{\cal W}^{s,p}}=\left\Vert f\right\Vert _{p}+\left\Vert f\right\Vert _{s,p}<\infty\text{,}
\]
where 
\[
\left\Vert f\right\Vert _{s,p}=\left[\int_{\mathbb{R}^{d}}\int_{\mathbb{R}^{d}}\frac{\left|f\left(x\right)-f\left(y\right)\right|^{p}}{\left\Vert x-y\right\Vert _{2}^{d+sp}}\mathrm{d}x\mathrm{d}y\right]^{1/p}\text{.}
\]
Note that the fractional Sobolev space ${\cal W}^{s,p}$ is equal
to the Besov space ${\cal B}_{p}^{s,p}$ (cf. \citealp{Leoni2023FractionalSobolev}).
When a statement below is written for $s\in\left(0,1\right]$, the
endpoint notation ${\cal W}^{s,p}$ at $s=1$ is understood as ${\cal W}^{1,p}$.

\paragraph{Classes of probability density functions}

Let 
$${\cal P}=\left\{ f:\mathbb{R}^{d}\to\mathbb{R}_{\ge0}\text{ | }f\text{ is Lebesgue-measurable and }\int_{\mathbb{R}^{d}}f\left(x\right)\mathrm{d}x=1\right\} $$
denote the set of all probability density functions (PDFs) on $\mathbb{R}^{d}$.
Our target density $f_{0}$ belongs to ${\cal P}$. Next, let ${\cal Q}\subset{\cal P}$
be an arbitrary class of PDFs. We then define the set of $m$-component
mixtures (convex combinations of $m$ vertices) generated by ${\cal Q}$
by 
\[
\text{co}_{m}\left({\cal Q}\right)=\left\{ \sum_{j=1}^{m}\pi_{j}q_{j}\text{ | }q_{j}\in{\cal Q}\text{ and }\pi_{j}\ge0\text{ for all }j\in\left[m\right]\text{, }\sum_{j=1}^{m}\pi_{j}=1\right\} \text{,}
\]
The class of all finite mixtures generated by ${\cal Q}$ is written
as 
\[
\text{co}\left({\cal Q}\right)=\bigcup_{m\in\mathbb{N}}\text{co}_{m}\left({\cal Q}\right)\text{,}
\]
whose ${\cal L}_{p}$ closure $\overline{\text{co}\left({\cal Q}\right)}^{p}$
consists of all ${\cal L}_{p}$ limits of sequences $\left(f_{n}\right)_{n}\subset\text{co}\left({\cal Q}\right)$.

For an arbitrary base PDF $\varphi\in{\cal P}$, we define its dilation
(in the language of \citealp{CheneyLight2000Approximation}) at scale
$\nu>0$ by $x\mapsto\varphi_{\nu}\left(x\right)=\nu^{d}\varphi\left(\nu x\right)$,
and subsequently we define the location--scale class of $\varphi$
by 
\[
{\cal P}_{\varphi}=\left\{ x\mapsto\varphi_{\nu}\left(x-\mu\right)\text{ | }\mu\in\mathbb{R}^{d}\text{ and }\nu>0\right\} \text{.}
\]
The location--scale mixtures of $\varphi$ are then $\text{co}\left({\cal P}_{\varphi}\right)$.
We will also find use for fixed-scale location classes of $\varphi$,
in which case we write ${\cal P}_{\varphi,\nu}=\left\{ x\mapsto\varphi_{\nu}\left(x-\mu\right)\text{ | }\mu\in\mathbb{R}^{d}\right\} $.
More generally, if $h_{\nu}\left(x\right)=\nu^{d}h\left(\nu x\right)$,
then 
\[
\left\Vert h_{\nu}\right\Vert _{p}=\nu^{d\left(1-1/p\right)}\left\Vert h\right\Vert _{p}\qquad\text{for }p\in\left[1,\infty\right)\text{,}
\]
and $\left\Vert h_{\nu}\right\Vert _{\infty}=\nu^{d}\left\Vert h\right\Vert _{\infty}$
whenever $h\in{\cal L}_{\infty}$.

\paragraph{Convolutions}

For two $\lambda$-measurable functions $f,g:\mathbb{R}^{d}\to\mathbb{R}$,
we write their convolution as 
\[
x\mapsto\left(f*g\right)\left(x\right)=\int_{\mathbb{R}^{d}}f\left(y\right)g\left(x-y\right)\mathrm{d}y=\int_{\mathbb{R}^{d}}f\left(x-y\right)g\left(y\right)\mathrm{d}y\text{,}
\]
where the convolution is understood as an a.e. defined function whenever
$\int_{\mathbb{R}^{d}}\left|f\left(y\right)g\left(x-y\right)\right|\mathrm{d}y<\infty$
for $\lambda$-a.e. $x$. If $\mu$ is a finite signed measure and
$g$ is $\mu$-measurable, we also write 
\[
x\mapsto\left(g*\mu\right)\left(x\right)=\int_{\mathbb{R}^{d}}g\left(x-y\right)\mu\left(\mathrm{d}y\right)
\]
whenever the integral exists. In particular, if $f\in{\cal L}_{1}$
and $g\in{\cal L}_{p}$ for some $p\ge1$, then the convolution exists and
$\left\Vert f*g\right\Vert _{p}\le\left\Vert f\right\Vert _{1}\left\Vert g\right\Vert _{p}$.
If, in addition, $g\in{\cal L}_{1}$, then
$\int_{\mathbb{R}^{d}}\left(f*g\right)\left(x\right)\mathrm{d}x=\int_{\mathbb{R}^{d}}f\left(x\right)\mathrm{d}x\int_{\mathbb{R}^{d}}g\left(y\right)\mathrm{d}y$.

\paragraph{Fourier analysis and finite signed measures}

We use the unitary Fourier transform 
\[
\widehat{g}\left(\xi\right)=\left(2\pi\right)^{-d/2}\int_{\mathbb{R}^{d}}e^{-i\left\langle x,\xi\right\rangle }g\left(x\right)\mathrm{d}x
\]
for $g\in{\cal L}_{1}\left(\mathbb{R}^{d}\right)$, and the same notation
for the Fourier transform of a finite signed measure $\mu$: 
\[
\widehat{\mu}\left(\xi\right)=\left(2\pi\right)^{-d/2}\int_{\mathbb{R}^{d}}e^{-i\left\langle x,\xi\right\rangle }\mu\left(\mathrm{d}x\right)\text{.}
\]
If $g\in{\cal L}_{1}$ and $\widehat{g}\in{\cal L}_{1}$, then Fourier
inversion yields 
\[
g\left(x\right)=\left(2\pi\right)^{-d/2}\int_{\mathbb{R}^{d}}e^{i\left\langle x,\xi\right\rangle }\widehat{g}\left(\xi\right)\mathrm{d}\xi
\]
for $\lambda$-a.e. $x$, and the right-hand side defines a continuous
version of $g$. The Fourier transform extends uniquely to an isometry
on ${\cal L}_{2}$, again denoted by $g\mapsto\widehat{g}$, and Plancherel's
theorem gives 
\[
\left\Vert g\right\Vert _{2}=\left\Vert \widehat{g}\right\Vert _{2}\qquad\text{for every }g\in{\cal L}_{2}\text{.}
\]
If $g,h\in{\cal L}_{1}$, then 
\[
\widehat{g*h}=\left(2\pi\right)^{d/2}\widehat{g}\,\widehat{h}\text{,}
\]
and, more generally, if $g\in{\cal L}_{1}$ and $\mu$ is a finite
signed measure, then 
\[
\widehat{g*\mu}=\left(2\pi\right)^{d/2}\widehat{g}\,\widehat{\mu}\text{.}
\]
For a finite signed measure $\mu$, we write $\left\Vert \mu\right\Vert _{\mathrm{TV}}$
for its total variation norm, and $\delta_{x}$ denotes the Dirac
mass at $x\in\mathbb{R}^{d}$. Finally, for the dilation $\varphi_{\nu}\left(x\right)=\nu^{d}\varphi\left(\nu x\right)$,
\[
\widehat{\varphi_{\nu}}\left(\xi\right)=\widehat{\varphi}\left(\xi/\nu\right)\text{.}
\]

\paragraph{Empirical processes}

Let $\left(\Omega,\mathfrak{F},\text{P}\right)$ be a probability
space with typical element $\omega$ and expectation operator $\text{E}$,
supporting the random vector $X:\left(\Omega,\mathfrak{F}\right)\to\left(\mathbb{R}^{d},\mathfrak{B}\left(\mathbb{R}^{d}\right)\right)$,
where $\mathfrak{B}\left(\mathbb{R}^{d}\right)$ is the Borel $\sigma$-algebra
of $\mathbb{R}^{d}$, and let $\left(X_{n}\right)_{n}$ be an independent
and identically distributed sample, where $X_{1}$ has the same distribution
as $X$. Let $\ell^{\infty}\left(\mathbb{R}^{d}\right)=\left\{ f:\mathbb{R}^{d}\to\mathbb{R}:\sup_{x\in\mathbb{R}^{d}}\left|f\left(x\right)\right|<\infty\right\} $
denote the space of bounded functions on $\mathbb{R}^{d}$, and for
$f\in\ell^{\infty}\left(\mathbb{R}^{d}\right)$, we define the operations
$f\mapsto P_{n}f=n^{-1}\sum_{i=1}^{n}f\left(X_{i}\right)$ and $f\mapsto Pf=\text{E}f\left(X\right)$.
For a set ${\cal F}\subset\ell^{\infty}\left(\mathbb{R}^{d}\right)$,
we say that it is a VC-subgraph class if the family of sets ${\cal C}=\left\{ \mathbb{G}_{f}\text{ | }f\in\mathcal{F}\right\} $
has finite Vapnik--Chervonenkis (VC) dimension, and we write this
number as $\operatorname{VC}\left({\cal C}\right)$, where $\mathbb{G}_{f}=\left\{ \left(x,t\right)\in\mathbb{R}^{d}\times\mathbb{R}\text{ | }t\le f\left(x\right)\right\} $
(cf. \citealp[Sec. 3.6]{GineNickl2016}) and with some abuse of nomenclature,
also say that $\operatorname{VC}\left({\cal C}\right)$ is the VC
dimension of ${\cal F}$: $\operatorname{VC}\left({\cal F}\right)$.
We say that ${\cal F}$ has an envelope $\bar{F}:\mathbb{R}^{d}\to\mathbb{R}$
if $\left|f\left(x\right)\right|\le\bar{F}\left(x\right)$ for every
$f\in{\cal F}$ and $x\in\mathbb{R}^{d}$, and we write its ${\cal L}_{p}\left(P\right)$
norm as $\left\Vert \bar{F}\right\Vert _{p}$. We shall also say that
${\cal F}$ is pointwise measurable if there exists a countable subclass
${\cal F}_{0}\subset{\cal F}$ such that for each $f\in{\cal F}$ there
is a sequence $\left(f_{n}\right)_{n}\subset{\cal F}_{0}$ with $f_{n}\left(x\right)\to f\left(x\right)$
for every $x\in\mathbb{R}^{d}$ (cf. \citealp[Ex.~2.3.4]{VanderVaartWellner2023}).
For a possibly nonmeasurable nonnegative random variable $Z$, $\mathrm{E}^{*}Z$
denotes its outer expectation.

\section{Approximation rates}

\label{sec:Approximation-rates-for}

The main tool of our analysis is the following fixed-scale quantization
lemma. The mixture model approximation theorem then follows by combining
this fixed-scale quantization bound with a smoothing-bias estimate
and optimizing in $\nu$. 
\begin{lem}
\label{lem:FixedScaleCompression}For $p\in\left(1,\infty\right)$,
let $\varphi,f_{0}\in{\cal L}_{p}\cap{\cal P}$. Then, for every $\nu>0$,
there exists a sequence $\left(f_{m,\nu}\right)_{m}$ with $f_{m,\nu}\in\mathrm{co}_{m}\left({\cal P}_{\varphi,\nu}\right)$,
for each $m\in\mathbb{N}$, such that 
\[
\left\Vert f_{m,\nu}-\varphi_{\nu}*f_{0}\right\Vert _{p}\le\begin{cases}
3\nu^{d/q}\left\Vert \varphi\right\Vert _{p}C_{p}m^{-1/q} & \text{if }1<p<2\text{,}\\
3\nu^{d/q}\left\Vert \varphi\right\Vert _{p}C_{p}m^{-1/2} & \text{if }2\le p<\infty\text{,}
\end{cases}
\]
where $1/p+1/q=1$ and $C_{p}<\infty$ is the constant from Lemma~\ref{lem:convex-approximation}. 
\end{lem}

\begin{proof}
The proof proceeds in two steps.

\paragraph{Step 1 (Membership in the fixed-scale convex closure)}

Fix $\nu>0$ and write 
\[
g=\varphi_{\nu}*f_{0}\text{.}
\]
We claim that $g\in\overline{\mathrm{co}\left({\cal P}_{\varphi,\nu}\right)}^{p}$.
Let $\epsilon>0$. By Lemma~\ref{lem:translation-continuity-Lp},
there exists $\delta>0$ such that 
\[
\left\Vert \varphi_{\nu}\left(\cdot-z\right)-\varphi_{\nu}\right\Vert _{p}<\frac{\epsilon}{2}
\]
whenever $\left\Vert z\right\Vert _{2}<\delta$. By translation-invariance
of the ${\cal L}_{p}$ norm, it follows that 
\[
\left\Vert \varphi_{\nu}\left(\cdot-s\right)-\varphi_{\nu}\left(\cdot-t\right)\right\Vert _{p}<\frac{\epsilon}{2}
\]
whenever $\left\Vert s-t\right\Vert _{2}<\delta$.

Next, choose $R>0$ such that, with 
\[
Q_{R}=\left[-R,R\right]^{d},
\]
one has 
\[
\int_{\mathbb{R}^{d}\setminus Q_{R}}f_{0}\left(s\right)\mathrm{d}s<\frac{\epsilon}{4\left\Vert \varphi_{\nu}\right\Vert _{p}}\text{.}
\]
Partition $Q_{R}$ into finitely many half-open cubes $A_{1},\dots,A_{N}$,
each having diameter strictly less than $\delta$, and choose $y_{j}\in A_{j}$
for each $j\in\left[N\right]$. Also set 
\[
A_{0}=\mathbb{R}^{d}\setminus Q_{R},
\]
and choose any $y_{0}\in A_{0}$. Define 
\[
\pi_{j}=\int_{A_{j}}f_{0}\left(s\right)\mathrm{d}s,\qquad j=0,1,\dots,N\text{.}
\]
Then $\pi_{j}\ge0$ for all $j$ and 
\[
\sum_{j=0}^{N}\pi_{j}=\int_{\mathbb{R}^{d}}f_{0}\left(s\right)\mathrm{d}s=1\text{.}
\]
Hence the function 
\[
g_{\epsilon}=\sum_{j=0}^{N}\pi_{j}\varphi_{\nu}\left(\cdot-y_{j}\right)
\]
belongs to $\mathrm{co}_{N+1}\left({\cal P}_{\varphi,\nu}\right)\subset\mathrm{co}\left({\cal P}_{\varphi,\nu}\right)$.

For $\lambda$-a.e. $x\in\mathbb{R}^{d}$, 
\begin{align*}
g\left(x\right)-g_{\epsilon}\left(x\right) & =\int_{\mathbb{R}^{d}}\varphi_{\nu}\left(x-s\right)f_{0}\left(s\right)\mathrm{d}s-\sum_{j=0}^{N}\pi_{j}\varphi_{\nu}\left(x-y_{j}\right)\\
 & =\sum_{j=0}^{N}\int_{A_{j}}\left[\varphi_{\nu}\left(x-s\right)-\varphi_{\nu}\left(x-y_{j}\right)\right]f_{0}\left(s\right)\mathrm{d}s\text{.}
\end{align*}
Therefore, by Minkowski's integral inequality (Lemma \ref{lem:minkowski-integral-inequality}), 
\begin{align*}
\left\Vert g-g_{\epsilon}\right\Vert _{p} & \le\sum_{j=0}^{N}\int_{A_{j}}f_{0}\left(s\right)\left\Vert \varphi_{\nu}\left(\cdot-s\right)-\varphi_{\nu}\left(\cdot-y_{j}\right)\right\Vert _{p}\mathrm{d}s\\
 & =\sum_{j=1}^{N}\int_{A_{j}}f_{0}\left(s\right)\left\Vert \varphi_{\nu}\left(\cdot-s\right)-\varphi_{\nu}\left(\cdot-y_{j}\right)\right\Vert _{p}\mathrm{d}s\\
 & \quad+\int_{A_{0}}f_{0}\left(s\right)\left\Vert \varphi_{\nu}\left(\cdot-s\right)-\varphi_{\nu}\left(\cdot-y_{0}\right)\right\Vert _{p}\mathrm{d}s\text{.}
\end{align*}
If $s\in A_{j}$ for some $j\in\left[N\right]$, then $\left\Vert s-y_{j}\right\Vert _{2}<\delta$,
so 
\[
\left\Vert \varphi_{\nu}\left(\cdot-s\right)-\varphi_{\nu}\left(\cdot-y_{j}\right)\right\Vert _{p}<\frac{\epsilon}{2}\text{.}
\]
Hence 
\[
\sum_{j=1}^{N}\int_{A_{j}}f_{0}\left(s\right)\left\Vert \varphi_{\nu}\left(\cdot-s\right)-\varphi_{\nu}\left(\cdot-y_{j}\right)\right\Vert _{p}\mathrm{d}s\le\frac{\epsilon}{2}\sum_{j=1}^{N}\pi_{j}\le\frac{\epsilon}{2}\text{.}
\]
For the tail cell $A_{0}$, the triangle inequality gives 
\[
\left\Vert \varphi_{\nu}\left(\cdot-s\right)-\varphi_{\nu}\left(\cdot-y_{0}\right)\right\Vert _{p}\le2\left\Vert \varphi_{\nu}\right\Vert _{p}\text{,}
\]
whence 
\[
\int_{A_{0}}f_{0}\left(s\right)\left\Vert \varphi_{\nu}\left(\cdot-s\right)-\varphi_{\nu}\left(\cdot-y_{0}\right)\right\Vert _{p}\mathrm{d}s\le2\left\Vert \varphi_{\nu}\right\Vert _{p}\int_{\mathbb{R}^{d}\setminus Q_{R}}f_{0}\left(s\right)\mathrm{d}s<\frac{\epsilon}{2}\text{.}
\]
Combining the last two displays yields 
\[
\left\Vert g-g_{\epsilon}\right\Vert _{p}<\epsilon\text{.}
\]
Since $\epsilon>0$ was arbitrary, it follows that 
\[
\varphi_{\nu}*f_{0}\in\overline{\mathrm{co}\left({\cal P}_{\varphi,\nu}\right)}^{p}\text{.}
\]

\paragraph{Step 2 (Fixed-scale quantitative quantization)}

Let $g\in\mathrm{co}\left({\cal P}_{\varphi,\nu}\right)$. Then there
exist $m\in\mathbb{N}$, $\mu_{1},\dots,\mu_{m}\in\mathbb{R}^{d}$,
and $\pi_{1},\dots,\pi_{m}\ge0$ with $\sum_{i=1}^{m}\pi_{i}=1$ such
that 
\[
g=\sum_{i=1}^{m}\pi_{i}\varphi_{\nu}\left(\cdot-\mu_{i}\right)\text{.}
\]
By Minkowski's inequality and translation-invariance of the ${\cal L}_{p}$
norm, 
\[
\left\Vert g\right\Vert _{p}\le\sum_{i=1}^{m}\pi_{i}\left\Vert \varphi_{\nu}\left(\cdot-\mu_{i}\right)\right\Vert _{p}=\left\Vert \varphi_{\nu}\right\Vert _{p}\text{.}
\]
Then, Young's inequality gives 
\[
\left\Vert \varphi_{\nu}*f_{0}\right\Vert _{p}\le\left\Vert \varphi_{\nu}\right\Vert _{p}\left\Vert f_{0}\right\Vert _{1}=\left\Vert \varphi_{\nu}\right\Vert _{p}\text{,}
\]
so 
\[
\left\Vert g-\varphi_{\nu}*f_{0}\right\Vert _{p}\le2\left\Vert \varphi_{\nu}\right\Vert _{p}\text{.}
\]
A change of variables shows that 
\[
\left\Vert \varphi_{\nu}\right\Vert _{p}=\left[\int_{\mathbb{R}^{d}}\left[\nu^{d}\varphi\left(\nu x\right)\right]^{p}\mathrm{d}x\right]^{1/p}=\nu^{d/q}\left\Vert \varphi\right\Vert _{p}\text{.}
\]
Since $\varphi_{\nu}*f_{0}\in\overline{\mathrm{co}\left({\cal P}_{\varphi,\nu}\right)}^{p}$,
Lemma~\ref{lem:convex-approximation} implies that, taking $B=2\nu^{d/q}\left\Vert \varphi\right\Vert _{p}$
and $\epsilon=\nu^{d/q}\left\Vert \varphi\right\Vert _{p}$, there
exists a sequence $\left(f_{m,\nu}\right)_{m}$ with $f_{m,\nu}\in\mathrm{co}_{m}\left({\cal P}_{\varphi,\nu}\right)$
such that 
\[
\left\Vert f_{m,\nu}-\varphi_{\nu}*f_{0}\right\Vert _{p}\le\begin{cases}
3\nu^{d/q}\left\Vert \varphi\right\Vert _{p}C_{p}m^{-1/q} & \text{if }1<p<2\text{,}\\
3\nu^{d/q}\left\Vert \varphi\right\Vert _{p}C_{p}m^{-1/2} & \text{if }2\le p<\infty\text{,}
\end{cases}
\]
as claimed. 
\end{proof}
\begin{rem}
\label{rem:Step-1-approx}The argument in Step~1 of the proof of
Lemma~\ref{lem:FixedScaleCompression} is a constructive version
of the separation argument used in \citet{nguyen2023approximation}.
Rather than invoking Hahn--Banach separation, it quantizes the target
density $f_{0}$ directly by a finite mixture of approximation kernels. 
\end{rem}

\subsection{Finite mixture approximation}

The next lemma isolates the smoothing-bias term that must be balanced
against the quantization error from Lemma~\ref{lem:FixedScaleCompression}. 
\begin{lem}
\label{lem:SmoothingBias}For $p\in\left(1,\infty\right)$, let $\varphi,f_{0}\in{\cal L}_{p}\cap{\cal P}$
and assume that for some finite $\alpha,K_{1},K_{2}>0$, 
\begin{equation}
\int_{\mathbb{R}^{d}}\left\Vert x\right\Vert _{2}^{\alpha}\varphi\left(x\right)\mathrm{d}x\le K_{1}\text{,}\label{eq:moment-assumption}
\end{equation}
and 
\begin{equation}
\left[\int_{\mathbb{R}^{d}}\left|f_{0}\left(x-y\right)-f_{0}\left(x\right)\right|^{p}\mathrm{d}x\right]^{1/p}\le K_{2}\left\Vert y\right\Vert _{2}^{\alpha}\qquad\text{for every }y\in\mathbb{R}^{d}\text{.}\label{eq:smoothness-assumption}
\end{equation}
Then, for every $\nu>0$, 
\[
\left\Vert \varphi_{\nu}*f_{0}-f_{0}\right\Vert _{p}\le K_{1}K_{2}\nu^{-\alpha}\text{.}
\]
\end{lem}

\begin{proof}
Fix $\nu>0$. The convolution identity below holds for $\lambda$-a.e.
$x$. Since $\varphi_{\nu}\in{\cal P}$, Jensen's inequality
followed by Minkowski's integral inequality (Lemma~\ref{lem:minkowski-integral-inequality})
gives 
\begin{align*}
\left\Vert \varphi_{\nu}*f_{0}-f_{0}\right\Vert _{p} & =\left[\int_{\mathbb{R}^{d}}\left|\int_{\mathbb{R}^{d}}\varphi_{\nu}\left(y\right)\left[f_{0}\left(x-y\right)-f_{0}\left(x\right)\right]\mathrm{d}y\right|^{p}\mathrm{d}x\right]^{1/p}\\
 & \le\left[\int_{\mathbb{R}^{d}}\left[\int_{\mathbb{R}^{d}}\varphi_{\nu}\left(y\right)\left|f_{0}\left(x-y\right)-f_{0}\left(x\right)\right|\mathrm{d}y\right]^{p}\mathrm{d}x\right]^{1/p}\\
 & \le\int_{\mathbb{R}^{d}}\varphi_{\nu}\left(y\right)\left[\int_{\mathbb{R}^{d}}\left|f_{0}\left(x-y\right)-f_{0}\left(x\right)\right|^{p}\mathrm{d}x\right]^{1/p}\mathrm{d}y\text{.}
\end{align*}
The translation modulus assumption (\ref{eq:smoothness-assumption})
and the change of variables $z=\nu y$ therefore yield 
\begin{align*}
\left\Vert \varphi_{\nu}*f_{0}-f_{0}\right\Vert _{p} & \le K_{2}\int_{\mathbb{R}^{d}}\left\Vert y\right\Vert _{2}^{\alpha}\varphi_{\nu}\left(y\right)\mathrm{d}y\\
 & =K_{2}\int_{\mathbb{R}^{d}}\left\Vert y\right\Vert _{2}^{\alpha}\nu^{d}\varphi\left(\nu y\right)\mathrm{d}y\\
 & =\frac{K_{2}}{\nu^{\alpha}}\int_{\mathbb{R}^{d}}\left\Vert z\right\Vert _{2}^{\alpha}\varphi\left(z\right)\mathrm{d}z\\
 & \le K_{1}K_{2}\nu^{-\alpha}\text{,}
\end{align*}
which proves the claim. 
\end{proof}
The following corollary is our main approximation theoretic outcome
and is obtained by balancing the fixed-scale quantization bound from
Lemma~\ref{lem:FixedScaleCompression} with the smoothing-bias bound
of Lemma~\ref{lem:SmoothingBias}. 
\begin{cor}
\label{cor:GeneralApproxRate}For $p\in\left(1,\infty\right)$, let
$\varphi,f_{0}\in{\cal L}_{p}\cap{\cal P}$ and assume that conditions
(\ref{eq:moment-assumption}) and (\ref{eq:smoothness-assumption})
hold for some finite $\alpha,K_{1},K_{2}>0$. Then, there exists a
sequence $\left(f_{m}\right)_{m}$ with $f_{m}\in\mathrm{co}_{m}\left({\cal P}_{\varphi}\right)$,
for each $m\in\mathbb{N}$, such that 
\[
\left\Vert f_{m}-f_{0}\right\Vert _{p}\le\begin{cases}
Km^{-\frac{\alpha}{\alpha q+d}} & \text{if }1<p<2\text{,}\\
Km^{-\frac{\alpha q}{2\left(\alpha q+d\right)}} & \text{if }2\le p<\infty\text{,}
\end{cases}
\]
where $1/p+1/q=1$ and $K$ is a function of $\alpha,K_{1},K_{2},p,d,\varphi$. 
\end{cor}

\begin{proof}
Fix $m\in\mathbb{N}$ and $\nu>0$. By Lemma~\ref{lem:FixedScaleCompression},
there exists $f_{m,\nu}\in\mathrm{co}_{m}\left({\cal P}_{\varphi,\nu}\right)\subset\mathrm{co}_{m}\left({\cal P}_{\varphi}\right)$
such that 
\[
\left\Vert f_{m,\nu}-f_{0}\right\Vert _{p}\le\left\Vert f_{m,\nu}-\varphi_{\nu}*f_{0}\right\Vert _{p}+\left\Vert \varphi_{\nu}*f_{0}-f_{0}\right\Vert _{p}\text{.}
\]
Combining Lemmas~\ref{lem:FixedScaleCompression} and~\ref{lem:SmoothingBias}
yields 
\[
\left\Vert f_{m,\nu}-f_{0}\right\Vert _{p}\le\begin{cases}
3\nu^{d/q}\left\Vert \varphi\right\Vert _{p}C_{p}m^{-1/q}+K_{1}K_{2}\nu^{-\alpha} & \text{if }1<p<2\text{,}\\
3\nu^{d/q}\left\Vert \varphi\right\Vert _{p}C_{p}m^{-1/2}+K_{1}K_{2}\nu^{-\alpha} & \text{if }2\le p<\infty\text{.}
\end{cases}
\]

If $1<p<2$, choose $\nu=\nu_{m}\asymp m^{1/\left(\alpha q+d\right)}$.
Then both terms on the right-hand side are of order $m^{-\alpha/\left(\alpha q+d\right)}$.
If $2\le p<\infty$, choose $\nu=\nu_{m}\asymp m^{q/\left(2\left(\alpha q+d\right)\right)}$.
Then both terms are of order $m^{-\alpha q/\left(2\left(\alpha q+d\right)\right)}$.
Absorbing the constants into $K$ proves the claim. 
\end{proof}

\subsection{Approximation in Sobolev spaces}

Condition (\ref{eq:smoothness-assumption}) is a smoothness assumption
enforced on $f_{0}$. For instance, if $f_{0}\in{\cal W}^{1,p}$,
it holds with $\alpha=1$ and 
\[
K_{2}=\left[\int_{\mathbb{R}^{d}}\left\Vert \nabla f_{0}\left(x\right)\right\Vert _{2}^{p}\mathrm{d}x\right]^{1/p}\text{,}
\]
and if $f_{0}\in{\cal W}^{s,p}$ with $0<s<1$, then the condition
holds with $\alpha=s$ and $K_{2}=C_{d,s,p}\left\Vert f_{0}\right\Vert _{s,p}$
for a finite constant $C_{d,s,p}>0$.

Indeed, inspecting the dependence of the constant $K$, only $\alpha$
and $K_{2}$ depend explicitly on $f_{0}$, while $K_{1}$ is determined
once $\varphi$ is chosen. Since $\varphi$ can be chosen in light
of the smoothness of $f_{0}$ so as to satisfy condition (\ref{eq:moment-assumption}),
and since a PDF $\varphi$ with compact support satisfies (\ref{eq:moment-assumption})
for every $\alpha>0$, Corollary~\ref{cor:GeneralApproxRate} yields
the following concrete result. 
\begin{cor}
\label{cor:ConcreteApproxRate}For $p\in\left(1,\infty\right)$ and
$s\in\left(0,1\right]$, if $\varphi\in{\cal L}_{p}\cap{\cal P}$,
$f_{0}\in{\cal W}^{s,p}\cap{\cal P}$, and 
\[
\int_{\mathbb{R}^{d}}\left\Vert x\right\Vert _{2}^{s}\varphi\left(x\right)\mathrm{d}x\le K_{1}\text{,}
\]
then there exists a sequence $\left(f_{m}\right)_{m}$ with $f_{m}\in\mathrm{co}_{m}\left({\cal P}_{\varphi}\right)$,
for each $m\in\mathbb{N}$, such that 
\[
\left\Vert f_{m}-f_{0}\right\Vert _{p}\le\begin{cases}
Km^{-\frac{s}{sq+d}} & \text{if }1<p<2\text{,}\\
Km^{-\frac{sq}{2\left(sq+d\right)}} & \text{if }2\le p<\infty\text{,}
\end{cases}
\]
where $K$ is a function of $s,K_{1},p,d,\varphi$ and depends on
$f_{0}$ only via $\left\Vert f_{0}\right\Vert _{s,p}$ when $s<1$,
and via $\left[\int_{\mathbb{R}^{d}}\left\Vert \nabla f_{0}\left(x\right)\right\Vert _{2}^{p}\mathrm{d}x\right]^{1/p}$
when $s=1$. 
\end{cor}

\begin{proof}
Apply Corollary~\ref{cor:GeneralApproxRate} with $\alpha=s$. When
$s\in\left(0,1\right)$, condition (\ref{eq:smoothness-assumption})
follows from Lemma~\ref{lem:smoothness-Wsp} with $K_{2}=C_{d,s,p}\left\Vert f_{0}\right\Vert _{s,p}$.
When $s=1$, condition (\ref{eq:smoothness-assumption}) follows from
Lemma~\ref{lem:smoothness-W1p} with 
\[
K_{2}=\left[\int_{\mathbb{R}^{d}}\left\Vert \nabla f_{0}\left(x\right)\right\Vert _{2}^{p}\mathrm{d}x\right]^{1/p}\text{.}
\]
The stated dependence of $K$ is therefore immediate from Corollary~\ref{cor:GeneralApproxRate}. 
\end{proof}
It is natural to ask whether the result above can be extended to targets
$f_{0}$ of greater Sobolev smoothness. The present argument, however,
is tailored to the regime $s\in\left(0,1\right]$, because it controls
the smoothing bias only through the first-order translation modulus
in condition (\ref{eq:smoothness-assumption}). Standard higher-order
kernel approximation results, such as the Calderón--Zygmund approach from Propositions 4.1.5 and 4.3.33
of \citet{GineNickl2016}, suggest a partial extension provided that
the base kernel has the required vanishing moments. In particular,
to obtain an approximation order strictly larger than one, one needs
additional moment-cancellation conditions of the form 
\[
\int_{\mathbb{R}^{d}}x^{\beta}\varphi\left(x\right)\mathrm{d}x=0
\]
for suitable multi-indices $\beta$. For $1<s\le2$, this requires
only vanishing first moments, which is compatible with centered symmetric
densities such as the Gaussian kernel. By contrast, once second-order
cancellation is required, the probability-mixture framework imposes
an obstruction: if $\varphi\in{\cal P}$ and 
\[
\int_{\mathbb{R}^{d}}\left\Vert x\right\Vert _{2}^{2}\varphi\left(x\right)\mathrm{d}x=0\text{,}
\]
then $\left\Vert x\right\Vert _{2}^{2}\varphi\left(x\right)=0$ for
Lebesgue-a.e. $x$, so $\varphi=0$ almost everywhere on $\mathbb{R}^{d}\setminus\left\{ 0\right\} $.
This is impossible for a probability density with respect to Lebesgue
measure. Consequently, the Calderón--Zygmund higher-order kernel
method cannot produce approximation orders beyond two within the present
class of nonnegative probability kernels.

\section{Least-squares estimation}

\label{sec:Finite-Gaussian-mixture-estimation}

In this section, we study empirical least-squares estimation over
finite mixture sieves. For a set ${\cal F}\subset{\cal L}_{2}\cap{\cal P}$
and $\epsilon\ge0$, we say that $\hat{f}_{n,\epsilon}\in{\cal F}$
is an $\epsilon$-minimizer of the least-squares criterion if 
\begin{equation}
\int_{\mathbb{R}^{d}}\hat{f}_{n,\epsilon}^{2}\left(x\right)\mathrm{d}x-2P_{n}\hat{f}_{n,\epsilon}\le\inf_{f\in{\cal F}}\left[\int_{\mathbb{R}^{d}}f^{2}\left(x\right)\mathrm{d}x-2P_{n}f\right]+\epsilon\text{.}\label{eq:sample-eps-minimizer-def}
\end{equation}
Observe that the population version of this criterion differs from
$\left\Vert f-f_{0}\right\Vert _{2}^{2}$ only by the constant $\int_{\mathbb{R}^{d}}f_{0}^{2}\left(x\right){\rm d}x$.
Whenever expectations are taken below, the chosen $\epsilon$-minimizer
is understood to be measurable with respect to the sample. Throughout
this section and the next, the observations $X_{1},\dots,X_{n}$ are
assumed to be i.i.d. with common density $f_{0}$ with respect to
Lebesgue measure. In particular, 
\[
Pf=\int_{\mathbb{R}^{d}}f\left(x\right)f_{0}\left(x\right)\mathrm{d}x
\]
whenever the integral is well-defined.

\subsection{A fast least-squares rate for kernels with exponentially decaying
Fourier transform}

We use the following assumption throughout this subsection.
\begin{assumption}[Exponential Fourier decay]
\label{ass:ExponentialFourierKernel}
The kernel $\varphi\in{\cal L}_{2}\cap{\cal L}_{\infty}\cap{\cal P}$
has an exponentially decaying Fourier transform: there exist constants
$c,C,\beta>0$ such that
\begin{equation}
\left|\widehat{\varphi}\left(\xi\right)\right|\le C\exp\left(-c\left\Vert \xi\right\Vert _{2}^{\beta}\right)\qquad\text{for all }\xi\in\mathbb{R}^{d}\text{.}\label{eq:ExpFourierCondition}
\end{equation}
Such kernels are often called supersmooth in the deconvolution literature;
see \citet[Sec.~2.4.3]{meister2009deconvolution}. 
\end{assumption}
\begin{rem}
\label{rem:ExponentialFourierRepresentative}Under Assumption~\ref{ass:ExponentialFourierKernel},
the exponential right-hand side is integrable on $\mathbb{R}^{d}$, so
$\widehat{\varphi}\in{\cal L}_{1}$. Therefore Fourier inversion yields
a continuous version of $\varphi$. Throughout this subsection, whenever
Assumption~\ref{ass:ExponentialFourierKernel} is imposed, we work
with that representative.
\end{rem}
The key point is a direct verification of the local continuity modulus
required by Theorem~\ref{thm:MassartNedelecLeastSquares} through
a Fourier decomposition that exploits Assumption~\ref{ass:ExponentialFourierKernel}.
Our approach is reminiscent of the Fourier-analytic upper-bound arguments
used by \citet{kim2022minimax} for Gaussian location mixtures, although
the present proof is tailored to empirical-risk minimization
rather than to a kernel estimator.

For $m\in\mathbb{N}$ and $\nu>0$, write 
\[
{\cal G}_{m,\nu}=\mathrm{co}_{m}\left({\cal P}_{\varphi,\nu}\right)\text{.}
\]
In the bounded-target setting considered below, we also write 
\[
M_{\nu}=\left\Vert f_{0}\right\Vert _{\infty}\vee\left\Vert \varphi_{\nu}\right\Vert _{\infty}=\left\Vert f_{0}\right\Vert _{\infty}\vee\nu^{d}\left\Vert \varphi\right\Vert _{\infty}\text{.}
\]

\begin{lem}
\label{lem:GaussianSieveSeparability}Under Assumption~\ref{ass:ExponentialFourierKernel},
for every $m\in\mathbb{N}$ and every $\nu>0$,
there exists a countable subclass ${\cal G}_{m,\nu}^{0}\subset{\cal G}_{m,\nu}$
such that every $f\in{\cal G}_{m,\nu}$ admits a sequence $\left(f_{k}\right)_{k}\subset{\cal G}_{m,\nu}^{0}$
with 
\[
f_{k}\left(x\right)\to f\left(x\right)\text{ for every }x\in\mathbb{R}^{d}\qquad\text{and}\qquad\left\Vert f_{k}-f\right\Vert _{2}\to0\text{.}
\]
\end{lem}

\begin{proof}
See the Appendix. 
\end{proof}
\begin{lem}
\label{lem:GaussianLocalizedOscillation}Under Assumption~\ref{ass:ExponentialFourierKernel},
let $f_{0}\in{\cal L}_{\infty}\cap{\cal P}$,
assume that $X_{1},\dots,X_{n}$ are i.i.d. with common density $f_{0}$,
and set
\[
M_{\nu}=\left\Vert f_{0}\right\Vert _{\infty}\vee\left\Vert \varphi_{\nu}\right\Vert _{\infty}\text{.}
\]
There exists a finite constant
$K_{d,\varphi}\ge1$, depending only on $d$, $\left\Vert \varphi\right\Vert _{\infty}$,
and the constants in Assumption~\ref{ass:ExponentialFourierKernel}, such that if
\[
A_{n}=K_{d,\varphi}\left(1+\left(\log n\right)^{d/\left(2\beta\right)}\right)\text{,}\qquad n\ge3\text{,}
\]
then, for every $m\in\mathbb{N}$, every $\nu>0$, every $u\in{\cal G}_{m,\nu}$,
and every $\sigma>0$ satisfying $\sigma\ge A_{n}/\sqrt{n}$, 
\[
\sqrt{n}\,\mathrm{E}^{*}\sup_{\substack{f\in{\cal G}_{m,\nu}\\
\left\Vert f-u\right\Vert _{2}\le2\sqrt{M_{\nu}}\,\sigma
}
}\left|\left(P_{n}-P\right)\left(f-u\right)\right|\le2M_{\nu}A_{n}\sigma\text{.}
\]
If, in addition, 
\[
\operatorname{supp}\left(\widehat{\varphi}\right)\subseteq B\left(0,R_{0}\right)
\]
for some $R_{0}>0$, then there exists a finite constant $K_{d,\varphi,R_{0}}^{\mathrm{bl}}\ge1$,
depending only on $d$, $R_{0}$, and $\left\Vert \varphi\right\Vert _{\infty}$,
such that the same conclusion holds with 
\[
A_{n}=K_{d,\varphi,R_{0}}^{\mathrm{bl}}\text{,}\qquad n\ge3\text{.}
\]
\end{lem}

\begin{proof}
See the Appendix. 
\end{proof}
\begin{thm}
\label{thm:FastGaussianLeastSquares}Let $s\in\left(0,1\right]$,
let $\varphi\in{\cal L}_{2}\cap{\cal L}_{\infty}\cap{\cal P}$ satisfy
Assumption~\ref{ass:ExponentialFourierKernel}, assume that 
\[
\int_{\mathbb{R}^{d}}\left\Vert x\right\Vert _{2}^{s}\varphi\left(x\right)\mathrm{d}x\le K_{1}\text{,}
\]
and let $f_{0}\in{\cal W}^{s,2}\cap{\cal L}_{\infty}\cap{\cal P}$.
For each $n\ge3$, let $m_{n}\in\mathbb{N}$, let $\nu_{n}>0$, let
$\epsilon_{n}\ge0$, and let $\hat{f}_{n}\in{\cal G}_{m_{n},\nu_{n}}$
satisfy (\ref{eq:sample-eps-minimizer-def}) with ${\cal F}={\cal G}_{m_{n},\nu_{n}}$.
Then 
\[
\mathrm{E}\left\Vert \hat{f}_{n}-f_{0}\right\Vert _{2}^{2}\lesssim\frac{\nu_{n}^{d}}{m_{n}}+\nu_{n}^{-2s}+M_{\nu_{n}}\frac{\left(1+\log n\right)^{d/\beta}}{n}+\epsilon_{n}\text{.}
\]
In particular, fix $M,N>0$, and define 
\[
m_{n}=\left\lceil M\frac{n}{\left(1+\log n\right)^{d/\beta}}\right\rceil \text{ \qquad and \qquad}\nu_{n}=N\left(\frac{n}{\left(1+\log n\right)^{d/\beta}}\right)^{1/\left(2s+d\right)}\text{,}
\]
for all $n\ge3$. If 
\[
\epsilon_{n}\lesssim\left(\frac{\left(1+\log n\right)^{d/\beta}}{n}\right)^{\frac{2s}{2s+d}}\text{,}
\]
then 
\[
\mathrm{E}\left\Vert \hat{f}_{n}-f_{0}\right\Vert _{2}^{2}\lesssim\left(\frac{\left(1+\log n\right)^{d/\beta}}{n}\right)^{\frac{2s}{2s+d}}\text{.}
\]

If, in addition, 
\[
\operatorname{supp}\left(\widehat{\varphi}\right)\subseteq B\left(0,R_{0}\right)
\]
for some $R_{0}>0$, then one has the bandlimited bound 
\[
\mathrm{E}\left\Vert \hat{f}_{n}-f_{0}\right\Vert _{2}^{2}\lesssim\frac{\nu_{n}^{d}}{m_{n}}+\nu_{n}^{-2s}+M_{\nu_{n}}\frac{1}{n}+\epsilon_{n}\text{.}
\]
In particular, fix $M,N>0$, and if one chooses 
\[
m_{n}=\left\lceil Mn\right\rceil \text{ \qquad and \qquad}\nu_{n}=Nn^{1/\left(2s+d\right)}
\]
for all $n\ge3$, then the condition 
\[
\epsilon_{n}\lesssim n^{-\frac{2s}{2s+d}}
\]
implies the bandlimited rate 
\[
\mathrm{E}\left\Vert \hat{f}_{n}-f_{0}\right\Vert _{2}^{2}\lesssim n^{-\frac{2s}{2s+d}}\text{.}
\]
\end{thm}

\begin{proof}
Fix $n\ge3$ and write 
\[
{\cal F}_{n}={\cal G}_{m_{n},\nu_{n}}\text{.}
\]
By Lemmas~\ref{lem:GaussianSieveSeparability} and~\ref{lem:GaussianLocalizedOscillation},
the hypotheses of Theorem~\ref{thm:MassartNedelecLeastSquares} hold
on ${\cal F}_{n}$ with $M=M_{\nu_{n}}$ and 
\[
A_{n}=K_{d,\varphi}\left(1+\left(\log n\right)^{d/\left(2\beta\right)}\right)\text{.}
\]
Consequently, 
\[
\mathrm{E}\left\Vert \hat{f}_{n}-f_{0}\right\Vert _{2}^{2}\le2\inf_{f\in{\cal F}_{n}}\left\Vert f-f_{0}\right\Vert _{2}^{2}+2\epsilon_{n}+C\,M_{\nu_{n}}\frac{A_{n}^{2}}{n}
\]
for some finite universal constant $C>0$.

Next, Lemma~\ref{lem:FixedScaleCompression} with $p=2$ yields an
element $g_{n}\in{\cal F}_{n}$ such that 
\[
\left\Vert g_{n}-\varphi_{\nu_{n}}*f_{0}\right\Vert _{2}^{2}\lesssim\frac{\nu_{n}^{d}}{m_{n}}\text{.}
\]
Since $f_{0}\in{\cal W}^{s,2}$, condition (\ref{eq:smoothness-assumption})
holds with $\alpha=s$ and some finite constant $K_{2}>0$; indeed,
\[
K_{2}=\begin{cases}
C_{d,s,2}\left\Vert f_{0}\right\Vert _{s,2} & \text{if }0<s<1\text{,}\\
\left[\int_{\mathbb{R}^{d}}\left\Vert \nabla f_{0}\left(x\right)\right\Vert _{2}^{2}\mathrm{d}x\right]^{1/2} & \text{if }s=1\text{,}
\end{cases}
\]
is admissible by Lemma~\ref{lem:smoothness-Wsp} when $0<s<1$ and
by Lemma~\ref{lem:smoothness-W1p} when $s=1$. Since the present
theorem assumes condition (\ref{eq:moment-assumption}) for the same
order $s$, Lemma~\ref{lem:SmoothingBias} yields 
\[
\left\Vert \varphi_{\nu_{n}}*f_{0}-f_{0}\right\Vert _{2}^{2}\lesssim\nu_{n}^{-2s}\text{.}
\]
Therefore, 
\begin{align*}
\inf_{f\in{\cal F}_{n}}\left\Vert f-f_{0}\right\Vert _{2}^{2} & \le\left\Vert g_{n}-f_{0}\right\Vert _{2}^{2}\\
 & \le2\left\Vert g_{n}-\varphi_{\nu_{n}}*f_{0}\right\Vert _{2}^{2}+2\left\Vert \varphi_{\nu_{n}}*f_{0}-f_{0}\right\Vert _{2}^{2}\\
 & \lesssim\frac{\nu_{n}^{d}}{m_{n}}+\nu_{n}^{-2s}\text{.}
\end{align*}
Since $A_{n}^{2}\lesssim\left(1+\log n\right)^{d/\beta}$, this proves
the first displayed bound.

For the tuned choice under Assumption~\ref{ass:ExponentialFourierKernel}, write
\[
r_{n}=\frac{\left(1+\log n\right)^{d/\beta}}{n}\text{,}\qquad n\ge3\text{.}
\]
Then 
\[
M_{\nu_{n}}\le\left\Vert f_{0}\right\Vert _{\infty}+\left\Vert \varphi_{\nu_{n}}\right\Vert _{\infty}\lesssim1+\nu_{n}^{d}\text{,}
\]
so 
\[
M_{\nu_{n}}r_{n}\lesssim r_{n}+\nu_{n}^{d}r_{n}\text{.}
\]
Moreover, 
\[
\frac{\nu_{n}^{d}}{m_{n}}\asymp\nu_{n}^{d}r_{n}\asymp\nu_{n}^{-2s}\asymp r_{n}^{\frac{2s}{2s+d}}\text{.}
\]
Since $r_{n}\to0$ and $0<2s/\left(2s+d\right)<1$, we also have 
\[
r_{n}\lesssim r_{n}^{\frac{2s}{2s+d}}
\]
after enlarging the implicit constant, if necessary, to absorb the
finitely many values of $n\ge3$ for which $r_{n}>1$. Combining these
bounds with the first part of the theorem proves the claimed logarithmic
rate.

Now assume, in addition, that $\operatorname{supp}\left(\widehat{\varphi}\right)\subseteq B\left(0,R_{0}\right)$
for some $R_{0}>0$. By Lemmas~\ref{lem:GaussianSieveSeparability} and~\ref{lem:GaussianLocalizedOscillation}, the hypotheses of
Theorem~\ref{thm:MassartNedelecLeastSquares} also hold on ${\cal F}_{n}$
with $M=M_{\nu_{n}}$ and 
\[
A_{n}=K_{d,\varphi,R_{0}}^{\mathrm{bl}}\text{.}
\]
Therefore 
\[
\mathrm{E}\left\Vert \hat{f}_{n}-f_{0}\right\Vert _{2}^{2}\le2\inf_{f\in{\cal F}_{n}}\left\Vert f-f_{0}\right\Vert _{2}^{2}+2\epsilon_{n}+C\,M_{\nu_{n}}\frac{1}{n}\text{,}
\]
which, together with the approximation bound above, proves the bandlimited
result.

Now fix $M,N>0$ and take 
\[
m_{n}=\left\lceil Mn\right\rceil \text{ \qquad and \qquad}\nu_{n}=Nn^{1/\left(2s+d\right)}
\]
for all $n\ge3$. Since 
\[
M_{\nu_{n}}\le\left\Vert f_{0}\right\Vert _{\infty}+\left\Vert \varphi_{\nu_{n}}\right\Vert _{\infty}\lesssim1+\nu_{n}^{d}\text{,}
\]
the bandlimited bound yields 
\[
M_{\nu_{n}}\frac{1}{n}\lesssim\frac{1}{n}+\frac{\nu_{n}^{d}}{n}\text{.}
\]
Moreover, 
\[
\frac{\nu_{n}^{d}}{m_{n}}\asymp\frac{\nu_{n}^{d}}{n}\asymp\nu_{n}^{-2s}\asymp n^{-\frac{2s}{2s+d}}\text{.}
\]
Since $0<2s/\left(2s+d\right)<1$, we also have $n^{-1}\lesssim n^{-\frac{2s}{2s+d}}$.
Combining these estimates with the bandlimited display proves the
final claim. 
\end{proof}
\begin{rem}
\label{rem:BandlimitedMomentExamples}The bandlimited part of Theorem~\ref{thm:FastGaussianLeastSquares}
is compatible with the moment assumption from the theorem. Indeed,
let $\psi\in{\cal C}_{c}^{\infty}\left(\mathbb{R}^{d}\right)$ be
nonzero, and let $g$ denote the inverse Fourier transform of $\psi$,
so 
\[
\widehat{g}=\psi\text{.}
\]
Set 
\[
\varphi\left(x\right)=\frac{\left|g\left(x\right)\right|^{2}}{\int_{\mathbb{R}^{d}}\left|g\left(u\right)\right|^{2}\mathrm{d}u}\text{,}\qquad x\in\mathbb{R}^{d}\text{.}
\]
Then $\varphi\in{\cal P}\cap{\cal L}_{2}\cap{\cal L}_{\infty}$ and
$\varphi\ge0$. Since $\psi\in{\cal C}_{c}^{\infty}\left(\mathbb{R}^{d}\right)$,
its inverse Fourier transform $g$ belongs to ${\cal S}\left(\mathbb{R}^{d}\right)$.
Hence $\varphi\in{\cal S}\left(\mathbb{R}^{d}\right)$ as well (cf.
\citealp[Thm. 7.1.5]{Hormander2003}), so 
\[
\int_{\mathbb{R}^{d}}\left\Vert x\right\Vert _{2}^{t}\varphi\left(x\right)\mathrm{d}x<\infty\qquad\text{for every }t>0\text{.}
\]
Moreover, if $\widetilde{\psi}\left(\xi\right)=\overline{\psi\left(-\xi\right)}$,
then 
\[
\widehat{\varphi}=\left(2\pi\right)^{-d/2}\left\Vert g\right\Vert _{2}^{-2}\,\psi*\widetilde{\psi}\text{,}
\]
so $\widehat{\varphi}$ is compactly supported. Thus bandlimited kernels
satisfying the hypotheses of Theorem~\ref{thm:FastGaussianLeastSquares}
exist in abundance.

For a simple one-dimensional example, one may take 
\[
\varphi\left(x\right)=\frac{3}{4\pi}\left(\frac{\sin\left(x/2\right)}{x/2}\right)^{4}\text{,}\qquad x\in\mathbb{R}\text{.}
\]
Then $\varphi\in{\cal P}$ and 
\[
\widehat{\varphi}\left(\xi\right)=\left(2\pi\right)^{-1/2}\begin{cases}
1-\frac{3}{2}\left|\xi\right|^{2}+\frac{3}{4}\left|\xi\right|^{3} & \text{if }\left|\xi\right|\le1\text{,}\\
\frac{1}{4}\left(2-\left|\xi\right|\right)^{3} & \text{if }1\le\left|\xi\right|\le2\text{,}\\
0 & \text{if }\left|\xi\right|>2\text{,}
\end{cases}
\]
so $\widehat{\varphi}$ is compactly supported. Since $\varphi\left(x\right)\lesssim\left(1+\left|x\right|\right)^{-4}$,
we have 
\[
\int_{\mathbb{R}}\left|x\right|^{s}\varphi\left(x\right)\mathrm{d}x<\infty\qquad\text{for every }s\in\left(0,3\right)\text{,}
\]
and in particular for every $s\in\left(0,1\right]$. Tensor products
provide multivariate examples. 
\end{rem}

\begin{rem}
\label{rem:GaussianKernelSpecial}Choosing $\varphi$ to be the standard
Gaussian kernel 
\[
\varphi\left(x\right)=\left(2\pi\right)^{-d/2}\exp\left(-\left\Vert x\right\Vert _{2}^{2}/2\right)\text{,}\qquad x\in\mathbb{R}^{d}\text{,}
\]
we have 
\[
\widehat{\varphi}\left(\xi\right)=\left(2\pi\right)^{-d/2}\exp\left(-\left\Vert \xi\right\Vert _{2}^{2}/2\right)\text{,}
\]
so Assumption~\ref{ass:ExponentialFourierKernel} holds with $\beta=2$, $c=1/2$,
and $C=\left(2\pi\right)^{-d/2}$. Therefore Theorem~\ref{thm:FastGaussianLeastSquares}
specializes to 
\[
\mathrm{E}\left\Vert \hat{f}_{n}-f_{0}\right\Vert _{2}^{2}\lesssim\frac{\nu_{n}^{d}}{m_{n}}+\nu_{n}^{-2s}+M_{\nu_{n}}\frac{\left(1+\log n\right)^{d/2}}{n}+\epsilon_{n}\text{,}
\]
and, after the tuned choice of $m_{n}$ and $\nu_{n}$, to 
\[
\mathrm{E}\left\Vert \hat{f}_{n}-f_{0}\right\Vert _{2}^{2}\lesssim\left(\frac{\left(1+\log n\right)^{d/2}}{n}\right)^{\frac{2s}{2s+d}}\text{.}
\]
More generally, one-dimensional symmetric stable densities provide
examples with $\beta\in\left(0,2\right]$ (cf. \citealp{nolan2020univariate}).
Under the present unitary Fourier convention, their Fourier transforms
take the form 
\[
\left(2\pi\right)^{-1/2}\exp\left(-c\left|\xi\right|^{\beta}\right)
\]
for some $c>0$ and $\xi\in\mathbb{R}$. Tensor products of such one-dimensional
kernels then furnish multivariate examples satisfying Assumption~\ref{ass:ExponentialFourierKernel}.
For the full conclusion of Theorem~\ref{thm:FastGaussianLeastSquares}
at smoothness $s$, the order-$s$ moment condition must also hold;
for stable kernels this requires $s<\beta$. 
\end{rem}

\begin{rem}
\label{rem:FastGaussianMinimax}The classical minimax benchmark for
squared ${\cal L}_{2}$ risk over $d$-dimensional Sobolev balls of
smoothness $s$ is $n^{-2s/\left(2s+d\right)}$; see, for example,
\citet[Cor. 11.11]{Klemela2009Smoothing} and \citet[Ex. 6.3.2]{GineNickl2016}.
Therefore, under the order-$s$ moment hypothesis on $\varphi$ and
for bounded targets $f_{0}\in{\cal W}^{s,2}\cap{\cal L}_{\infty}\cap{\cal P}$,
the rate in Theorem~\ref{thm:FastGaussianLeastSquares} is the classical
Sobolev minimax benchmark multiplied by the logarithmic factor 
\[
\left(1+\log n\right)^{\frac{2sd}{\beta\left(2s+d\right)}}\text{.}
\]
If, in addition, $\widehat{\varphi}$ is compactly supported, then
the bandlimited part of Theorem~\ref{thm:FastGaussianLeastSquares}
recovers the classical Sobolev rate $n^{-2s/\left(2s+d\right)}$.
In particular, the Gaussian choice $\beta=2$ yields the factor $\left(1+\log n\right)^{sd/\left(2s+d\right)}$
in the supersmooth case, whereas the bandlimited case carries no logarithmic
loss. 
\end{rem}

\subsection{A slower rate under weaker assumptions}

\label{subsec:A-slower-rate}

The fast least-squares theorem above is specific to kernels whose
Fourier transforms decay exponentially. Under much weaker assumptions,
one can still obtain a slower least-squares benchmark of the same
order as the deterministic approximation bound from Section~\ref{sec:Approximation-rates-for}.
The only additional input is the standard least-squares basic inequality,
recorded as Lemma~\ref{lem:LeastSquaresBasicIneq} in the Appendix. 
\begin{thm}
\label{thm:SlowLeastSquaresBenchmark}Let $s\in\left(0,1\right]$,
let $\varphi\in\ell^{\infty}\left(\mathbb{R}^{d}\right)\cap{\cal L}_{2}\cap{\cal P}$,
let $f_{0}\in{\cal W}^{s,2}\cap{\cal P}$, assume that each fixed-scale
class ${\cal P}_{\varphi,\nu}$ is pointwise measurable, and suppose that $\operatorname{VC}\left({\cal P}_{\varphi}\right)<\infty$
and 
\[
\int_{\mathbb{R}^{d}}\left\Vert x\right\Vert _{2}^{s}\varphi\left(x\right)\mathrm{d}x\le K_{1}\text{.}
\]
Then there exist constants $M,N>0$ and a sequence $\left(\epsilon_{n}\right)_{n}$
with $\epsilon_{n}\lesssim n^{-\frac{s}{2s+d}}$ such that, if 
\[
m_{n}=\left\lceil Mn^{1/2}\right\rceil \text{ \qquad and \qquad}\nu_{n}=Nn^{1/\left(2\left(2s+d\right)\right)}\text{,}
\]
and if, for each $n\in\mathbb{N}$, $\hat{f}_{n}\in\mathrm{co}_{m_{n}}\left({\cal P}_{\varphi,\nu_{n}}\right)\subset\mathrm{co}_{m_{n}}\left({\cal P}_{\varphi}\right)$
is measurable and satisfies (\ref{eq:sample-eps-minimizer-def}) with ${\cal F}=\mathrm{co}_{m_{n}}\left({\cal P}_{\varphi,\nu_{n}}\right)$
and $\epsilon=\epsilon_{n}$, then 
\[
\mathrm{E}\left\Vert \hat{f}_{n}-f_{0}\right\Vert _{2}^{2}\lesssim n^{-\frac{s}{2s+d}}\text{.}
\]
\end{thm}

\begin{proof}
See the Appendix. 
\end{proof}
\begin{rem}
\label{rem:The-VC-assumption}The VC assumption in Theorem~\ref{thm:SlowLeastSquaresBenchmark}
is needed only for the empirical-process bound and is separate from
the localized Fourier argument above. The pointwise measurability assumption
is mild: if $\varphi$ is continuous, then, for each fixed $\nu>0$,
${\cal P}_{\varphi,\nu}$ is pointwise measurable, since the countable subclass
indexed by $\mu\in\mathbb{Q}^{d}$ is pointwise dense.

A convenient sufficient condition for the VC-subgraph property is
that the map 
\[
\left(x,t,\mu,\nu\right)\mapsto\mathbf{1}\left\{ t\le\nu^{d}\varphi\left(\nu\left(x-\mu\right)\right)\right\} 
\]
be definable in an $o$-minimal expansion of the ordered field of
real numbers. In that case, the subgraph family of ${\cal P}_{\varphi}$
is a definable family with parameter vector $\left(\mu,\nu\right)\in\mathbb{R}^{d}\times\left(0,\infty\right)$,
and hence $\operatorname{VC}\left({\cal P}_{\varphi}\right)<\infty$
by \citet[Thm. 10.11]{vidyasagar2003learning}; see also \citet[Ch. 5]{van1998tame}.
In particular, the Gaussian kernel is definable in $\mathbb{R}_{\mathrm{an},\exp}$,
while the usual compactly supported polynomial kernels used in multivariate
kernel density estimation, such as the triangular, Epanechnikov, biweight,
and triweight kernels, are all definable in $\mathbb{R}_{\mathrm{alg}}$;
see \citet[Ch. 6]{scott2015multivariate}. Therefore these kernels
satisfy the qualitative VC hypothesis required for Theorem~\ref{thm:SlowLeastSquaresBenchmark}. 
\end{rem}

\section{Location-mixture targets}

\label{sec:LocationMixtureTargets}

We now turn to fixed-scale location-mixture model targets. For a kernel
$\varphi\in{\cal P}$ represented by a Borel density and a scale $\nu>0$,
write 
\[
{\cal M}_{\varphi,\nu}=\left\{ x\mapsto\int_{\mathbb{R}^{d}}\varphi_{\nu}\left(x-\mu\right)G\left(\mathrm{d}\mu\right)\ |\ G\text{ is a probability measure on }\mathbb{R}^{d}\right\} \text{,}
\]
and, for any class ${\cal F}\subset{\cal L}_{2}\cap{\cal P}$, let
\[
R_{n}\left({\cal F}\right)=\inf_{\tilde{f}_{n}}\sup_{f_{0}\in{\cal F}}\mathrm{E}_{f_{0}}\left\Vert \tilde{f}_{n}-f_{0}\right\Vert _{2}^{2}
\]
denote the squared ${\cal L}_{2}$ minimax risk over ${\cal F}$,
where $\tilde{f}_{n}$ ranges over all measurable estimators based
on $n$ i.i.d.\ observations, with the convention that the loss is infinite whenever $\tilde{f}_{n}$ is not ${\cal L}_{2}$-valued. Throughout this section, let 
\[
\gamma\left(x\right)=\left(2\pi\right)^{-d/2}\exp\left(-\left\Vert x\right\Vert _{2}^{2}/2\right)\text{,}\qquad x\in\mathbb{R}^{d}\text{,}
\]
denote the standard Gaussian density. The sieve classes ${\cal G}_{m,\nu}=\mathrm{co}_{m}\left({\cal P}_{\varphi,\nu}\right)$
studied above are natural finite-dimensional approximations to ${\cal M}_{\varphi,\nu}$.
Consequently, the fixed-scale problem isolates the statistical effect
of kernel Fourier decay without the additional smoothing bias $\nu^{-2s}$
that appears for general Sobolev targets. The upper bound below applies
throughout the Fourier-tail class. The subsequent lower bounds show
that, for the Gaussian kernel and for the tensor-product odd-degree
Student-$t$ family, this fixed-scale rate is minimax-sharp.

\subsection{A fast least-squares estimator for location-mixture targets}

The fixed-scale theorem below should be compared with the greedy least-squares
procedure of \citet{klemela2007density}, built from a fixed dictionary,
which, under a finite entropy-integral condition, can be proved to
obtain a global squared ${\cal L}_{2}$ risk bound of order $n^{-1/2}$
once the target density belongs to the ${\cal L}_{2}$-closure of
the corresponding convex hull. Here, we retain this global dictionary-driven
viewpoint. When the dictionary is ${\cal P}_{\varphi,\nu}$ and $\varphi$
satisfies Assumption~\ref{ass:ExponentialFourierKernel}, the localized continuity
modulus from Lemma~\ref{lem:GaussianLocalizedOscillation} yields
the following fast least-squares analogue. 
\begin{thm}
\label{thm:FastFixedScaleLeastSquares}Fix $\nu>0$, let $\varphi\in{\cal L}_{2}\cap{\cal L}_{\infty}\cap{\cal P}$
satisfy Assumption~\ref{ass:ExponentialFourierKernel}, and suppose that 
\[
f_{0}\in\overline{\mathrm{co}\left({\cal P}_{\varphi,\nu}\right)}^{2}\cap{\cal P}\text{.}
\]
For each $n\ge3$, let $m_{n}\in\mathbb{N}$, let $\epsilon_{n}\ge0$,
and let $\hat{f}_{n}\in{\cal G}_{m_{n},\nu}$ be a measurable selection satisfying
(\ref{eq:sample-eps-minimizer-def}) with ${\cal F}={\cal G}_{m_{n},\nu}$.
Then 
\[
\mathrm{E}\left\Vert \hat{f}_{n}-f_{0}\right\Vert _{2}^{2}\lesssim\frac{\nu^{d}}{m_{n}}+\nu^{d}\frac{\left(1+\log n\right)^{d/\beta}}{n}+\epsilon_{n}\text{.}
\]
The implicit constants in this display do not depend on $f_{0}$. Consequently,
for each $n$, fix any measurable $\epsilon_{n}$-minimizer selection $\hat{f}_{n}$
over ${\cal G}_{m_{n},\nu}$. If $M>0$ is fixed and 
\[
m_{n}=\left\lceil M\frac{n}{\left(1+\log n\right)^{d/\beta}}\right\rceil 
\]
for all $n\ge3$, and if 
\[
\epsilon_{n}\lesssim\frac{\left(1+\log n\right)^{d/\beta}}{n}\text{,}
\]
then 
\[
\sup_{f_{0}\in{\cal M}_{\varphi,\nu}}\mathrm{E}\left\Vert \hat{f}_{n}-f_{0}\right\Vert _{2}^{2}\lesssim\nu^{d}\frac{\left(1+\log n\right)^{d/\beta}}{n}\text{.}
\]

If, in addition, 
\[
\operatorname{supp}\left(\widehat{\varphi}\right)\subseteq B\left(0,R_{0}\right)
\]
for some $R_{0}>0$, then one has the bandlimited bound 
\[
\mathrm{E}\left\Vert \hat{f}_{n}-f_{0}\right\Vert _{2}^{2}\lesssim\frac{\nu^{d}}{m_{n}}+\nu^{d}\frac{1}{n}+\epsilon_{n}\text{.}
\]
Consequently, with the same convention of a fixed measurable selection, if $M>0$ is fixed and 
\[
m_{n}=\left\lceil Mn\right\rceil 
\]
for all $n\ge3$, and if 
\[
\epsilon_{n}\lesssim\frac{1}{n}\text{,}
\]
then 
\[
\sup_{f_{0}\in{\cal M}_{\varphi,\nu}}\mathrm{E}\left\Vert \hat{f}_{n}-f_{0}\right\Vert _{2}^{2}\lesssim\frac{\nu^{d}}{n}\text{.}
\]
The implicit constants in the two fixed-scale ``in particular'' bounds
may depend on the fixed scale $\nu$, but not on $n$ or on the particular
$f_{0}\in{\cal M}_{\varphi,\nu}$.
\end{thm}

\begin{proof}
See the Appendix. 
\end{proof}
\begin{rem}
\label{rem:GaussianKimMinimax}For $\nu>0$, write 
\[
{\cal M}_{\gamma,\nu}=\left\{ x\mapsto\int_{\mathbb{R}^{d}}\gamma_{\nu}\left(x-\mu\right)G\left(\mathrm{d}\mu\right)\text{ | }G\text{ is a probability measure on }\mathbb{R}^{d}\right\} \text{.}
\]
Choosing $\varphi=\gamma$, $\nu=1$,
\[
m_{n}=\left\lceil M\frac{n}{\left(1+\log n\right)^{d/2}}\right\rceil
\]
for a fixed $M>0$, and $\epsilon_{n}\lesssim\left(1+\log n\right)^{d/2}/n$, we have $\beta=2$, so Theorem~\ref{thm:FastFixedScaleLeastSquares} yields 
\[
\sup_{f_{0}\in{\cal M}_{\gamma,1}}\mathrm{E}\left\Vert \hat{f}_{n}-f_{0}\right\Vert _{2}^{2}\lesssim\frac{\left(1+\log n\right)^{d/2}}{n}\text{.}
\]
Thus Theorem~\ref{thm:FastFixedScaleLeastSquares} matches the minimax order $n^{-1}\left(\log n\right)^{d/2}$
proved by \citet{kim2022minimax}. In this sense, Theorem~\ref{thm:FastFixedScaleLeastSquares}
furnishes a fast localized least-squares analogue of \citet[Thm. 2]{klemela2007density}.
The next subsection shows that the same lower-bound order persists
on many Gaussian convolution submodels of the Gaussian location-mixture
class. 
\end{rem}

\subsection{A lower bound on a Gaussian convolution subclass}

Theorem~\ref{thm:FastFixedScaleLeastSquares} gives an upper bound
throughout the Fourier-tail class. As discussed in Section~\ref{sec:Discussion},
no matching lower bound can hold uniformly at that level of generality.
The next result nevertheless identifies a Gaussian convolution family,
including the Gaussian kernel itself, on which the lower bound of
\citet[Thm.~2.1]{kim2022minimax} transfers by a bounded deconvolution
argument. As we report below, the corresponding model ${\cal M}_{\varphi,\nu}$
is always contained in ${\cal M}_{\gamma,\nu}$, and it is a strict
submodel whenever the auxiliary factor $\eta$ is not a point mass. 
\begin{thm}
\label{thm:GaussianCoreLowerBound}Fix $\nu>0$, and suppose that
\[
\varphi=\gamma*\eta
\]
for some Borel probability measure $\eta$ on $\mathbb{R}^{d}$ satisfying
\[
a_{\eta}=\inf_{\xi\in\mathbb{R}^{d}}\left|\widehat{\eta}\left(\xi\right)\right|>0\text{.}
\]
Then, for every $n\ge3$, 
\[
R_{n}\left({\cal M}_{\varphi,\nu}\right)\gtrsim a_{\eta}^{2}\nu^{d}\frac{\left(\log n\right)^{d/2}}{n}\text{.}
\]
Moreover, every such $\varphi$ belongs to ${\cal L}_{2}\cap{\cal L}_{\infty}\cap{\cal P}$
and satisfies Assumption~\ref{ass:ExponentialFourierKernel} with $\beta=2$. 
\end{thm}

\begin{proof}
See the Appendix. 
\end{proof}

\begin{prop}
\label{prop:GaussianCoreSubmodelRelation}Fix $\nu>0$, and suppose
that 
\[
\varphi=\gamma*\eta
\]
for some Borel probability measure $\eta$ on $\mathbb{R}^{d}$. Then
${\cal M}_{\varphi,\nu}\subseteq{\cal M}_{\gamma,\nu}$. Furthermore,
${\cal M}_{\varphi,\nu}={\cal M}_{\gamma,\nu}$ if and only if $\eta$
is a point mass. 
\end{prop}

\begin{proof}
See the Appendix. 
\end{proof}
\begin{rem}
\label{rem:GaussianCoreExamples}Theorem~\ref{thm:GaussianCoreLowerBound}
contains the Gaussian model itself by taking $\eta=\delta_{0}$. Moreover,
Proposition~\ref{prop:GaussianCoreSubmodelRelation} shows that ${\cal M}_{\gamma*\eta,\nu}\subseteq{\cal M}_{\gamma,\nu}$
for every admissible $\eta$. Hence 
\[
\bigcup_{\eta:\,a_{\eta}>0}{\cal M}_{\gamma*\eta,\nu}={\cal M}_{\gamma,\nu}\text{.}
\]
At the same time, Proposition~\ref{prop:GaussianCoreSubmodelRelation}
shows that whenever $\eta$ is not a point mass, the corresponding
class ${\cal M}_{\gamma*\eta,\nu}$ is a strict subset of ${\cal M}_{\gamma,\nu}$.
Thus Theorem~\ref{thm:GaussianCoreLowerBound} shows that the lower-bound
order from \citet{kim2022minimax} persists on strict submodels
of their framework. 
\end{rem}

\subsection{A lower bound for coordinatewise odd-degree Student-$t$ location
mixtures}

We next report a counterpart of \citet[Thm. 2]{kim2022minimax} for
the odd-degree Student-$t$ family, obtained by tensorizing the one-dimensional
construction. For $r\in\mathbb{N}\cup\left\{ 0\right\} $, let 
\[
\tau_{r}\left(x\right)=\frac{\Gamma\left(r+1\right)}{\sqrt{\pi}\,\Gamma\left(r+\frac{1}{2}\right)}\left(1+x^{2}\right)^{-r-1}\text{,}\qquad x\in\mathbb{R}\text{.}
\]
Thus $\tau_{r}$ is the odd-degree Student-$t$ kernel corresponding
to $2r+1$ degrees of freedom, written in the rescaled form for which
the quadratic term is $1+x^{2}$. If $t_{2r+1}$ denotes the usual
standard Student-$t$ density with $2r+1$ degrees of freedom, then
\[
t_{2r+1}\left(x\right)=\frac{1}{\sqrt{2r+1}}\tau_{r}\left(\frac{x}{\sqrt{2r+1}}\right)=\left(\tau_{r}\right)_{1/\sqrt{2r+1}}\left(x\right)\text{,}\qquad x\in\mathbb{R}\text{.}
\]
For $d\in\mathbb{N}$, define the coordinatewise $d$-variate odd-degree
Student-$t$ kernel by 
\[
\tau_{r,d}\left(x\right)=\prod_{k=1}^{d}\tau_{r}\left(x_{k}\right)\text{,}\qquad x=\left(x_{1},\dots,x_{d}\right)\in\mathbb{R}^{d}\text{.}
\]

\begin{thm}
\label{thm:OddStudentLowerBound}Fix $d\in\mathbb{N}$, $r\in\mathbb{N}\cup\left\{ 0\right\} $,
and $\nu>0$. Then there exists a constant $\kappa_{r,d}>0$ such
that 
\[
R_{n}\left({\cal M}_{\tau_{r,d},\nu}\right)\ge\kappa_{r,d}\nu^{d}\frac{\left(\log n\right)^{d}}{n}
\]
for all sufficiently large $n$. 
\end{thm}

\begin{proof}
See the Appendix. 
\end{proof}
\begin{rem}
\label{rem:OddStudentUpperBound}For each fixed $r$ and $d$, Lemma~\ref{lem:OddStudentBesselFactor}
gives
\[
\widehat{\tau_{r}}(t)=\left(2\pi\right)^{-1/2}e^{-\left|t\right|}q_{r}\left(\left|t\right|\right),
\]
where $q_{r}$ is the degree-$r$ polynomial defined there. Therefore
$\tau_{r,d}\in{\cal L}_{2}\cap{\cal L}_{\infty}\cap{\cal P}$ and 
\[
\left|\widehat{\tau_{r,d}}\left(\xi\right)\right|\lesssim\exp\left(-\left\Vert \xi\right\Vert _{2}/2\right)\text{,}\qquad\xi\in\mathbb{R}^{d}\text{,}
\]
because $\prod_{k=1}^{d}q_{r}\left(\left|\xi_{k}\right|\right)\lesssim\left(1+\left\Vert \xi\right\Vert _{1}\right)^{rd}$
and, for every integer $m\ge0$, there exists a finite constant $C_{m}>0$
such that 
\[
\left(1+t\right)^{m}\le C_{m}e^{t/2}\qquad\text{for all }t\ge0\text{.}
\]
Applying this with $t=\left\Vert \xi\right\Vert _{1}$ and $m=rd$
gives 
\[
\prod_{k=1}^{d}q_{r}\left(\left|\xi_{k}\right|\right)e^{-\left\Vert \xi\right\Vert _{1}}\lesssim e^{-\left\Vert \xi\right\Vert _{1}/2}\le e^{-\left\Vert \xi\right\Vert _{2}/2}
\]
for all $\xi\in\mathbb{R}^{d}$. Hence $\tau_{r,d}$ satisfies Assumption~\ref{ass:ExponentialFourierKernel}
with $\beta=1$, so Theorem~\ref{thm:FastFixedScaleLeastSquares}, with
\[
m_{n}=\left\lceil M\frac{n}{\left(1+\log n\right)^{d}}\right\rceil
\]
for a fixed $M>0$ and $\epsilon_{n}\lesssim\left(1+\log n\right)^{d}/n$, yields 
\[
\sup_{f_{0}\in{\cal M}_{\tau_{r,d},\nu}}\mathrm{E}\left\Vert \hat{f}_{n}-f_{0}\right\Vert _{2}^{2}\lesssim\nu^{d}\frac{\left(1+\log n\right)^{d}}{n}\text{.}
\]
Therefore Theorem~\ref{thm:OddStudentLowerBound} matches the fixed-scale
upper-bound order. 
\end{rem}

\begin{rem}
\label{rem:OddStudentIsotropicOpen}The theorem above proves the lower
bound for the natural tensor-product odd-degree Student-$t$ kernel
$\tau_{r,d}$. The corresponding rotationally invariant $d$-variate
Student-$t$ kernel described in \citet[Sec.~5.5]{kotz2004multivariate}
does not factor coordinatewise and is therefore not covered by our
proof. Thus, we do not claim the analogous sharp lower bound for that
rotationally invariant model here. 
\end{rem}

\section{Discussion}

\label{sec:Discussion}

\paragraph{Lower bounds and limitations for the Fourier-tail class}

Theorem~\ref{thm:GaussianCoreLowerBound} and Theorem~\ref{thm:OddStudentLowerBound}
show that the fixed-scale upper bound from Theorem~\ref{thm:FastFixedScaleLeastSquares}
is minimax-sharp on two canonical families: Gaussian location mixtures (more generally,
the Gaussian convolution submodels from Theorem~\ref{thm:GaussianCoreLowerBound})
and the tensor-product odd-degree Student-$t$ family. In the Gaussian
case, \citet{kim2022minimax} provide the base lower bound and Theorem~\ref{thm:GaussianCoreLowerBound}
transfers it through a bounded deconvolution operator whenever the
extra Fourier factor is bounded away from zero.

A natural question is whether one can obtain a matching lower bound
over the class defined by Assumption~\ref{ass:ExponentialFourierKernel}.
With the current assumptions, the answer is negative. Indeed, Assumption~\ref{ass:ExponentialFourierKernel}
imposes only an upper bound on $\left|\widehat{\varphi}\right|$ and therefore
also permits bandlimited probability kernels. For example, Remark~\ref{rem:BandlimitedMomentExamples}
gives the one-dimensional kernel 
\[
\varphi\left(x\right)=\frac{3}{4\pi}\left(\frac{\sin\left(x/2\right)}{x/2}\right)^{4}\text{,}\qquad x\in\mathbb{R}\text{,}
\]
whose Fourier transform is compactly supported and for which 
\[
\int_{\mathbb{R}}\left|x\right|^{s}\varphi\left(x\right)\mathrm{d}x<\infty\qquad\text{for every }s\in\left(0,3\right)\text{.}
\]
Tensor products therefore provide $d$-dimensional bandlimited kernels
that still satisfy the moment hypothesis from Theorem~\ref{thm:FastGaussianLeastSquares}.
This example already rules out any uniform lower bound with the logarithmic
factor over the whole class defined by Assumption~\ref{ass:ExponentialFourierKernel}.
More generally, if $\operatorname{supp}\left(\widehat{\varphi}\right)\subseteq B\left(0,R_{0}\right)$,
then in the proof of Lemma~\ref{lem:GaussianLocalizedOscillation}
one has $\widehat{h}_{>R}=0$ for every $R>\nu R_{0}$. The high-frequency
tail estimate therefore vanishes, so the localized modulus of continuity
becomes linear in $\delta$, with no logarithmic correction. This
is exactly the mechanism behind the bandlimited parts of Theorem~\ref{thm:FastGaussianLeastSquares}
and Theorem~\ref{thm:FastFixedScaleLeastSquares}. In particular,
the bandlimited part of Theorem~\ref{thm:FastFixedScaleLeastSquares}
gives rate of order $\nu^{d}/n$ over the corresponding fixed-scale
location-mixture class. Thus the minimax risk for that model is at
most of the same order, which is strictly smaller than $\nu^{d}\left(\log n\right)^{d/\beta}/n$
for large $n$. Consequently, no lower bound with the logarithmic
factor can hold uniformly over all kernels satisfying only Assumption~\ref{ass:ExponentialFourierKernel}.

If one strengthens Assumption~\ref{ass:ExponentialFourierKernel} to a two-sided
supersmooth condition (cf. \citealt[Sec. 2.4.3]{meister2009deconvolution})
of the form 
\[
c_{1}\exp\left(-C_{1}\left\Vert \xi\right\Vert _{2}^{\beta}\right)\le\left|\widehat{\varphi}\left(\xi\right)\right|\le C_{2}\exp\left(-c_{2}\left\Vert \xi\right\Vert _{2}^{\beta}\right)
\]
for all sufficiently large $\left\Vert \xi\right\Vert _{2}$, then
compact Fourier support is excluded. In this restricted setting, the
difficulty is not the transfer step. Indeed, once a canonical $\beta$-supersmooth
base kernel $\psi_{\beta}$ is available with a lower bound of order
$\nu^{d}\left(\log n\right)^{d/\beta}/n$ over the corresponding fixed-scale
location-mixture class, any representation 
\[
\varphi=\psi_{\beta}*\eta
\]
with $\eta$ a Borel probability measure on $\mathbb{R}^{d}$ satisfying
\[
\inf_{\xi\in\mathbb{R}^{d}}|\widehat{\eta}(\xi)|>0
\]
would transfer that lower bound exactly as in Theorem~\ref{thm:GaussianCoreLowerBound}.
Thus the unresolved step is the base theorem itself. For $\beta=2$,
the Gaussian kernel provides such a base model and the Hermite-based
Assouad construction of \citet{kim2022minimax} yields the required
lower bound. At exponential order, Theorem~\ref{thm:OddStudentLowerBound}
shows that the tensor-product odd-degree Student-$t$ family also
provides a base model because its Fourier multiplier differs from
the product Cauchy multiplier only by a fixed polynomial factor. Beyond
these structured examples, however, the missing proof step is an appropriate
$\beta$-supersmooth mixture class with a sharp lower bound of order
$\nu^{d}\left(\log n\right)^{d/\beta}/n$. In the tensor-product odd-degree
Student-$t$ case, the Fourier transform factorizes as 
\[
\widehat{\tau_{r,d}}\left(\xi\right)=\left(2\pi\right)^{-d/2}e^{-\left\Vert \xi\right\Vert _{1}}\prod_{k=1}^{d}q_{r}\left(\left|\xi_{k}\right|\right)\text{,}
\]
so dividing by $\prod_{k=1}^{d}q_{r}\left(\left|\xi_{k}\right|\right)$
permits diagonalization. The inverse multiplier remains tractable
due to Lemma~\ref{lem:OddStudentBesselFactor}(ii). This mechanism
is specific to the odd-degree Student-$t$ kernel. We do not currently
know of any approaches that allow us to prove lower bounds, even in
apparently basic cases such as even-degree Student-$t$ kernels or
$\alpha$ stable kernels.

\paragraph{Beyond ${\cal L}_{p}$ spaces}

One immediate extension of the approximation corollaries from Section~\ref{sec:Approximation-rates-for}
is to consider approximation with respect to norms of the form $f\mapsto\left\Vert f\right\Vert =\sum_{p\in\mathbb{P}}\left\Vert f\right\Vert _{p}$
for some finite set $\mathbb{P}\subset\left(1,\infty\right)$. For
such norms, one expects the slowest rate among the constituent exponents
from Corollary~\ref{cor:GeneralApproxRate} to govern. Establishing
this within the present framework would require choosing a single
approximating sequence that works simultaneously for all $p\in\mathbb{P}$.

Less trivially, we can consider two natural extensions of ${\cal L}_{p}$
spaces. The first extension is to consider the weighted Lebesgue spaces
${\cal L}_{p,w}$, whose norm is determined by a positive function
$w:\mathbb{R}^{d}\to\mathbb{R}_{\ge0}$ via the expression $f\mapsto\left\Vert f\right\Vert _{p,w}=\left[\int_{\mathbb{R}^{d}}\left|f\left(x\right)\right|^{p}w\left(x\right)\mathrm{d}x\right]^{1/p}$
for each $p\in\left[1,\infty\right)$, as discussed in \citet[Sec. 3.8]{castillo2016introductory}
and appearing often in the approximation theory literature (see, e.g.,
\citealp{Mhaskar1996}). The functional-analytic ingredients needed
in this setting are the weighted counterparts of the duality and integral
inequalities used in the unweighted case, together with the standard
theory of approximate identities in weighted ${\cal L}_{p}$ spaces.
Under hypotheses on the weight that guarantee these tools, it is plausible
that analogues of the approximation results from Section~\ref{sec:Approximation-rates-for}
can be proved by adapting the same strategy.

The second extension is to consider variable exponent Lebesgue spaces
${\cal L}_{p\left(\cdot\right)}$ of functions $f$ satisfying $\int_{\mathbb{R}^{d}}\left|f\left(x\right)\right|^{p\left(x\right)}\mathrm{d}x<\infty$,
for $p:\mathbb{R}^{d}\to\left[1,\infty\right]$ and whose norm naturally
extends the usual ${\cal L}_{p}$ norms when $p\in\left[1,\infty\right]$
is constant (cf. \citealp{diening2011lebesgue} and \citealp{cruz2013variable}).
Under appropriate assumptions, these spaces support convolution-based
approximation and suggest a possible variable-exponent analogue of
the present results.

Lastly, qualitative convolution-based approximation theorems are also
available for Sobolev spaces (cf. \citealp{Leoni2017SobolevSpaces}
and \citealp{adams2003sobolev}). After the corresponding modifications
to the arguments of \citet{nguyen2023approximation}, these results
should yield qualitative mixture-approximation theorems in Sobolev
settings. Quantitative analogues of Corollary~\ref{cor:GeneralApproxRate},
however, appear to require additional technical input beyond the scope
of the present manuscript.

\section*{Appendix}

\subsection*{Auxiliary technical results}
\begin{lem}
\label{lem:translation-continuity-Lp}Let $p\in\left[1,\infty\right)$
and $f\in{\cal L}_{p}$. Then 
\[
\lim_{\left\Vert y\right\Vert _{2}\to0}\left\Vert f\left(\cdot-y\right)-f\right\Vert _{p}=0\text{.}
\]
\end{lem}

\begin{proof}
Fix $\epsilon>0$. By the density of ${\cal C}_{c}^{\infty}\left(\mathbb{R}^{d}\right)$
in ${\cal L}_{p}$, choose $g\in{\cal C}_{c}^{\infty}\left(\mathbb{R}^{d}\right)$
such that 
\[
\left\Vert f-g\right\Vert _{p}<\frac{\epsilon}{3}\text{.}
\]
Choose $R>0$ so large that $\operatorname{supp}\left(g\right)\subseteq\left[-R,R\right]^{d}$.
Since $g$ is uniformly continuous, there exists $\delta\in\left(0,1\right)$
such that 
\[
\sup_{x\in\mathbb{R}^{d}}\left|g\left(x-y\right)-g\left(x\right)\right|<\frac{\epsilon}{3\left(2R+2\right)^{d/p}}
\]
whenever $\left\Vert y\right\Vert _{2}<\delta$. For such $y$, the
function $x\mapsto g\left(x-y\right)-g\left(x\right)$ is supported
in $\left[-R-1,R+1\right]^{d}$, so 
\[
\left\Vert g\left(\cdot-y\right)-g\right\Vert _{p}\le\left(2R+2\right)^{d/p}\sup_{x\in\mathbb{R}^{d}}\left|g\left(x-y\right)-g\left(x\right)\right|<\frac{\epsilon}{3}\text{.}
\]
By translation-invariance of the ${\cal L}_{p}$ norm, 
\[
\left\Vert f\left(\cdot-y\right)-f\right\Vert _{p}\le\left\Vert f\left(\cdot-y\right)-g\left(\cdot-y\right)\right\Vert _{p}+\left\Vert g\left(\cdot-y\right)-g\right\Vert _{p}+\left\Vert g-f\right\Vert _{p}<\epsilon
\]
whenever $\left\Vert y\right\Vert _{2}<\delta$. This proves the claim. 
\end{proof}
\begin{lem}
\label{lem:finite-supported-weak-density}For every Borel probability
measure $G$ on $\mathbb{R}^{d}$, there exists a sequence $\left(G_{k}\right)_{k}$
of finitely supported Borel probability measures on $\mathbb{R}^{d}$
such that $G_{k}$ converges weakly to $G$. 
\end{lem}

\begin{proof}
For each $k\in\mathbb{N}$, let $Q_{k}=\left[-k,k\right]^{d}$. Partition
$Q_{k}$ into finitely many half-open cubes $A_{k,1},\dots,A_{k,N_{k}}$,
each having diameter strictly less than $1/k$, and choose $y_{k,j}\in A_{k,j}$
for each $j\in\left[N_{k}\right]$. Set 
\[
A_{k,0}=\mathbb{R}^{d}\setminus Q_{k},
\]
choose any $y_{k,0}\in\mathbb{R}^{d}$, and define 
\[
G_{k}=\sum_{j=0}^{N_{k}}G\left(A_{k,j}\right)\delta_{y_{k,j}}\text{.}
\]
Then each $G_{k}$ is a finitely supported Borel probability measure.

Let $h:\mathbb{R}^{d}\to\mathbb{R}$ be bounded and continuous, and
fix $\epsilon>0$. If $\left\Vert h\right\Vert _{\infty}=0$, then
the claim is immediate. So assume $\left\Vert h\right\Vert _{\infty}>0$
and choose $R>0$ such that 
\[
G\left(\mathbb{R}^{d}\setminus Q_{R}\right)<\frac{\epsilon}{4\left\Vert h\right\Vert _{\infty}}\text{.}
\]
Since $h$ is uniformly continuous on the compact set $Q_{R+1}$,
there exists $\delta>0$ such that 
\[
\left|h\left(x\right)-h\left(y\right)\right|<\frac{\epsilon}{2}
\]
whenever $x,y\in Q_{R+1}$ and $\left\Vert x-y\right\Vert _{2}<\delta$.
Choose $k$ so large that $k\ge R+1$ and $1/k<\delta$. Then 
\begin{align*}
\left|\int_{\mathbb{R}^{d}}h\left(x\right)G_{k}\left(\mathrm{d}x\right)-\int_{\mathbb{R}^{d}}h\left(x\right)G\left(\mathrm{d}x\right)\right| & =\left|\sum_{j=0}^{N_{k}}\int_{A_{k,j}}\left[h\left(y_{k,j}\right)-h\left(x\right)\right]G\left(\mathrm{d}x\right)\right|\\
 & \le\sum_{j=0}^{N_{k}}\int_{A_{k,j}\cap Q_{R}}\left|h\left(y_{k,j}\right)-h\left(x\right)\right|G\left(\mathrm{d}x\right)\\
 & \quad+\sum_{j=0}^{N_{k}}\int_{A_{k,j}\cap\left(\mathbb{R}^{d}\setminus Q_{R}\right)}\left|h\left(y_{k,j}\right)-h\left(x\right)\right|G\left(\mathrm{d}x\right)\text{.}
\end{align*}
If $x\in A_{k,j}\cap Q_{R}$, then $y_{k,j}\in Q_{R+1}$ and $\left\Vert x-y_{k,j}\right\Vert _{2}<\delta$,
so the first sum is at most $\epsilon/2$. The second sum is bounded
by 
\[
2\left\Vert h\right\Vert _{\infty}G\left(\mathbb{R}^{d}\setminus Q_{R}\right)<\frac{\epsilon}{2}\text{.}
\]
Hence 
\[
\left|\int_{\mathbb{R}^{d}}h\left(x\right)G_{k}\left(\mathrm{d}x\right)-\int_{\mathbb{R}^{d}}h\left(x\right)G\left(\mathrm{d}x\right)\right|<\epsilon\text{.}
\]
Since $h$ was an arbitrary bounded continuous function, it follows
that $G_{k}$ converges weakly to $G$. 
\end{proof}
\begin{lem}[Convex approximation theorem]
\label{lem:convex-approximation}Let ${\cal A}\subset{\cal L}_{p}$
for $p\in\left(1,\infty\right)$ and suppose that $f\in\overline{\mathrm{co}\left({\cal A}\right)}^{p}$.
Further, assume that there exists a $B>0$ such that for all $g\in{\cal A}$,
$\left\Vert g-f\right\Vert _{p}\le B$. Then, for each $\epsilon>0$,
there exists a sequence $\left(g_{m}\right)_{m}$ with $g_{m}\in\mathrm{co}_{m}\left({\cal A}\right)$,
for each $m\in\mathbb{N}$, such that 
\[
\left\Vert g_{m}-f\right\Vert _{p}\le\begin{cases}
\left(B+\epsilon\right)C_{p}m^{-1/q} & \text{if }1<p<2\text{,}\\
\left(B+\epsilon\right)C_{p}m^{-1/2} & \text{if }2\le p<\infty\text{,}
\end{cases}
\]
where $C_{p}<\infty$ is a fixed constant depending only on $p$ and
$1/p+1/q=1$.
\end{lem}

\begin{lem}[Minkowski's integral inequality]
\label{lem:minkowski-integral-inequality}Let $\left(\mathbb{X},\mathfrak{X},\mu_{1}\right)$
and $\left(\mathbb{Y},\mathfrak{Y},\mu_{2}\right)$ be $\sigma$-finite
measure spaces, and let $f:\mathbb{X}\times\mathbb{Y}\to\bar{\mathbb{R}}$
be $\mathfrak{X}\otimes\mathfrak{Y}$-measurable. Then, for each $p\in\left[1,\infty\right)$,
\[
\left[\int_{\mathbb{X}}\left[\int_{\mathbb{Y}}\left|f\left(x,y\right)\right|\mu_{2}\left(\mathrm{d}y\right)\right]^{p}\mu_{1}\left(\mathrm{d}x\right)\right]^{1/p}\le\int_{\mathbb{Y}}\left[\int_{\mathbb{X}}\left|f\left(x,y\right)\right|^{p}\mu_{1}\left(\mathrm{d}x\right)\right]^{1/p}\mu_{2}\left(\mathrm{d}y\right)\text{.}
\]
\end{lem}

\begin{lem}
\label{lem:smoothness-W1p}Let $p\in\left[1,\infty\right)$ and $f\in{\cal W}^{1,p}$.
Then, for all $y\in\mathbb{R}^{d}$, 
\[
\left[\int_{\mathbb{R}^{d}}\left|f\left(x-y\right)-f\left(x\right)\right|^{p}\mathrm{d}x\right]^{1/p}\le\left[\int_{\mathbb{R}^{d}}\left\Vert \nabla f\left(x\right)\right\Vert _{2}^{p}\mathrm{d}x\right]^{1/p}\left\Vert y\right\Vert _{2}\text{.}
\]
\end{lem}

\begin{lem}
\label{lem:smoothness-Wsp}Let $p\in\left[1,\infty\right)$, $s\in\left(0,1\right)$,
and $f\in{\cal W}^{s,p}$. Then there exists a finite constant
$C_{d,s,p}>0$ such that, for any $y\in\mathbb{R}^{d}$, 
\[
\left[\int_{\mathbb{R}^{d}}\left|f\left(x-y\right)-f\left(x\right)\right|^{p}\mathrm{d}x\right]^{1/p}\le C_{d,s,p}\left\Vert f\right\Vert _{s,p}\left\Vert y\right\Vert _{2}^{s}\text{.}
\]
\end{lem}

\begin{lem}
\label{lem:VCBoundMaximalIneq}Let ${\cal F}\subset\ell^{\infty}\left(\mathbb{R}^{d}\right)$
be a pointwise measurable VC-subgraph class with measurable envelope
$\bar{F}$ satisfying $P\bar{F}^{2}<\infty$. Then, 
\[
\mathrm{E}\sup_{f\in{\cal F}}\left|P_{n}f-Pf\right|\lesssim\sqrt{\frac{\operatorname{VC}\left({\cal F}\right)}{n}}\left[P\bar{F}^{2}\right]^{1/2}\text{.}
\]
\end{lem}

\begin{lem}[Assouad's lemma]
\label{lem:Assouad}Let ${\cal F}\subset{\cal L}_{2}\cap{\cal P}$
be a class of densities on $\mathbb{R}^{d}$. Suppose that $\left\{ f_{\alpha}:\alpha\in\left\{ 0,1\right\} ^{N}\right\} \subseteq{\cal F}$.
Write 
\[
\Upsilon\left(\alpha,\beta\right)=\sum_{k=1}^{N}\mathbf{1}\left\{ \alpha_{k}\neq\beta_{k}\right\} 
\]
for the Hamming distance. Then 
\[
R_{n}\left({\cal F}\right)\ge\frac{N}{8}\min_{\alpha\neq\beta}\frac{\left\Vert f_{\alpha}-f_{\beta}\right\Vert _{2}^{2}}{\Upsilon\left(\alpha,\beta\right)}\min_{\Upsilon\left(\alpha,\beta\right)=1}\left(1-\sqrt{\frac{n}{2}\chi^{2}\left(f_{\alpha}\|f_{\beta}\right)}\right)
\]
where 
\[
\chi^{2}\left(f\|g\right)=\int_{\mathbb{R}^{d}}\frac{\left(f\left(x\right)-g\left(x\right)\right)^{2}}{g\left(x\right)}\mathrm{d}x\text{.}
\]
\end{lem}

\begin{lem}[Least-squares basic inequality]
\label{lem:LeastSquaresBasicIneq}Let ${\cal F}\subset{\cal L}_{2}\cap{\cal P}$,
let $f_{0}\in{\cal L}_{2}\cap{\cal P}$, let $\epsilon\ge0$, and
assume that $X_{1},\dots,X_{n}$ are i.i.d.\ with common density
$f_{0}$. Suppose that $\hat{f}_{n,\epsilon}\in{\cal F}$ satisfies
(\ref{eq:sample-eps-minimizer-def}). Then, for every $f^{*}\in{\cal F}$,
\[
\left\Vert \hat{f}_{n,\epsilon}-f_{0}\right\Vert _{2}^{2}\le\left\Vert f^{*}-f_{0}\right\Vert _{2}^{2}+2\left(P_{n}-P\right)\left(\hat{f}_{n,\epsilon}-f^{*}\right)+\epsilon\text{.}
\]
\end{lem}

\begin{proof}
Since 
\[
\left\Vert f-f_{0}\right\Vert _{2}^{2}=\int_{\mathbb{R}^{d}}f^{2}\left(x\right)\mathrm{d}x-2Pf+\int_{\mathbb{R}^{d}}f_{0}^{2}\left(x\right)\mathrm{d}x\text{,}
\]
we obtain 
\begin{align*}
\left\Vert \hat{f}_{n,\epsilon}-f_{0}\right\Vert _{2}^{2}-\left\Vert f^{*}-f_{0}\right\Vert _{2}^{2} & =\left[\int_{\mathbb{R}^{d}}\hat{f}_{n,\epsilon}^{2}\left(x\right)\mathrm{d}x-2P\hat{f}_{n,\epsilon}\right]-\left[\int_{\mathbb{R}^{d}}f^{*2}\left(x\right)\mathrm{d}x-2Pf^{*}\right]\\
 & =\left[\int_{\mathbb{R}^{d}}\hat{f}_{n,\epsilon}^{2}\left(x\right)\mathrm{d}x-2P_{n}\hat{f}_{n,\epsilon}\right]-\left[\int_{\mathbb{R}^{d}}f^{*2}\left(x\right)\mathrm{d}x-2P_{n}f^{*}\right]\\
 & \qquad+2\left(P_{n}-P\right)\left(\hat{f}_{n,\epsilon}-f^{*}\right)\\
 & \le\epsilon+2\left(P_{n}-P\right)\left(\hat{f}_{n,\epsilon}-f^{*}\right)\text{,}
\end{align*}
where the last step is exactly (\ref{eq:sample-eps-minimizer-def}).
This proves the claim. 
\end{proof}
\begin{thm}[Oracle inequality for bounded least-squares density estimation]
\label{thm:MassartNedelecLeastSquares}Let ${\cal F}\subset{\cal L}_{2}\cap{\cal L}_{\infty}\cap{\cal P}$
and assume that there exists a countable subclass ${\cal F}^{0}\subset{\cal F}$
such that every $f\in{\cal F}$ admits a sequence $\left(f_{k}\right)_{k}\subset{\cal F}^{0}$
with 
\[
f_{k}\left(x\right)\to f\left(x\right)\text{ for every }x\in\mathbb{R}^{d}\qquad\text{and}\qquad\left\Vert f_{k}-f\right\Vert _{2}\to0\text{.}
\]
Let $f_{0}\in{\cal L}_{\infty}\cap{\cal P}$, assume that $X_{1},\dots,X_{n}$
are i.i.d.\ with common density $f_{0}$, and set 
\[
M=\left\Vert f_{0}\right\Vert _{\infty}\vee\sup_{f\in{\cal F}}\left\Vert f\right\Vert _{\infty}<\infty\text{.}
\]
Suppose that $\hat{f}_{n,\epsilon}\in{\cal F}$ satisfies (\ref{eq:sample-eps-minimizer-def}),
and assume that there exists $A\ge1$ such that, for every $u\in{\cal F}$
and every $\sigma\ge A/\sqrt{n}$, 
\[
\sqrt{n}\,\mathrm{E}^{*}\sup_{\substack{f\in{\cal F}\\
\left\Vert f-u\right\Vert _{2}\le2\sqrt{M}\,\sigma
}
}\left|\left(P_{n}-P\right)\left(f-u\right)\right|\le2MA\sigma\text{.}
\]
Then there exists a finite universal constant $C>0$ such that 
\[
\mathrm{E}\left\Vert \hat{f}_{n,\epsilon}-f_{0}\right\Vert _{2}^{2}\le2\inf_{f\in{\cal F}}\left\Vert f-f_{0}\right\Vert _{2}^{2}+2\epsilon+CM\frac{A^{2}}{n}\text{.}
\]
The requirement $A\ge1$ is harmless, since enlarging $A$ preserves
all of the hypotheses. 
\end{thm}

\begin{proof}
Define, for $f\in{\cal F}\cup\left\{ f_{0}\right\} $, 
\[
\widetilde{\gamma}_{f}\left(x\right)=\frac{1}{4M}\left[\int_{\mathbb{R}^{d}}f^{2}\left(y\right)\mathrm{d}y-2f\left(x\right)+2M\right]\text{.}
\]
We verify, one by one, the hypotheses of Theorem~2 of \citet{massart2006risk},
taking the ambient index class to be ${\cal T}={\cal F}\cup\left\{ f_{0}\right\} $
and the model class to be ${\cal S}_{0}={\cal F}$.

First, since $0\le f\left(x\right)\le M$ and $\left\Vert f\right\Vert _{2}^{2}\le\left\Vert f\right\Vert _{\infty}\left\Vert f\right\Vert _{1}\le M$
for every $f\in{\cal F}\cup\left\{ f_{0}\right\} $, we have 
\[
0\le\widetilde{\gamma}_{f}\left(x\right)\le\frac{M+2M}{4M}=\frac{3}{4}\le1\qquad\text{for all }x\in\mathbb{R}^{d}\text{.}
\]
Thus the contrast is bounded in $\left[0,1\right]$.

Second, for every $f\in{\cal F}$, 
\[
P\widetilde{\gamma}_{f}-P\widetilde{\gamma}_{f_{0}}=\frac{1}{4M}\left\Vert f-f_{0}\right\Vert _{2}^{2}\text{.}
\]
Hence $f_{0}$ is the population minimizer of the contrast over the
ambient index class ${\cal T}$, and $\hat{f}_{n,\epsilon}$ is a
$\rho$-empirical risk minimizer for this contrast with 
\[
\rho=\frac{\epsilon}{4M}\text{.}
\]

Third, define the pseudo-distance 
\[
d\left(f,g\right)=\frac{\left\Vert f-g\right\Vert _{2}}{2\sqrt{M}}\text{.}
\]
Then 
\[
d\left(f_{0},f\right)=\sqrt{P\widetilde{\gamma}_{f}-P\widetilde{\gamma}_{f_{0}}}\text{,}
\]
so the margin requirement in Theorem 2 of \citet{massart2006risk}
holds with $w\left(x\right)=x$. Moreover, 
\[
\widetilde{\gamma}_{f}\left(x\right)-\widetilde{\gamma}_{g}\left(x\right)=\frac{\left\Vert f\right\Vert _{2}^{2}-\left\Vert g\right\Vert _{2}^{2}}{4M}-\frac{f\left(x\right)-g\left(x\right)}{2M}\text{,}
\]
so the first term on the right-hand side is deterministic and therefore
\begin{align*}
\operatorname{Var}\left[\widetilde{\gamma}_{f}\left(X\right)-\widetilde{\gamma}_{g}\left(X\right)\right] & =\operatorname{Var}\left[\frac{f\left(X\right)-g\left(X\right)}{2M}\right]\\
 & \le\frac{1}{4M^{2}}\mathrm{E}\left[\left(f\left(X\right)-g\left(X\right)\right)^{2}\right]\\
 & =\frac{1}{4M^{2}}\int_{\mathbb{R}^{d}}\left(f\left(x\right)-g\left(x\right)\right)^{2}f_{0}\left(x\right)\mathrm{d}x\\
 & \le\frac{\left\Vert f-g\right\Vert _{2}^{2}}{4M}=d\left(f,g\right)^{2}\text{.}
\end{align*}
Thus the variance condition of \citet[Thm.~2]{massart2006risk} is
satisfied.

Fourth, we check the separability condition $(M)$ from \citet{massart2006risk}
for the model ${\cal F}$. If $f\in{\cal F}$ and $f_{k}\in{\cal F}^{0}$
converges pointwise to $f$ with $\left\Vert f_{k}-f\right\Vert _{2}\to0$,
then also $\left\Vert f_{k}\right\Vert _{2}^{2}\to\left\Vert f\right\Vert _{2}^{2}$,
so for every $x\in\mathbb{R}^{d}$, 
\[
\widetilde{\gamma}_{f_{k}}\left(x\right)\to\widetilde{\gamma}_{f}\left(x\right)\text{.}
\]
Therefore the contrast class indexed by the model ${\cal F}$ satisfies
condition $(M)$ with countable core ${\cal F}^{0}$.

Finally, let $u\in{\cal F}^{0}$ and let $\sigma\ge A/\sqrt{n}$.
Since ${\cal F}^{0}$ is countable, the localized supremum over ${\cal F}^{0}$
that appears in condition (20) of \citet[Thm.~2]{massart2006risk}
is measurable. Moreover, it is bounded by the corresponding supremum
over ${\cal F}$ from the hypothesis of the present theorem. Since
\[
\left(P_{n}-P\right)\widetilde{\gamma}_{f}-\left(P_{n}-P\right)\widetilde{\gamma}_{u}=-\frac{1}{2M}\left(P_{n}-P\right)\left(f-u\right)\text{,}
\]
and since $d\left(f,u\right)\le\sigma$ is equivalent to $\left\Vert f-u\right\Vert _{2}\le2\sqrt{M}\,\sigma$,
the assumed bound implies 
\begin{align*}
\sqrt{n}\,\mathrm{E}\sup_{\substack{f\in{\cal F}^{0}\\
d\left(f,u\right)\le\sigma
}
}\left[\left(P_{n}-P\right)\widetilde{\gamma}_{u}-\left(P_{n}-P\right)\widetilde{\gamma}_{f}\right] & \le\frac{\sqrt{n}}{2M}\,\mathrm{E}^{*}\sup_{\substack{f\in{\cal F}\\
\left\Vert f-u\right\Vert _{2}\le2\sqrt{M}\,\sigma
}
}\left|\left(P_{n}-P\right)\left(f-u\right)\right|\\
 & \le A\sigma\text{.}
\end{align*}
Thus condition (20) of \citet[Thm.~2]{massart2006risk} holds with
$\phi\left(\sigma\right)=A\sigma$. The equation defining their $\epsilon_{*}$
becomes 
\[
\sqrt{n}\,\epsilon_{*}^{2}=\phi\left(w\left(\epsilon_{*}\right)\right)=A\epsilon_{*}\text{,}
\]
so $\epsilon_{*}=A/\sqrt{n}$.

Applying \citet[Thm.~2]{massart2006risk} therefore yields 
\[
\frac{1}{4M}\mathrm{E}\left\Vert \hat{f}_{n,\epsilon}-f_{0}\right\Vert _{2}^{2}\le2\left(\rho+\frac{1}{4M}\inf_{f\in{\cal F}}\left\Vert f-f_{0}\right\Vert _{2}^{2}+k\frac{A^{2}}{n}\right)
\]
for a finite universal constant $k>0$. Multiplying by $4M$ proves
the claim with $C=8k$. 
\end{proof}

\subsection*{Technical proofs}
\begin{proof}[Proof of Lemma~\ref{lem:GaussianSieveSeparability}]
Let ${\cal G}_{m,\nu}^{0}$ be the subclass of mixtures 
\[
f\left(x\right)=\sum_{j=1}^{m}\pi_{j}\varphi_{\nu}\left(x-\mu_{j}\right)
\]
for which each $\mu_{j}$ belongs to $\mathbb{Q}^{d}$ and the weight
vector $\left(\pi_{1},\dots,\pi_{m}\right)$ belongs to the rational
simplex 
\[
\Delta_{m}^{\mathbb{Q}}=\left\{ \left(\pi_{1},\dots,\pi_{m}\right)\in\mathbb{Q}^{m}\text{ | }\pi_{j}\ge0\text{ and }\sum_{j=1}^{m}\pi_{j}=1\right\} \text{.}
\]
Then ${\cal G}_{m,\nu}^{0}$ is countable. Fix any 
\[
f=\sum_{j=1}^{m}\pi_{j}\varphi_{\nu}\left(\cdot-\mu_{j}\right)\in{\cal G}_{m,\nu}\text{.}
\]
Choose sequences $\left(\pi_{1,k},\dots,\pi_{m,k}\right)\in\Delta_{m}^{\mathbb{Q}}$
and $\mu_{1,k},\dots,\mu_{m,k}\in\mathbb{Q}^{d}$ such that $\pi_{j,k}\to\pi_{j}$
and $\mu_{j,k}\to\mu_{j}$ for every $j\in\left[m\right]$, and set
\[
f_{k}=\sum_{j=1}^{m}\pi_{j,k}\varphi_{\nu}\left(\cdot-\mu_{j,k}\right)\in{\cal G}_{m,\nu}^{0}\text{.}
\]
Since $\varphi$ is continuous, $f_{k}\left(x\right)\to f\left(x\right)$
for every $x\in\mathbb{R}^{d}$. Moreover, 
\begin{align*}
\left\Vert f_{k}-f\right\Vert _{2} & \le\sum_{j=1}^{m}\left|\pi_{j,k}-\pi_{j}\right|\left\Vert \varphi_{\nu}\right\Vert _{2}+\sum_{j=1}^{m}\pi_{j}\left\Vert \varphi_{\nu}\left(\cdot-\mu_{j,k}\right)-\varphi_{\nu}\left(\cdot-\mu_{j}\right)\right\Vert _{2}\text{.}
\end{align*}
The first sum tends to zero because $\pi_{j,k}\to\pi_{j}$. The second
sum tends to zero because translations are continuous in ${\cal L}_{2}$
and $\mu_{j,k}\to\mu_{j}$ for every $j$. Hence $\left\Vert f_{k}-f\right\Vert _{2}\to0$,
as required. 
\end{proof}
\begin{proof}[Proof of Lemma~\ref{lem:GaussianLocalizedOscillation}]
We first prove that there exist finite constants
$C_{d,\varphi}>0$, depending only on $d$, $\left\Vert \varphi\right\Vert _{\infty}$,
and the constants in Assumption~\ref{ass:ExponentialFourierKernel}, and, in the
bandlimited case, $C_{d,\varphi,R_{0}}^{\mathrm{bl}}>0$ such that,
for every $m\in\mathbb{N}$, every $\nu>0$, every $u\in{\cal G}_{m,\nu}$,
and every $\delta>0$, 
\[
\sqrt{n}\,\mathrm{E}^{*}\sup_{\substack{f\in{\cal G}_{m,\nu}\\
\left\Vert f-u\right\Vert _{2}\le\delta
}
}\left|\left(P_{n}-P\right)\left(f-u\right)\right|\le C_{d,\varphi}\sqrt{M_{\nu}}\,\delta\left[1+\left(1+\log_{+}\frac{C_{d,\varphi}M_{\nu}}{\delta^{2}}\right)^{d/\left(2\beta\right)}\right]\text{,}
\]
and, if $\operatorname{supp}\left(\widehat{\varphi}\right)\subseteq B\left(0,R_{0}\right)$,
\[
\sqrt{n}\,\mathrm{E}^{*}\sup_{\substack{f\in{\cal G}_{m,\nu}\\
\left\Vert f-u\right\Vert _{2}\le\delta
}
}\left|\left(P_{n}-P\right)\left(f-u\right)\right|\le C_{d,\varphi,R_{0}}^{\mathrm{bl}}\sqrt{M_{\nu}}\,\delta\text{.}
\]
Taking for granted these bounds for the moment, the stated lemma follows by
taking $\delta=2\sqrt{M_{\nu}}\,\sigma$. Indeed, the general case
then gives 
\begin{align*}
 & \sqrt{n}\,\mathrm{E}^{*}\sup_{\substack{f\in{\cal G}_{m,\nu}\\
\left\Vert f-u\right\Vert _{2}\le2\sqrt{M_{\nu}}\,\sigma
}
}\left|\left(P_{n}-P\right)\left(f-u\right)\right|\\
 & \qquad\le2C_{d,\varphi}M_{\nu}\sigma\left[1+\left(1+\log_{+}\frac{C_{d,\varphi}}{4\sigma^{2}}\right)^{d/\left(2\beta\right)}\right]\text{.}
\end{align*}
If $\sigma\ge A_{n}/\sqrt{n}$, then 
\[
\log_{+}\frac{C_{d,\varphi}}{4\sigma^{2}}\le\log_{+}\left(\frac{C_{d,\varphi}n}{4A_{n}^{2}}\right)\le C_{d,\varphi}^{\prime}\log n
\]
for some finite constant $C_{d,\varphi}^{\prime}>0$ depending only
on $d$ and the constants in Assumption~\ref{ass:ExponentialFourierKernel}. Hence
\[
1+\left(1+\log_{+}\frac{C_{d,\varphi}}{4\sigma^{2}}\right)^{d/\left(2\beta\right)}\le C_{d,\varphi}^{\prime\prime}\left(1+\left(\log n\right)^{d/\left(2\beta\right)}\right)
\]
for some finite constant $C_{d,\varphi}^{\prime\prime}>0$ depending
only on $d$ and the constants in Assumption~\ref{ass:ExponentialFourierKernel}.
Choosing $K_{d,\varphi}\ge1\vee2C_{d,\varphi}C_{d,\varphi}^{\prime\prime}$
proves the first assertion.

If, in addition, $\operatorname{supp}\left(\widehat{\varphi}\right)\subseteq B\left(0,R_{0}\right)$,
then the bandlimited estimate above yields 
\begin{align*}
 & \sqrt{n}\,\mathrm{E}^{*}\sup_{\substack{f\in{\cal G}_{m,\nu}\\
\left\Vert f-u\right\Vert _{2}\le2\sqrt{M_{\nu}}\,\sigma
}
}\left|\left(P_{n}-P\right)\left(f-u\right)\right|\\
 & \qquad\le2C_{d,\varphi,R_{0}}^{\mathrm{bl}}M_{\nu}\sigma\text{.}
\end{align*}
Thus the same conclusion holds with $K_{d,\varphi,R_{0}}^{\mathrm{bl}}=1\vee2C_{d,\varphi,R_{0}}^{\mathrm{bl}}$.

We now verify the displayed bounds. Fix $m\in\mathbb{N}$,
$\nu>0$, $u\in{\cal G}_{m,\nu}$, and $\delta>0$. Set 
\[
{\cal H}_{\delta}\left(u\right)=\left\{ h=f-u\text{ | }f\in{\cal G}_{m,\nu}\text{ and }\left\Vert h\right\Vert _{2}\le\delta\right\} \text{.}
\]
Let ${\cal G}_{m,\nu}^{0}\subset{\cal G}_{m,\nu}$ be the countable
subclass from Lemma~\ref{lem:GaussianSieveSeparability}. For $\eta>0$,
define 
\[
{\cal H}_{\delta,\eta}^{0}\left(u\right)=\left\{ h=f-u\text{ | }f\in{\cal G}_{m,\nu}^{0}\text{ and }\left\Vert h\right\Vert _{2}\le\delta+\eta\right\} \text{.}
\]
We first reduce the proof to the countable class ${\cal H}_{\delta,\eta}^{0}\left(u\right)$.
Take any $h=f-u\in{\cal H}_{\delta}\left(u\right)$. By the explicit
construction in the proof of Lemma~\ref{lem:GaussianSieveSeparability},
we may write 
\[
f=\sum_{j=1}^{m}\pi_{j}\varphi_{\nu}\left(\cdot-\mu_{j}\right)
\]
and choose $f_{k}\in{\cal G}_{m,\nu}^{0}$ of the form 
\[
f_{k}=\sum_{j=1}^{m}\pi_{j,k}\varphi_{\nu}\left(\cdot-\mu_{j,k}\right)
\]
with $\pi_{j,k}\to\pi_{j}$ and $\mu_{j,k}\to\mu_{j}$ for every $j\in\left[m\right]$.
Then $f_{k}\left(x\right)\to f\left(x\right)$ for every $x\in\mathbb{R}^{d}$
and $\left\Vert f_{k}-f\right\Vert _{2}\to0$. In particular, 
\[
\left\Vert f_{k}-u\right\Vert _{2}\le\left\Vert f_{k}-f\right\Vert _{2}+\left\Vert f-u\right\Vert _{2}\to\left\Vert f-u\right\Vert _{2}\le\delta,
\]
so $h_{k}=f_{k}-u$ belongs to ${\cal H}_{\delta,\eta}^{0}\left(u\right)$
for all sufficiently large $k$. Moreover, for every $\xi\in\mathbb{R}^{d}$,
\[
\widehat{f_{k}}\left(\xi\right)=\widehat{\varphi_{\nu}}\left(\xi\right)\sum_{j=1}^{m}\pi_{j,k}e^{-i\left\langle \mu_{j,k},\xi\right\rangle }\to\widehat{\varphi_{\nu}}\left(\xi\right)\sum_{j=1}^{m}\pi_{j}e^{-i\left\langle \mu_{j},\xi\right\rangle }=\widehat{f}\left(\xi\right),
\]
and therefore $\widehat{h_{k}}\left(\xi\right)\to\widehat{h}\left(\xi\right)$
pointwise. We will see below that all such Fourier transforms satisfy
a common exponential bound of the form 
\[
\left|\widehat{h}\left(\xi\right)\right|\le C_{d,\varphi}\exp\left[-c\frac{\left\Vert \xi\right\Vert _{2}^{\beta}}{\nu^{\beta}}\right]\qquad\text{for all }\xi\in\mathbb{R}^{d}
\]
with a finite constant $C_{d,\varphi}>0$ independent of $h$. Consequently,
for every $R>0$, dominated convergence on the Fourier inversion formula
gives 
\[
(h_{k})_{\le R}\left(x\right)\to h_{\le R}\left(x\right)\qquad\text{and}\qquad(h_{k})_{>R}\left(x\right)\to h_{>R}\left(x\right)
\]
for every $x\in\mathbb{R}^{d}$, where the low- and high-frequency
pieces are defined below. Moreover, for each fixed $R>0$, the same
bound yields finite constants 
\[
B_{\le R}=\left(2\pi\right)^{-d/2}C_{d,\varphi}\int_{\left\{ \left\Vert \xi\right\Vert _{2}\le R\right\} }\exp\left[-c\frac{\left\Vert \xi\right\Vert _{2}^{\beta}}{\nu^{\beta}}\right]\mathrm{d}\xi
\]
and 
\[
B_{>R}=\left(2\pi\right)^{-d/2}C_{d,\varphi}\int_{\left\{ \left\Vert \xi\right\Vert _{2}>R\right\} }\exp\left[-c\frac{\left\Vert \xi\right\Vert _{2}^{\beta}}{\nu^{\beta}}\right]\mathrm{d}\xi
\]
such that, for every $x\in\mathbb{R}^{d}$ and every $k$, 
\[
\left|(h_{k})_{\le R}\left(x\right)\right|,\left|h_{\le R}\left(x\right)\right|\le B_{\le R}\qquad\text{and}\qquad\left|(h_{k})_{>R}\left(x\right)\right|,\left|h_{>R}\left(x\right)\right|\le B_{>R}\text{.}
\]
Since $|\zeta_{n}|$ is a finite measure, dominated convergence with
respect to $|\zeta_{n}|$ yields 
\[
\zeta_{n}\left((h_{k})_{\le R}\right)\to\zeta_{n}\left(h_{\le R}\right)\qquad\text{and}\qquad\zeta_{n}\left((h_{k})_{>R}\right)\to\zeta_{n}\left(h_{>R}\right)
\]
almost surely. Hence, for every $\eta>0$ and every $R>0$, 
\[
\sup_{h\in{\cal H}_{\delta}\left(u\right)}\left|\zeta_{n}\left(h_{\le R}\right)\right|\le\sup_{g\in{\cal H}_{\delta,\eta}^{0}\left(u\right)}\left|\zeta_{n}\left(g_{\le R}\right)\right|
\]
and similarly for the high-frequency part. Indeed, for each fixed
$h\in{\cal H}_{\delta}\left(u\right)$ the sequence constructed above
satisfies $h_{k}\in{\cal H}_{\delta,\eta}^{0}\left(u\right)$ for all
sufficiently large $k$. Thus the corresponding values
$\left|\zeta_{n}\left((h_{k})_{\le R}\right)\right|$ are bounded by the
right-hand supremum in the preceding display, and their limit
$\left|\zeta_{n}\left(h_{\le R}\right)\right|$ is bounded by the same
quantity. Taking the supremum over $h\in{\cal H}_{\delta}\left(u\right)$
gives the displayed comparison; the argument for $h_{>R}$ is identical.
Since the right-hand side of the desired inequality is continuous in
$\delta$, it therefore suffices to prove the bound with
${\cal H}_{\delta}\left(u\right)$ replaced by the countable class
${\cal H}_{\delta,\eta}^{0}\left(u\right)$ and $\delta$ replaced by
$\delta+\eta$, and then let $\eta\downarrow0$.
In what follows, fix $\eta>0$, write 
\[
\delta_{\eta}=\delta+\eta,
\]
and suppress the superscript $0$ from the local class. All suprema
below are then measurable.

Let $c,C,\beta>0$ be the constants from Assumption~\ref{ass:ExponentialFourierKernel}.
Take any $h\in{\cal H}_{\delta_{\eta}}\left(u\right)$. Since both
$f$ and $u$ belong to $\mathrm{co}\left({\cal P}_{\varphi,\nu}\right)$,
there exist probability measures $G$ and $H$ on $\mathbb{R}^{d}$
such that 
\[
f=\varphi_{\nu}*G\text{ \qquad and \qquad}u=\varphi_{\nu}*H\text{.}
\]
Hence 
\[
h=\varphi_{\nu}*\Delta\text{,}\qquad\Delta=G-H\text{,}
\]
where $\Delta$ is a finite signed measure with total variation norm
\[
\left\Vert \Delta\right\Vert _{\mathrm{TV}}\le\left\Vert G\right\Vert _{\mathrm{TV}}+\left\Vert H\right\Vert _{\mathrm{TV}}=2\text{.}
\]
Since $\widehat{\varphi_{\nu}}\left(\xi\right)=\widehat{\varphi}\left(\xi/\nu\right)$,
Assumption~\ref{ass:ExponentialFourierKernel} gives 
\[
\left|\widehat{\varphi_{\nu}}\left(\xi\right)\right|\le C\exp\left[-c\frac{\left\Vert \xi\right\Vert _{2}^{\beta}}{\nu^{\beta}}\right]\qquad\text{for all }\xi\in\mathbb{R}^{d}\text{.}
\]
Moreover, 
\[
\left|\widehat{\Delta}\left(\xi\right)\right|\le\left(2\pi\right)^{-d/2}\left\Vert \Delta\right\Vert _{\mathrm{TV}}\le2\left(2\pi\right)^{-d/2}\text{,}
\]
so, after increasing the constant if necessary, there exists $C_{d,\varphi}>0$
depending only on $d$, $\left\Vert \varphi\right\Vert _{\infty}$,
and the constants in Assumption~\ref{ass:ExponentialFourierKernel} such that 
\[
\left|\widehat{h}\left(\xi\right)\right|=\left|\left(2\pi\right)^{d/2}\widehat{\varphi_{\nu}}\left(\xi\right)\widehat{\Delta}\left(\xi\right)\right|\le C_{d,\varphi}\exp\left[-c\frac{\left\Vert \xi\right\Vert _{2}^{\beta}}{\nu^{\beta}}\right]\qquad\text{for all }\xi\in\mathbb{R}^{d}\text{.}
\]

Fix $R\ge\nu$ and define $h_{\le R}$ and $h_{>R}$ through their
Fourier transforms, 
\[
\widehat{h_{\le R}}=\widehat{h}\mathbf{1}_{\left\{ \left\Vert \xi\right\Vert _{2}\le R\right\} }\qquad\text{and}\qquad\widehat{h_{>R}}=\widehat{h}\mathbf{1}_{\left\{ \left\Vert \xi\right\Vert _{2}>R\right\} }\text{.}
\]
The preceding exponential bound implies that $\widehat{h_{\le R}}\in{\cal L}_{1}\cap{\cal L}_{2}$
and $\widehat{h_{>R}}\in{\cal L}_{1}$, so Fourier inversion defines
bounded continuous functions $h_{\le R}$ and $h_{>R}$ with $h=h_{\le R}+h_{>R}$.

Let $\zeta_{n}=P_{n}-P$ and define its Fourier transform by 
\[
\widehat{\zeta}_{n}\left(\xi\right)=\left(2\pi\right)^{-d/2}\int_{\mathbb{R}^{d}}e^{-i\left\langle x,\xi\right\rangle }\zeta_{n}\left(\mathrm{d}x\right)\text{.}
\]
Since $\left\Vert \zeta_{n}\right\Vert _{\mathrm{TV}}\le2$, for either
choice $\sharp\in\left\{ \le R,>R\right\} $ we have 
\[
\int_{\mathbb{R}^{d}}\left|\widehat{h_{\sharp}}\left(\xi\right)\right|\left|\widehat{\zeta}_{n}\left(\xi\right)\right|\mathrm{d}\xi\le\left(2\pi\right)^{-d/2}\left\Vert \zeta_{n}\right\Vert _{\mathrm{TV}}\int_{\mathbb{R}^{d}}\left|\widehat{h_{\sharp}}\left(\xi\right)\right|\mathrm{d}\xi<\infty\text{.}
\]
Therefore Fubini's theorem yields the pairing identity 
\[
\zeta_{n}\left(h_{\sharp}\right)=\int_{\mathbb{R}^{d}}\widehat{h_{\sharp}}\left(\xi\right)\overline{\widehat{\zeta}_{n}\left(\xi\right)}\mathrm{d}\xi\text{.}
\]
In particular, 
\begin{align*}
\left|\zeta_{n}\left(h_{>R}\right)\right| & =\left|\int_{\left\{ \left\Vert \xi\right\Vert _{2}>R\right\} }\widehat{h}\left(\xi\right)\overline{\widehat{\zeta}_{n}\left(\xi\right)}\mathrm{d}\xi\right|\\
 & \le C_{d,\varphi}\int_{\left\{ \left\Vert \xi\right\Vert _{2}>R\right\} }\exp\left[-c\frac{\left\Vert \xi\right\Vert _{2}^{\beta}}{\nu^{\beta}}\right]\left|\widehat{\zeta}_{n}\left(\xi\right)\right|\mathrm{d}\xi\text{.}
\end{align*}
For each fixed $\xi\in\mathbb{R}^{d}$, 
\[
\widehat{\zeta}_{n}\left(\xi\right)=\frac{\left(2\pi\right)^{-d/2}}{n}\sum_{i=1}^{n}\left[e^{-i\left\langle X_{i},\xi\right\rangle }-\mathrm{E}e^{-i\left\langle X,\xi\right\rangle }\right]\text{,}
\]
so $\mathrm{E}\left|\widehat{\zeta}_{n}\left(\xi\right)\right|^{2}\le C_{d}/n$
and therefore $\sqrt{n}\,\mathrm{E}\left|\widehat{\zeta}_{n}\left(\xi\right)\right|\le C_{d}$
by the Cauchy--Schwarz inequality. Hence 
\[
\sqrt{n}\,\mathrm{E}\sup_{h\in{\cal H}_{\delta_{\eta}}\left(u\right)}\left|\left(P_{n}-P\right)h_{>R}\right|\le C_{d,\varphi}\int_{\left\{ \left\Vert \xi\right\Vert _{2}>R\right\} }\exp\left[-c\frac{\left\Vert \xi\right\Vert _{2}^{\beta}}{\nu^{\beta}}\right]\mathrm{d}\xi\text{.}
\]
Using polar coordinates and the elementary bound $t^{d-1}e^{-ct^{\beta}}\le C_{d,\varphi}e^{-ct^{\beta}/2}$
for $t\ge1$, we obtain 
\[
\int_{\left\{ \left\Vert \xi\right\Vert _{2}>R\right\} }\exp\left[-c\frac{\left\Vert \xi\right\Vert _{2}^{\beta}}{\nu^{\beta}}\right]\mathrm{d}\xi\le C_{d,\varphi}\nu^{d}\exp\left[-\frac{cR^{\beta}}{2\nu^{\beta}}\right]\text{,}\qquad R\ge\nu\text{.}
\]
Choose 
\[
R_{\delta_{\eta}}=\nu\max\left\{ 1,\left(\frac{2}{c}\log_{+}\left(\frac{4C_{d,\varphi}\nu^{d/2}}{\delta_{\eta}}\right)\right)^{1/\beta}\right\} \text{.}
\]
Then $R_{\delta_{\eta}}\ge\nu$ and, since $M_{\nu}\ge\left\Vert \varphi_{\nu}\right\Vert _{\infty}=\nu^{d}\left\Vert \varphi\right\Vert _{\infty}$,
\[
\sqrt{n}\,\mathrm{E}\sup_{h\in{\cal H}_{\delta_{\eta}}\left(u\right)}\left|\left(P_{n}-P\right)h_{>R_{\delta_{\eta}}}\right|\le C_{d,\varphi}\nu^{d}\exp\left[-\frac{cR_{\delta_{\eta}}^{\beta}}{2\nu^{\beta}}\right]\le C_{d,\varphi}\sqrt{M_{\nu}}\,\delta_{\eta}\text{.}
\]

For the low-frequency part, the same pairing identity and the Cauchy--Schwarz inequality
give 
\begin{align*}
\left|\zeta_{n}\left(h_{\le R}\right)\right| & \le\left\Vert \widehat{h}\mathbf{1}_{\left\{ \left\Vert \xi\right\Vert _{2}\le R\right\} }\right\Vert _{2}\left\Vert \widehat{\zeta}_{n}\right\Vert _{{\cal L}_{2}\left(B\left(0,R\right)\right)}\\
 & \le\left\Vert h\right\Vert _{2}\left\Vert \widehat{\zeta}_{n}\right\Vert _{{\cal L}_{2}\left(B\left(0,R\right)\right)}\\
 & \le\delta_{\eta}\left\Vert \widehat{\zeta}_{n}\right\Vert _{{\cal L}_{2}\left(B\left(0,R\right)\right)}\text{,}
\end{align*}
where we used Plancherel's theorem in the second step. Therefore 
\[
\sup_{h\in{\cal H}_{\delta_{\eta}}\left(u\right)}\left|\zeta_{n}\left(h_{\le R}\right)\right|\le\delta_{\eta}\left\Vert \widehat{\zeta}_{n}\right\Vert _{{\cal L}_{2}\left(B\left(0,R\right)\right)}\text{.}
\]
Now, for every fixed $\xi\in\mathbb{R}^{d}$, we still have $\mathrm{E}\left|\widehat{\zeta}_{n}\left(\xi\right)\right|^{2}\le C_{d}/n$,
and hence 
\begin{align*}
\mathrm{E}\left\Vert \widehat{\zeta}_{n}\right\Vert _{{\cal L}_{2}\left(B\left(0,R\right)\right)}^{2} & =\int_{B\left(0,R\right)}\mathrm{E}\left|\widehat{\zeta}_{n}\left(\xi\right)\right|^{2}\mathrm{d}\xi\\
 & \le\frac{C_{d}}{n}\left|B\left(0,R\right)\right|\\
 & \le\frac{C_{d}R^{d}}{n}\text{.}
\end{align*}
By Jensen's inequality, 
\[
\sqrt{n}\,\mathrm{E}\left\Vert \widehat{\zeta}_{n}\right\Vert _{{\cal L}_{2}\left(B\left(0,R\right)\right)}\le C_{d}R^{d/2}\text{.}
\]
Applying this with $R=R_{\delta_{\eta}}$ yields 
\[
\sqrt{n}\,\mathrm{E}\sup_{h\in{\cal H}_{\delta_{\eta}}\left(u\right)}\left|\left(P_{n}-P\right)h_{\le R_{\delta_{\eta}}}\right|\le C_{d}\delta_{\eta}R_{\delta_{\eta}}^{d/2}\text{.}
\]

Combining the low- and high-frequency bounds gives 
\[
\sqrt{n}\,\mathrm{E}\sup_{h\in{\cal H}_{\delta_{\eta}}\left(u\right)}\left|\left(P_{n}-P\right)h\right|\le C_{d,\varphi}\delta_{\eta}\left[\sqrt{M_{\nu}}+R_{\delta_{\eta}}^{d/2}\right]\text{.}
\]
Since $M_{\nu}\ge\left\Vert \varphi_{\nu}\right\Vert _{\infty}=\nu^{d}\left\Vert \varphi\right\Vert _{\infty}$,
we have $\nu^{d/2}\le C_{d,\varphi}\sqrt{M_{\nu}}$. By the definition
of $R_{\delta_{\eta}}$, 
\[
R_{\delta_{\eta}}^{d/2}\le C_{d,\varphi}\nu^{d/2}\left[1+\left(\log_{+}\frac{C_{d,\varphi}\nu^{d}}{\delta_{\eta}^{2}}\right)^{d/\left(2\beta\right)}\right]\le C_{d,\varphi}\sqrt{M_{\nu}}\left[1+\left(\log_{+}\frac{C_{d,\varphi}M_{\nu}}{\delta_{\eta}^{2}}\right)^{d/\left(2\beta\right)}\right]\text{.}
\]
Substituting this into the previous display yields 
\[
\sqrt{n}\,\mathrm{E}\sup_{h\in{\cal H}_{\delta_{\eta}}\left(u\right)}\left|\left(P_{n}-P\right)h\right|\le C_{d,\varphi}\sqrt{M_{\nu}}\,\delta_{\eta}\left[1+\left(1+\log_{+}\frac{C_{d,\varphi}M_{\nu}}{\delta_{\eta}^{2}}\right)^{d/\left(2\beta\right)}\right]\text{.}
\]
Letting $\eta\downarrow0$ proves the general stronger estimate.

Now assume, in addition, that $\operatorname{supp}\left(\widehat{\varphi}\right)\subseteq B\left(0,R_{0}\right)$
for some $R_{0}>0$. Then 
\[
\widehat{h}\left(\xi\right)=0\qquad\text{whenever }\left\Vert \xi\right\Vert _{2}>\nu R_{0}\text{,}
\]
because $\widehat{h}\left(\xi\right)=\left(2\pi\right)^{d/2}\widehat{\varphi}\left(\xi/\nu\right)\widehat{\Delta}\left(\xi\right)$.
Hence, with 
\[
R_{\mathrm{bl}}=\nu\max\left\{ 1,R_{0}\right\} \text{,}
\]
we have $h_{>R_{\mathrm{bl}}}=0$ for every $h\in{\cal H}_{\delta_{\eta}}\left(u\right)$.
The low-frequency estimate obtained above therefore gives 
\[
\sqrt{n}\,\mathrm{E}\sup_{h\in{\cal H}_{\delta_{\eta}}\left(u\right)}\left|\left(P_{n}-P\right)h\right|\le C_{d}\delta_{\eta}R_{\mathrm{bl}}^{d/2}\text{.}
\]
Since $M_{\nu}\ge\left\Vert \varphi_{\nu}\right\Vert _{\infty}=\nu^{d}\left\Vert \varphi\right\Vert _{\infty}$,
we obtain 
\[
R_{\mathrm{bl}}^{d/2}\le\max\left\{ 1,R_{0}\right\} ^{d/2}\nu^{d/2}\le\frac{\max\left\{ 1,R_{0}\right\} ^{d/2}}{\left\Vert \varphi\right\Vert _{\infty}^{1/2}}\sqrt{M_{\nu}}\text{.}
\]
Therefore 
\[
\sqrt{n}\,\mathrm{E}\sup_{h\in{\cal H}_{\delta_{\eta}}\left(u\right)}\left|\left(P_{n}-P\right)h\right|\le C_{d,\varphi,R_{0}}^{\mathrm{bl}}\sqrt{M_{\nu}}\,\delta_{\eta}
\]
for some finite constant $C_{d,\varphi,R_{0}}^{\mathrm{bl}}>0$ depending
only on $d$, $R_{0}$, and $\left\Vert \varphi\right\Vert _{\infty}$.
Letting $\eta\downarrow0$ proves the result, with the displayed
uncountable suprema interpreted through outer expectation. The countable
classes used in the reduction give ordinary measurable suprema at the
level needed in the proof of Theorem~\ref{thm:MassartNedelecLeastSquares}. 
\end{proof}
\begin{proof}[Proof of Theorem~\ref{thm:SlowLeastSquaresBenchmark}]
Since $f_{0}\in{\cal W}^{s,2}$, condition (\ref{eq:smoothness-assumption})
holds with $\alpha=s$ and some finite constant $K_{2}>0$; namely,
by Lemmas~\ref{lem:smoothness-W1p} and~\ref{lem:smoothness-Wsp},
we may take 
\[
K_{2}=\begin{cases}
C_{d,s,2}\left\Vert f_{0}\right\Vert _{s,2} & \text{if }0<s<1\text{,}\\
\left[\int_{\mathbb{R}^{d}}\left\Vert \nabla f_{0}\left(x\right)\right\Vert _{2}^{2}\mathrm{d}x\right]^{1/2} & \text{if }s=1\text{.}
\end{cases}
\]
For $m,n\in\mathbb{N}$, $\nu>0$, and $\epsilon>0$, let $\hat{f}_{m,n,\epsilon,\nu}\in\mathrm{co}_{m}\left({\cal P}_{\varphi,\nu}\right)$
be a measurable $\epsilon$-minimizer satisfying (\ref{eq:sample-eps-minimizer-def}).
By Lemma~\ref{lem:LeastSquaresBasicIneq}, for every $f_{m,\nu}\in\mathrm{co}_{m}\left({\cal P}_{\varphi,\nu}\right)$,
\[
\left\Vert \hat{f}_{m,n,\epsilon,\nu}-f_{0}\right\Vert _{2}^{2}\le\left\Vert f_{m,\nu}-f_{0}\right\Vert _{2}^{2}+2\left(P_{n}-P\right)\left(\hat{f}_{m,n,\epsilon,\nu}-f_{m,\nu}\right)+\epsilon\text{.}
\]
Taking expectations and using 
\[
\left|\left(P_{n}-P\right)\left(\hat{f}_{m,n,\epsilon,\nu}-f_{m,\nu}\right)\right|\le2\sup_{f\in\mathrm{co}_{m}\left({\cal P}_{\varphi,\nu}\right)}\left|P_{n}f-Pf\right|
\]
gives 
\[
\mathrm{E}\left\Vert \hat{f}_{m,n,\epsilon,\nu}-f_{0}\right\Vert _{2}^{2}\le\left\Vert f_{m,\nu}-f_{0}\right\Vert _{2}^{2}+4\mathrm{E}\sup_{f\in\mathrm{co}_{m}\left({\cal P}_{\varphi,\nu}\right)}\left|P_{n}f-Pf\right|+\epsilon\text{.}
\]
Now choose $f_{m,\nu}$ as in Lemma~\ref{lem:FixedScaleCompression}
with $p=2$. Since $\left(a+b\right)^{2}\le2a^{2}+2b^{2}$, Lemmas~\ref{lem:FixedScaleCompression}
and~\ref{lem:SmoothingBias} yield 
\[
\left\Vert f_{m,\nu}-f_{0}\right\Vert _{2}^{2}\lesssim\frac{\nu^{d}}{m}+\nu^{-2s}\text{.}
\]
Further, if 
\[
f=\sum_{j=1}^{m}\pi_{j}g_{j}\in\mathrm{co}_{m}\left({\cal P}_{\varphi,\nu}\right)
\]
with $g_{j}\in{\cal P}_{\varphi,\nu}$ and $\pi_{j}\ge0$ summing
to one, then 
\[
\left|P_{n}f-Pf\right|\le\sum_{j=1}^{m}\pi_{j}\left|P_{n}g_{j}-Pg_{j}\right|\le\sup_{g\in{\cal P}_{\varphi,\nu}}\left|P_{n}g-Pg\right|\text{.}
\]
Taking the supremum over $f\in\mathrm{co}_{m}\left({\cal P}_{\varphi,\nu}\right)$
and observing that ${\cal P}_{\varphi,\nu}\subset\mathrm{co}_{m}\left({\cal P}_{\varphi,\nu}\right)$,
we obtain 
\[
\sup_{f\in\mathrm{co}_{m}\left({\cal P}_{\varphi,\nu}\right)}\left|P_{n}f-Pf\right|=\sup_{f\in{\cal P}_{\varphi,\nu}}\left|P_{n}f-Pf\right|\text{.}
\]
Since $\varphi\in\ell^{\infty}\left(\mathbb{R}^{d}\right)$, for each
$f\in{\cal P}_{\varphi,\nu}$ and $x\in\mathbb{R}^{d}$, 
\[
\left|f\left(x\right)\right|\le\nu^{d}\left\Vert \varphi\right\Vert _{\infty}\text{,}
\]
so $\bar{F}_{\nu}=\nu^{d}\left\Vert \varphi\right\Vert _{\infty}$
is an envelope for ${\cal P}_{\varphi,\nu}$. Since ${\cal P}_{\varphi,\nu}$
is pointwise measurable, ${\cal P}_{\varphi,\nu}\subset{\cal P}_{\varphi}$,
and VC dimension is monotone under inclusion, Lemma~\ref{lem:VCBoundMaximalIneq}
implies that 
\[
\mathrm{E}\sup_{f\in\mathrm{co}_{m}\left({\cal P}_{\varphi,\nu}\right)}\left|P_{n}f-Pf\right|\lesssim\nu^{d}n^{-1/2}\text{.}
\]
Consequently, 
\[
\mathrm{E}\left\Vert \hat{f}_{m,n,\epsilon,\nu}-f_{0}\right\Vert _{2}^{2}\lesssim\frac{\nu^{d}}{m}+\nu^{-2s}+\nu^{d}n^{-1/2}+\epsilon\text{.}
\]
Now fix $M>0$ and set $m_{n}=\left\lceil Mn^{1/2}\right\rceil $.
Then $m_{n}^{-1}\asymp n^{-1/2}$. Choose $N>0$ and set $\nu_{n}=Nn^{1/\left(2\left(2s+d\right)\right)}$;
then $\nu_{n}^{d}m_{n}^{-1}\asymp\nu_{n}^{d}n^{-1/2}\asymp\nu_{n}^{-2s}\asymp n^{-\frac{s}{2s+d}}$.
Finally, choose $\epsilon_{n}\lesssim n^{-s/\left(2s+d\right)}$ and
let $\hat{f}_{n}=\hat{f}_{m_{n},n,\epsilon_{n},\nu_{n}}$. The preceding
bound then gives 
\[
\mathrm{E}\left\Vert \hat{f}_{n}-f_{0}\right\Vert _{2}^{2}\lesssim n^{-\frac{s}{2s+d}}\text{,}
\]
as required. 
\end{proof}
\begin{proof}[Proof of Theorem~\ref{thm:FastFixedScaleLeastSquares}]
Fix $n\ge3$ and write 
\[
{\cal F}_{n}={\cal G}_{m_{n},\nu}\text{ \qquad and \qquad}A_{n}=K_{d,\varphi}\left(1+\left(\log n\right)^{d/\left(2\beta\right)}\right)\text{,}
\]
where $K_{d,\varphi}$ denotes the constant from Lemma~\ref{lem:GaussianLocalizedOscillation}.
Choose a sequence $\left(u_{k}\right)_{k}\subset\mathrm{co}\left({\cal P}_{\varphi,\nu}\right)$
with $\left\Vert u_{k}-f_{0}\right\Vert _{2}\to0$. Passing to a subsequence,
we may assume that $u_{k}\left(x\right)\to f_{0}\left(x\right)$ for
$\lambda$-a.e. $x\in\mathbb{R}^{d}$. Since each $u_{k}$ is a probability
density and satisfies 
\[
0\le u_{k}\left(x\right)\le\left\Vert \varphi_{\nu}\right\Vert _{\infty}\qquad\text{for every }x\in\mathbb{R}^{d}\text{,}
\]
it follows that 
\[
0\le f_{0}\left(x\right)\le\left\Vert \varphi_{\nu}\right\Vert _{\infty}\qquad\text{for \ensuremath{\lambda}-a.e. }x\in\mathbb{R}^{d}\text{.}
\]
Hence $f_{0}\in{\cal L}_{\infty}$ and 
\[
M_{\nu}=\left\Vert f_{0}\right\Vert _{\infty}\vee\left\Vert \varphi_{\nu}\right\Vert _{\infty}=\left\Vert \varphi_{\nu}\right\Vert _{\infty}\asymp\nu^{d}\text{.}
\]
By Lemmas~\ref{lem:GaussianSieveSeparability} and~\ref{lem:GaussianLocalizedOscillation},
the hypotheses of Theorem~\ref{thm:MassartNedelecLeastSquares} hold
on ${\cal F}_{n}$ with $M=M_{\nu}$ and $A=A_{n}$. Consequently,
\[
\mathrm{E}\left\Vert \hat{f}_{n}-f_{0}\right\Vert _{2}^{2}\le2\inf_{f\in{\cal F}_{n}}\left\Vert f-f_{0}\right\Vert _{2}^{2}+2\epsilon_{n}+C\,M_{\nu}\frac{A_{n}^{2}}{n}
\]
for some finite universal constant $C>0$.

Next, fix any $g\in{\cal P}_{\varphi,\nu}$. Since translations preserve
${\cal L}_{2}$ norms, we have $\left\Vert g\right\Vert _{2}=\left\Vert \varphi_{\nu}\right\Vert _{2}$.
Also, because $u_{k}\to f_{0}$ in ${\cal L}_{2}$ and $\left\Vert u_{k}\right\Vert _{2}\le\left\Vert \varphi_{\nu}\right\Vert _{2}$
for every $k$, we obtain $\left\Vert f_{0}\right\Vert _{2}\le\left\Vert \varphi_{\nu}\right\Vert _{2}$.
Hence 
\[
\left\Vert g-f_{0}\right\Vert _{2}\le\left\Vert g\right\Vert _{2}+\left\Vert f_{0}\right\Vert _{2}\le2\left\Vert \varphi_{\nu}\right\Vert _{2}\text{.}
\]
Therefore Lemma~\ref{lem:convex-approximation} with $p=2$, ${\cal A}={\cal P}_{\varphi,\nu}$,
$B=2\left\Vert \varphi_{\nu}\right\Vert _{2}$, and, say, $\epsilon=\left\Vert \varphi_{\nu}\right\Vert _{2}$,
yields an element $g_{n}\in{\cal F}_{n}$ such that 
\[
\left\Vert g_{n}-f_{0}\right\Vert _{2}\lesssim\left\Vert \varphi_{\nu}\right\Vert _{2}m_{n}^{-1/2}\text{.}
\]
Therefore, 
\[
\inf_{f\in{\cal F}_{n}}\left\Vert f-f_{0}\right\Vert _{2}^{2}\le\left\Vert g_{n}-f_{0}\right\Vert _{2}^{2}\lesssim\frac{\nu^{d}}{m_{n}}\text{.}
\]
Since $A_{n}^{2}\lesssim\left(1+\log n\right)^{d/\beta}$ and $M_{\nu}\asymp\nu^{d}$,
this proves the first display.

Now assume, in addition, that $\operatorname{supp}\left(\widehat{\varphi}\right)\subseteq B\left(0,R_{0}\right)$
for some $R_{0}>0$. By Lemmas~\ref{lem:GaussianSieveSeparability} and~\ref{lem:GaussianLocalizedOscillation}, the hypotheses of
Theorem~\ref{thm:MassartNedelecLeastSquares} also hold on ${\cal F}_{n}$
with $M=M_{\nu}$ and 
\[
A=K_{d,\varphi,R_{0}}^{\mathrm{bl}}\text{.}
\]
Consequently, 
\[
\mathrm{E}\left\Vert \hat{f}_{n}-f_{0}\right\Vert _{2}^{2}\le2\inf_{f\in{\cal F}_{n}}\left\Vert f-f_{0}\right\Vert _{2}^{2}+2\epsilon_{n}+C\,M_{\nu}\frac{1}{n}
\]
for some finite universal constant $C>0$. Together with the same
approximation bound and $M_{\nu}\asymp\nu^{d}$, this proves the bandlimited
display.

For the fixed-scale model assertions, fix $f_{0}\in{\cal M}_{\varphi,\nu}$
and write 
\[
f_{0}\left(x\right)=\int_{\mathbb{R}^{d}}\varphi_{\nu}\left(x-\mu\right)G\left(\mathrm{d}\mu\right)
\]
for some probability measure $G$ on $\mathbb{R}^{d}$. By Lemma~\ref{lem:finite-supported-weak-density},
we may choose finitely supported probability measures $G_{k}$ on
$\mathbb{R}^{d}$ that converge weakly to $G$, and set 
\[
f_{k}\left(x\right)=\int_{\mathbb{R}^{d}}\varphi_{\nu}\left(x-\mu\right)G_{k}\left(\mathrm{d}\mu\right)\in\mathrm{co}\left({\cal P}_{\varphi,\nu}\right)\text{.}
\]
If $\widehat{G}_{k}$ and $\widehat{G}$ denote the Fourier transforms
of the measures $G_{k}$ and $G$, then 
\[
\widehat{f_{k}}\left(\xi\right)=\left(2\pi\right)^{d/2}\widehat{\varphi_{\nu}}\left(\xi\right)\widehat{G}_{k}\left(\xi\right)\text{ \qquad and \qquad}\widehat{f_{0}}\left(\xi\right)=\left(2\pi\right)^{d/2}\widehat{\varphi_{\nu}}\left(\xi\right)\widehat{G}\left(\xi\right)\text{.}
\]
Because weak convergence of finite measures implies pointwise convergence
of their Fourier transforms against the bounded continuous test function
$x\mapsto\exp\left(-i\left\langle \xi,x\right\rangle \right)$ (cf.
\citealt[Thm. 2.13]{vandervaart2000asymptotic}), we have $\widehat{G}_{k}\left(\xi\right)\to\widehat{G}\left(\xi\right)$
for every $\xi\in\mathbb{R}^{d}$. Moreover, 
\[
\left|\widehat{G}_{k}\left(\xi\right)-\widehat{G}\left(\xi\right)\right|\le2\left(2\pi\right)^{-d/2}\qquad\text{for all }\xi\in\mathbb{R}^{d}\text{ and all }k\text{.}
\]
Hence 
\[
\left|\widehat{f_{k}}\left(\xi\right)-\widehat{f_{0}}\left(\xi\right)\right|^{2}\le4\left|\widehat{\varphi_{\nu}}\left(\xi\right)\right|^{2}\qquad\text{for every }\xi\in\mathbb{R}^{d}\text{ and all }k\text{,}
\]
and the right-hand side is integrable because $\widehat{\varphi_{\nu}}\in{\cal L}_{2}$.
Dominated convergence and Plancherel's theorem therefore yield 
\[
\left\Vert f_{k}-f_{0}\right\Vert _{2}^{2}=\left\Vert \widehat{f_{k}}-\widehat{f_{0}}\right\Vert _{2}^{2}\to0\text{.}
\]
Thus $f_{0}\in\overline{\mathrm{co}\left({\cal P}_{\varphi,\nu}\right)}^{2}$.
Moreover, 
\[
0\le f_{0}\left(x\right)\le\left\Vert \varphi_{\nu}\right\Vert _{\infty}\qquad\text{for every }x\in\mathbb{R}^{d}\text{.}
\]
Applying the general bound above with the tuned choice of $m_{n}$
and $\epsilon_{n}$ proves the second display. If, in addition, $\operatorname{supp}\left(\widehat{\varphi}\right)\subseteq B\left(0,R_{0}\right)$,
then the bandlimited display with $m_{n}=\left\lceil Mn\right\rceil $
and $\epsilon_{n}\lesssim n^{-1}$ proves the final assertion. 
\end{proof}
\begin{proof}[Proof of Theorem~\ref{thm:GaussianCoreLowerBound}]
For $h\in{\cal L}_{2}$ and $\nu>0$, define the scaling operator
\[
\left(S_{\nu}h\right)\left(x\right)=\nu^{d}h\left(\nu x\right)\text{,}\qquad x\in\mathbb{R}^{d}\text{.}
\]
Then 
\[
\left\Vert S_{\nu}h\right\Vert _{2}^{2}=\nu^{d}\left\Vert h\right\Vert _{2}^{2}\text{,}
\]
and 
\[
\left(S_{\nu}^{-1}h\right)\left(x\right)=\nu^{-d}h\left(x/\nu\right)\text{,}\qquad x\in\mathbb{R}^{d}\text{,}
\]
with $S_{\nu}\left({\cal M}_{\gamma,1}\right)={\cal M}_{\gamma,\nu}$.
Fix any estimator $\hat{f}_{n}$ based on $n$ observations from ${\cal M}_{\gamma,\nu}$.
If $Y_{1},\dots,Y_{n}$ are i.i.d.\ with density $g\in{\cal M}_{\gamma,1}$,
then $X_{i}=Y_{i}/\nu$ has density $S_{\nu}g$ for each $i\in\left[n\right]$.
Therefore 
\[
\hat{g}_{n}\left(y_{1},\dots,y_{n}\right)=S_{\nu}^{-1}\left(\hat{f}_{n}\left(y_{1}/\nu,\dots,y_{n}/\nu\right)\right)
\]
is an estimator on ${\cal M}_{\gamma,1}$ satisfying 
\[
\mathrm{E}_{g}\left\Vert \hat{g}_{n}-g\right\Vert _{2}^{2}=\nu^{-d}\mathrm{E}_{S_{\nu}g}\left\Vert \hat{f}_{n}-S_{\nu}g\right\Vert _{2}^{2}\text{.}
\]
Taking first the supremum over $g\in{\cal M}_{\gamma,1}$ and then
the infimum over $\hat{f}_{n}$ gives 
\[
R_{n}\left({\cal M}_{\gamma,\nu}\right)\ge\nu^{d}R_{n}\left({\cal M}_{\gamma,1}\right)\text{.}
\]
By \citet[Thm.~2.1]{kim2022minimax}, 
\[
R_{n}\left({\cal M}_{\gamma,1}\right)\gtrsim\frac{\left(\log n\right)^{d/2}}{n}\text{,}
\]
whence 
\[
R_{n}\left({\cal M}_{\gamma,\nu}\right)\gtrsim\nu^{d}\frac{\left(\log n\right)^{d/2}}{n}\text{.}
\]

Next, let $\eta_{\nu}$ denote the image of $\eta$ under the map
$x\mapsto x/\nu$. Then 
\[
\widehat{\eta_{\nu}}\left(\xi\right)=\widehat{\eta}\left(\xi/\nu\right)
\]
in the sense of the Appendix, and so 
\[
\inf_{\xi\in\mathbb{R}^{d}}\left|\widehat{\eta_{\nu}}\left(\xi\right)\right|\ge a_{\eta}\text{.}
\]
Moreover, since $\varphi=\gamma*\eta$, we have 
\[
\varphi_{\nu}=\gamma_{\nu}*\eta_{\nu}\text{.}
\]
Fix any estimator $\hat{f}_{n}$ based on $n$ observations from ${\cal M}_{\varphi,\nu}$,
and fix $g\in{\cal M}_{\gamma,\nu}$. Let $X_{1},\dots,X_{n}$ be
i.i.d.\ with density $g$, let $Z_{1},\dots,Z_{n}$ be i.i.d.\ with
law $\eta_{\nu}$, independent of $X_{1},\dots,X_{n}$, and set $Y_{i}=X_{i}+Z_{i}$
for $i\in\left[n\right]$. If 
\[
g\left(x\right)=\int_{\mathbb{R}^{d}}\gamma_{\nu}\left(x-\mu\right)G\left(\mathrm{d}\mu\right)
\]
for some probability measure $G$ on $\mathbb{R}^{d}$, then 
\[
\left(g*\eta_{\nu}\right)\left(x\right)=\int_{\mathbb{R}^{d}}\varphi_{\nu}\left(x-\mu\right)G\left(\mathrm{d}\mu\right)\in{\cal M}_{\varphi,\nu}\text{.}
\]
Thus $Y_{1},\dots,Y_{n}$ are i.i.d.\ with density $g*\eta_{\nu}$.
Define $T_{\eta,\nu}^{-1}:{\cal L}_{2}\to{\cal L}_{2}$ to be the
bounded linear operator with Fourier multiplier 
\[
\widehat{T_{\eta,\nu}^{-1}h}\left(\xi\right)=\frac{\widehat{h}\left(\xi\right)}{\left(2\pi\right)^{d/2}\widehat{\eta_{\nu}}\left(\xi\right)}\text{,}\qquad h\in{\cal L}_{2}\text{.}
\]
Since $\left|\widehat{\eta_{\nu}}\left(\xi\right)\right|\ge a_{\eta}$
for all $\xi$, Plancherel's theorem yields 
\[
\left\Vert T_{\eta,\nu}^{-1}h\right\Vert _{2}\le\left(2\pi\right)^{-d/2}a_{\eta}^{-1}\left\Vert h\right\Vert _{2}\qquad\text{for all }h\in{\cal L}_{2}\text{.}
\]
Now define the randomized estimator 
\[
\tilde{g}_{n}=T_{\eta,\nu}^{-1}\left(\hat{f}_{n}\left(Y_{1},\dots,Y_{n}\right)\right)\text{.}
\]
Because $\widehat{g*\eta_{\nu}}=\left(2\pi\right)^{d/2}\widehat{g}\,\widehat{\eta_{\nu}}$,
we have $T_{\eta,\nu}^{-1}\left(g*\eta_{\nu}\right)=g$, and therefore
\[
\mathrm{E}_{g}\left\Vert \tilde{g}_{n}-g\right\Vert _{2}^{2}\le\left(2\pi\right)^{-d}a_{\eta}^{-2}\mathrm{E}_{g*\eta_{\nu}}\left\Vert \hat{f}_{n}-g*\eta_{\nu}\right\Vert _{2}^{2}\text{.}
\]
If the right-hand side is infinite, then there is nothing to prove.
Otherwise, $\tilde{g}_{n}$ is square-integrable as an ${\cal L}_{2}$-valued
random element. Since ${\cal L}_{2}$ is a separable Hilbert space,
its conditional expectation with respect to $\sigma\left(X_{1},\dots,X_{n}\right)$
is well-defined. Set 
\[
\bar{g}_{n}=\mathrm{E}\left[\tilde{g}_{n}\mid X_{1},\dots,X_{n}\right]\text{.}
\]
Then $\bar{g}_{n}$ depends only on the original sample $X_{1},\dots,X_{n}$,
and Jensen's inequality in the Hilbert space ${\cal L}_{2}$ gives
\[
\left\Vert \bar{g}_{n}-g\right\Vert _{2}^{2}\le\mathrm{E}\left[\left.\left\Vert \tilde{g}_{n}-g\right\Vert _{2}^{2}\right|X_{1},\dots,X_{n}\right]\text{.}
\]
Taking expectations therefore yields 
\[
\mathrm{E}_{g}\left\Vert \bar{g}_{n}-g\right\Vert _{2}^{2}\le\left(2\pi\right)^{-d}a_{\eta}^{-2}\mathrm{E}_{g*\eta_{\nu}}\left\Vert \hat{f}_{n}-g*\eta_{\nu}\right\Vert _{2}^{2}\text{.}
\]
Taking first the supremum over $g\in{\cal M}_{\gamma,\nu}$ and then
the infimum over $\hat{f}_{n}$ gives 
\[
R_{n}\left({\cal M}_{\gamma,\nu}\right)\le\left(2\pi\right)^{-d}a_{\eta}^{-2}R_{n}\left({\cal M}_{\varphi,\nu}\right)\text{.}
\]
Combining this with the Gaussian lower bound at scale $\nu$ proves
that 
\[
R_{n}\left({\cal M}_{\varphi,\nu}\right)\gtrsim a_{\eta}^{2}\nu^{d}\frac{\left(\log n\right)^{d/2}}{n}\text{.}
\]
Finally, since 
\[
\widehat{\varphi}\left(\xi\right)=\widehat{\gamma*\eta}\left(\xi\right)=\exp\left(-\left\Vert \xi\right\Vert _{2}^{2}/2\right)\widehat{\eta}\left(\xi\right)\text{,}
\]
we have 
\[
\left|\widehat{\varphi}\left(\xi\right)\right|\le\left(2\pi\right)^{-d/2}\exp\left(-\left\Vert \xi\right\Vert _{2}^{2}/2\right)\qquad\text{for all }\xi\in\mathbb{R}^{d}\text{.}
\]
Moreover, $\varphi\in{\cal P}$, $\left\Vert \varphi\right\Vert _{\infty}\le\left\Vert \gamma\right\Vert _{\infty}$,
and $\varphi\in{\cal L}_{2}$ by Plancherel's theorem. Thus $\varphi$
satisfies Assumption~\ref{ass:ExponentialFourierKernel} with $\beta=2$. 
\end{proof}

\subsection*{Proof of Theorem~\ref{thm:OddStudentLowerBound}}

We apply Assouad's lemma (Lemma~\ref{lem:Assouad}) to a tensorized
perturbation family built from Cauchy blocks.

For the remainder of this subsection, write 
\[
c\left(x\right)=\frac{1}{\pi\left(1+x^{2}\right)}\text{ \qquad and \qquad}c_{a}\left(x\right)=\frac{a}{\pi\left(a^{2}+x^{2}\right)}\text{,}\qquad x\in\mathbb{R}\text{ and }a>0\text{,}
\]
so $c=c_{1}$ and $c*c=c_{2}$. Let $\left(L_{k}\right)_{k\ge0}$
denote the Laguerre polynomials (see, e.g., \citealp[Sec. 18.5]{DLMF})
normalized by 
\[
\int_{0}^{\infty}e^{-s}L_{k}\left(s\right)L_{\ell}\left(s\right)\mathrm{d}s=\delta_{k\ell}\text{,}\qquad k,\ell\ge0\text{.}
\]
For $j\ge0$, define 
\[
w_{j}\left(s\right)=L_{2j}\left(s\right)-L_{2j+1}\left(s\right)\text{,}\qquad s\ge0\text{.}
\]
Since $L_{m}\left(0\right)=1$ for every $m\ge0$, we have $w_{j}\left(0\right)=0$.
Now define $\Gamma_{j}$ on $\mathbb{R}$ by 
\[
\widehat{\Gamma_{j}}\left(t\right)=e^{-2\left|t\right|}w_{j}\left(4\left|t\right|\right)\text{,}\qquad t\in\mathbb{R}\text{.}
\]

\begin{lem}
\label{lem:OddStudentCauchyBlocks}The family $\left\{ \Gamma_{j}:j\ge0\right\} $
satisfies the following properties. 
\begin{enumerate}
\item[(i)] $\left\{ \Gamma_{j}:j\ge0\right\} $ is orthonormal in ${\cal L}_{2}\left(\mathbb{R}\right)$. 
\item[(ii)] For every $j\ge0$, 
\[
\int_{\mathbb{R}}\frac{\Gamma_{j}\left(x\right)^{2}}{c_{2}\left(x\right)}\mathrm{d}x=4\pi\text{.}
\]
\end{enumerate}
\end{lem}

\begin{proof}
The Laguerre orthogonality relation gives 
\[
\int_{0}^{\infty}e^{-s}w_{j}\left(s\right)w_{\ell}\left(s\right)\mathrm{d}s=2\delta_{j\ell}\text{.}
\]
Hence, by Plancherel's theorem and the change of variables $s=4t$,
\begin{align*}
\int_{\mathbb{R}}\Gamma_{j}\left(x\right)\Gamma_{\ell}\left(x\right)\mathrm{d}x & =\int_{\mathbb{R}}\widehat{\Gamma_{j}}\left(t\right)\widehat{\Gamma_{\ell}}\left(t\right)\mathrm{d}t\\
 & =2\int_{0}^{\infty}e^{-4t}w_{j}\left(4t\right)w_{\ell}\left(4t\right)\mathrm{d}t\\
 & =\frac{1}{2}\int_{0}^{\infty}e^{-s}w_{j}\left(s\right)w_{\ell}\left(s\right)\mathrm{d}s\\
 & =\delta_{j\ell}\text{,}
\end{align*}
which proves (i).

For (ii), the Laplace transform identity 
\[
\int_{0}^{\infty}e^{-pt}L_{k}\left(qt\right)\mathrm{d}t=\frac{\left(p-q\right)^{k}}{p^{k+1}}\text{,}\qquad\Re p>0\text{,}
\]
with $p=2-ix$ and $q=4$ gives 
\[
\Gamma_{j}\left(x\right)=4\sqrt{\frac{2}{\pi}}\,\Re\left(\frac{\left(2+ix\right)^{2j}}{\left(2-ix\right)^{2j+2}}\right)\text{.}
\]
Writing $x=2\tan\theta$ with $\theta\in\left(-\pi/2,\pi/2\right)$,
we have $2\pm ix=2\sec\theta\,e^{\pm i\theta}$ and therefore 
\[
\Gamma_{j}\left(x\right)=\sqrt{\frac{2}{\pi}}\cos^{2}\theta\,\cos\left(\left(4j+2\right)\theta\right)\text{.}
\]
Since 
\[
c_{2}\left(x\right)=\frac{1}{2\pi}\cos^{2}\theta\qquad\text{and}\qquad c_{2}\left(x\right)\mathrm{d}x=\frac{\mathrm{d}\theta}{\pi}\text{,}
\]
it follows that 
\begin{align*}
\int_{\mathbb{R}}\frac{\Gamma_{j}\left(x\right)^{2}}{c_{2}\left(x\right)}\mathrm{d}x & =8\int_{-\pi/2}^{\pi/2}\cos^{2}\left(\left(4j+2\right)\theta\right)\mathrm{d}\theta\\
 & =4\pi\text{.}
\end{align*}
This proves (ii). 
\end{proof}
Fix $r\in\mathbb{N}\cup\left\{ 0\right\} $. Recall that 
\[
\tau_{r}\left(x\right)=\frac{\Gamma\left(r+1\right)}{\sqrt{\pi}\,\Gamma\left(r+\frac{1}{2}\right)}\left(1+x^{2}\right)^{-r-1}\text{,}\qquad x\in\mathbb{R}\text{.}
\]
Let $q_{r}$ be the degree-$r$ Bessel polynomial 
\[
q_{r}\left(u\right)=\sum_{k=0}^{r}\alpha_{r,k}u^{k}\text{,}\qquad\alpha_{r,k}=\frac{r!\left(2r-k\right)!2^{k}}{\left(2r\right)!\left(r-k\right)!k!}\text{,}\qquad u\in\mathbb{C}\text{.}
\]

\begin{lem}
\label{lem:OddStudentBesselFactor}The following hold. 
\begin{enumerate}
\item[(i)] The Fourier transform of $\tau_{r}$ is 
\[
\widehat{\tau_{r}}\left(t\right)=\left(2\pi\right)^{-1/2}e^{-\left|t\right|}q_{r}\left(\left|t\right|\right)\text{,}\qquad t\in\mathbb{R}\text{.}
\]
\item[(ii)] The polynomial $q_{r}$ has no zeros in the closed right half-plane.
Consequently, there exist a finite complex Borel measure $\mu_{r}$
on $\left[0,\infty\right)$ and a finite constant 
\[
V_{r}=\left\Vert \mu_{r}\right\Vert _{\mathrm{TV}}<\infty
\]
such that 
\[
\frac{1}{q_{r}\left(t\right)}=\int_{0}^{\infty}e^{-st}\mu_{r}\left(\mathrm{d}s\right)\text{,}\qquad t\ge0\text{.}
\]
\end{enumerate}
\end{lem}

\begin{proof}
For (i), \citet[(5)--(7)]{BergVignat2008} show that for the odd-degree
Student-$t$ density $f_{r+1/2}$ one has 
\[
\int_{\mathbb{R}}e^{ixt}f_{r+1/2}(x)\,\mathrm{d}x=e^{-|t|}q_{r}(|t|).
\]
Under the unitary Fourier convention of the manuscript, this yields
\[
\widehat{\tau_{r}}(t)=(2\pi)^{-1/2}e^{-|t|}q_{r}(|t|).
\]

For (ii), if $r=0$, then $q_{0}=1$. In this case we may take
\[
\mu_{0}=\delta_{0},
\]
so that $V_{0}=1$ and the Laplace representation is immediate. Assume
henceforth that $r\ge1$. \citet[p.~16]{BergVignat2008} state that
\[
\theta_{r}\left(u\right)=\frac{\left(2r\right)!}{r!\,2^{r}}q_{r}\left(u\right),
\]
where $\theta_{r}$ denotes the reverse Bessel polynomial. Thus $q_{r}$
differs from the reverse Bessel polynomial only by a positive multiplicative
constant. The standard representation of $\theta_{r}$ through the
modified Bessel function $K_{r+1/2}$ (see \citealp[Sec.~18.34]{DLMF})
gives 
\[
K_{r+1/2}\left(z\right)=\sqrt{\frac{\pi}{2}}\,z^{-r-1/2}e^{-z}\theta_{r}\left(z\right).
\]
Since the prefactor on the right-hand side never vanishes on $\mathbb{C}\setminus\left\{ 0\right\} $
and $q_{r}\left(0\right)=1$, the zeros of $q_{r}$ are exactly the
zeros of $K_{r+1/2}$. It is known that, for real order, every zero
of $K_{\nu}$ lies in the open left half-plane; see \citet[Sec.~15.7]{Watson1944Bessel}.
Hence every zero of $q_{r}$ lies in the open left half-plane, and
we may write the distinct zeros as $-\lambda_{\ell,r}$ with $\Re\lambda_{\ell,r}>0$.

Since $q_{r}(0)=1$, the rational function $1/q_{r}$ has no polynomial
part. Therefore a partial-fraction decomposition gives 
\[
\frac{1}{q_{r}\left(z\right)}=\sum_{\ell=1}^{J_{r}}\sum_{m=1}^{m_{\ell,r}}\frac{a_{\ell,m,r}}{\left(z+\lambda_{\ell,r}\right)^{m}}\text{,}
\]
where $m_{\ell,r}$ is the multiplicity of the zero $-\lambda_{\ell,r}$.
For $t\ge0$ and every $m\ge1$, 
\[
\left(t+\lambda\right)^{-m}=\frac{1}{\left(m-1\right)!}\int_{0}^{\infty}s^{m-1}e^{-s\left(t+\lambda\right)}\mathrm{d}s\qquad\text{whenever }\Re\lambda>0\text{.}
\]
Hence, for $t\ge0$, 
\[
\frac{1}{q_{r}\left(t\right)}=\int_{0}^{\infty}e^{-st}g_{r}(s)\,\mathrm{d}s\text{,}
\]
where 
\[
g_{r}(s)=\sum_{\ell=1}^{J_{r}}\sum_{m=1}^{m_{\ell,r}}a_{\ell,m,r}\frac{s^{m-1}}{\left(m-1\right)!}e^{-\lambda_{\ell,r}s}\text{,}\qquad s\ge0\text{.}
\]
Moreover, 
\[
\int_{0}^{\infty}|g_{r}(s)|\,\mathrm{d}s\le\sum_{\ell=1}^{J_{r}}\sum_{m=1}^{m_{\ell,r}}|a_{\ell,m,r}|\frac{1}{\left(m-1\right)!}\int_{0}^{\infty}s^{m-1}e^{-\Re\lambda_{\ell,r}s}\mathrm{d}s=\sum_{\ell=1}^{J_{r}}\sum_{m=1}^{m_{\ell,r}}\frac{|a_{\ell,m,r}|}{\left(\Re\lambda_{\ell,r}\right)^{m}}<\infty\text{.}
\]
Thus $\mu_{r}(\mathrm{d}s)=g_{r}(s)\,\mathrm{d}s$ is a finite complex
Borel measure on $\left[0,\infty\right)$, and the preceding display
gives the required Laplace representation. Its total variation equals
\[
V_{r}=\left\Vert \mu_{r}\right\Vert _{\mathrm{TV}}=\int_{0}^{\infty}|g_{r}(s)|\,\mathrm{d}s<\infty\text{.}
\]
\end{proof}
For $j\ge0$, define $v_{r,j}$ on $\mathbb{R}$ by 
\[
\widehat{v_{r,j}}\left(t\right)=e^{-\left|t\right|}\frac{w_{j}\left(4\left|t\right|\right)}{q_{r}\left(\left|t\right|\right)}\text{,}\qquad t\in\mathbb{R}\text{.}
\]
Since $w_{j}\left(0\right)=0$ and $q_{r}\left(0\right)=1$, we have
\[
\widehat{v_{r,j}}\left(0\right)=0\text{.}
\]
Also, by Lemma~\ref{lem:OddStudentBesselFactor}, 
\[
\left(2\pi\right)^{1/2}\widehat{\tau_{r}}\left(t\right)\widehat{v_{r,j}}\left(t\right)=e^{-2\left|t\right|}w_{j}\left(4\left|t\right|\right)=\widehat{\Gamma_{j}}\left(t\right)\text{.}
\]

\begin{lem}
\label{lem:OddStudentvrjBound}There exists a finite constant 
\[
U_{r}=4\sqrt{2\pi}\,V_{r}
\]
such that, for every $j\ge0$ and every $x\in\mathbb{R}$, 
\[
\left|v_{r,j}\left(x\right)\right|\le U_{r}9^{j}c\left(x\right)\text{.}
\]
\end{lem}

\begin{proof}
Because $\widehat{v_{r,j}}$ is even and real-valued, 
\[
v_{r,j}\left(x\right)=\sqrt{\frac{2}{\pi}}\int_{0}^{\infty}e^{-t}\frac{w_{j}\left(4t\right)}{q_{r}\left(t\right)}\cos\left(tx\right)\mathrm{d}t\text{.}
\]
Since $w_{j}$ is a polynomial and $\mu_{r}$ is finite, we have 
\[
\int_{0}^{\infty}\int_{0}^{\infty}e^{-(1+s)t}\left|w_{j}(4t)\right|\mathrm{d}t\,|\mu_{r}|(\mathrm{d}s)\le\left\Vert \mu_{r}\right\Vert _{\mathrm{TV}}\int_{0}^{\infty}e^{-t}\left|w_{j}(4t)\right|\mathrm{d}t<\infty\text{,}
\]
so Fubini's theorem applies to the Laplace representation from Lemma~\ref{lem:OddStudentBesselFactor}(ii).
Hence 
\[
v_{r,j}\left(x\right)=\sqrt{\frac{2}{\pi}}\,\Re\int_{0}^{\infty}\left[\int_{0}^{\infty}e^{-\left(1+s-ix\right)t}w_{j}\left(4t\right)\mathrm{d}t\right]\mu_{r}\left(\mathrm{d}s\right)\text{.}
\]
The Laplace transform identity for the Laguerre combination $w_{j}=L_{2j}-L_{2j+1}$
yields 
\[
\int_{0}^{\infty}e^{-pt}w_{j}\left(4t\right)\mathrm{d}t=4\frac{\left(p-4\right)^{2j}}{p^{2j+2}}\text{,}\qquad\Re p>0\text{.}
\]
Applying this with $p=1+s-ix$ gives 
\[
\left|\int_{0}^{\infty}e^{-\left(1+s-ix\right)t}w_{j}\left(4t\right)\mathrm{d}t\right|=4\frac{\left(\left(s-3\right)^{2}+x^{2}\right)^{j}}{\left(\left(s+1\right)^{2}+x^{2}\right)^{j+1}}\text{.}
\]
Now 
\[
\left(s-3\right)^{2}+x^{2}\le9\left(\left(s+1\right)^{2}+x^{2}\right)\text{ \qquad and \qquad}\frac{1}{\left(s+1\right)^{2}+x^{2}}\le\frac{1}{1+x^{2}}\text{,}
\]
so 
\[
\left|\int_{0}^{\infty}e^{-\left(1+s-ix\right)t}w_{j}\left(4t\right)\mathrm{d}t\right|\le\frac{4\,9^{j}}{1+x^{2}}\text{.}
\]
Taking absolute values and integrating against $\left\Vert \mu_{r}\right\Vert _{\mathrm{TV}}$
yields 
\[
\left|v_{r,j}\left(x\right)\right|\le4\sqrt{\frac{2}{\pi}}\,V_{r}\,\frac{9^{j}}{1+x^{2}}=4\sqrt{2\pi}\,V_{r}\,9^{j}c\left(x\right)\text{,}
\]
as claimed. 
\end{proof}
Since $c\in{\cal L}_{1}\left(\mathbb{R}\right)$, Lemma~\ref{lem:OddStudentvrjBound}
implies that $v_{r,j}\in{\cal L}_{1}\left(\mathbb{R}\right)$ for
every $j\ge0$. Therefore 
\[
\int_{\mathbb{R}}v_{r,j}\left(x\right)\mathrm{d}x=\left(2\pi\right)^{1/2}\widehat{v_{r,j}}\left(0\right)=0\text{.}
\]
Also, since $\tau_{r},v_{r,j}\in{\cal L}_{1}\left(\mathbb{R}\right)$,
the convolution $\tau_{r}*v_{r,j}$ is well defined and 
\[
\widehat{\tau_{r}*v_{r,j}}\left(t\right)=\left(2\pi\right)^{1/2}\widehat{\tau_{r}}\left(t\right)\widehat{v_{r,j}}\left(t\right)=\widehat{\Gamma_{j}}\left(t\right)\text{.}
\]
Because $\widehat{\Gamma_{j}}\in{\cal L}_{1}\left(\mathbb{R}\right)$,
Fourier inversion yields 
\[
\tau_{r}*v_{r,j}=\Gamma_{j}\text{.}
\]
For $d\in\mathbb{N}$ and $\mathbf{j}=\left(j_{1},\dots,j_{d}\right)\in\left(\mathbb{N}\cup\left\{ 0\right\} \right)^{d}$,
define 
\[
\Gamma_{\mathbf{j}}\left(x\right)=\prod_{k=1}^{d}\Gamma_{j_{k}}\left(x_{k}\right)\text{,}\qquad v_{r,\mathbf{j}}\left(x\right)=\prod_{k=1}^{d}v_{r,j_{k}}\left(x_{k}\right)\text{,}\qquad x=\left(x_{1},\dots,x_{d}\right)\in\mathbb{R}^{d}\text{,}
\]
and also write 
\[
c^{\otimes d}\left(x\right)=\prod_{k=1}^{d}c\left(x_{k}\right)\text{,}\qquad c_{2}^{\otimes d}\left(x\right)=\prod_{k=1}^{d}c_{2}\left(x_{k}\right)\text{,}\qquad\tau_{r,d}\left(x\right)=\prod_{k=1}^{d}\tau_{r}\left(x_{k}\right)\text{.}
\]

\begin{lem}
\label{lem:OddStudentTensorBlocks}For every $d\in\mathbb{N}$ and
every $\mathbf{j}\in\left(\mathbb{N}\cup\left\{ 0\right\} \right)^{d}$,
the following hold. 
\begin{enumerate}
\item[(i)] The family $\left\{ \Gamma_{\mathbf{j}}:\mathbf{j}\in\left(\mathbb{N}\cup\left\{ 0\right\} \right)^{d}\right\} $
is orthonormal in ${\cal L}_{2}\left(\mathbb{R}^{d}\right)$. 
\item[(ii)] One also has the weighted identity 
\[
\int_{\mathbb{R}^{d}}\frac{\Gamma_{\mathbf{j}}\left(x\right)^{2}}{c_{2}^{\otimes d}\left(x\right)}\mathrm{d}x=\left(4\pi\right)^{d}\text{.}
\]
\item[(iii)] $\tau_{r,d}*v_{r,\mathbf{j}}=\Gamma_{\mathbf{j}}$ and 
\[
\int_{\mathbb{R}^{d}}v_{r,\mathbf{j}}\left(x\right)\mathrm{d}x=0\text{.}
\]
\item[(iv)] For every $x\in\mathbb{R}^{d}$, 
\[
\left|v_{r,\mathbf{j}}\left(x\right)\right|\le U_{r}^{d}9^{\left\Vert \mathbf{j}\right\Vert _{1}}c^{\otimes d}\left(x\right)\text{.}
\]
\end{enumerate}
\end{lem}

\begin{proof}
Each claim follows by tensorization from Lemma~\ref{lem:OddStudentCauchyBlocks},
Lemma~\ref{lem:OddStudentvrjBound}, and the one-dimensional identities
displayed immediately above. Indeed, (i), (ii), and the integral identity
in (iii) follow from Fubini's theorem, while the convolution identity
in (iii) uses the product structure of $\tau_{r,d}$ and $v_{r,\mathbf{j}}$.
Finally, 
\[
\left|v_{r,\mathbf{j}}\left(x\right)\right|=\prod_{k=1}^{d}\left|v_{r,j_{k}}\left(x_{k}\right)\right|\le\prod_{k=1}^{d}\left(U_{r}9^{j_{k}}c\left(x_{k}\right)\right)=U_{r}^{d}9^{\left\Vert \mathbf{j}\right\Vert _{1}}c^{\otimes d}\left(x\right)\text{,}
\]
which proves (iv). 
\end{proof}
\begin{proof}[Proof of Theorem~\ref{thm:OddStudentLowerBound}]
For $h\in{\cal L}_{2}\left(\mathbb{R}^{d}\right)$ and $\nu>0$,
define 
\[
\left(S_{\nu}h\right)\left(x\right)=\nu^{d}h\left(\nu x\right)\text{,}\qquad x\in\mathbb{R}^{d}\text{.}
\]
Then 
\[
\left\Vert S_{\nu}h\right\Vert _{2}^{2}=\nu^{d}\left\Vert h\right\Vert _{2}^{2}\text{,}
\]
and the same rescaling argument as in the proof of Theorem~\ref{thm:GaussianCoreLowerBound}
shows that 
\[
R_{n}\left({\cal M}_{\tau_{r,d},\nu}\right)\ge\nu^{d}R_{n}\left({\cal M}_{\tau_{r,d},1}\right)\text{.}
\]
It therefore suffices to prove that there exists $\kappa_{r,d}>0$
such that 
\[
R_{n}\left({\cal M}_{\tau_{r,d},1}\right)\ge\kappa_{r,d}\frac{\left(\log n\right)^{d}}{n}
\]
for all sufficiently large $n$.

Fix 
\[
m=\left\lfloor \frac{\log n}{4d\log3}\right\rfloor \text{.}
\]
For all sufficiently large $n$, we have $m\ge1$ and $9^{md}\le\sqrt{n}$.
Let 
\[
{\cal J}_{m}=\left\{ 0,\dots,m-1\right\} ^{d}
\]
so that $\left|{\cal J}_{m}\right|=m^{d}$. Also let 
\[
p_{r}=\int_{-1}^{1}\tau_{r}\left(u\right)\mathrm{d}u>0
\]
and choose 
\[
\varepsilon_{r,d}=\min\left\{ U_{r}^{-d},\frac{1}{2}\left(\frac{p_{r}}{24\pi}\right)^{d/2}\right\} \text{.}
\]
Set 
\[
\epsilon=\varepsilon_{r,d}n^{-1/2}\text{.}
\]
For each $\alpha=\left(\alpha_{\mathbf{j}}\right)_{\mathbf{j}\in{\cal J}_{m}}\in\left\{ 0,1\right\} ^{{\cal J}_{m}}$,
define 
\[
g_{\alpha}\left(x\right)=c^{\otimes d}\left(x\right)+\epsilon\sum_{\mathbf{j}\in{\cal J}_{m}}\alpha_{\mathbf{j}}v_{r,\mathbf{j}}\left(x\right)\text{,}\qquad x\in\mathbb{R}^{d}\text{,}
\]
and let 
\[
f_{\alpha}=\tau_{r,d}*g_{\alpha}\text{.}
\]
Write 
\[
\Upsilon\left(\alpha,\beta\right)=\sum_{\mathbf{j}\in{\cal J}_{m}}\mathbf{1}\left\{ \alpha_{\mathbf{j}}\neq\beta_{\mathbf{j}}\right\} 
\]
for the Hamming distance on $\left\{ 0,1\right\} ^{{\cal J}_{m}}$.
By construction, each $g_{\alpha}$ integrates to one. Moreover, by
Lemma~\ref{lem:OddStudentTensorBlocks}, 
\begin{align*}
\epsilon\sum_{\mathbf{j}\in{\cal J}_{m}}\left|v_{r,\mathbf{j}}\left(x\right)\right| & \le\epsilon U_{r}^{d}\sum_{\mathbf{j}\in{\cal J}_{m}}9^{\left\Vert \mathbf{j}\right\Vert _{1}}c^{\otimes d}\left(x\right)\\
 & =\epsilon U_{r}^{d}\left(\sum_{j=0}^{m-1}9^{j}\right)^{d}c^{\otimes d}\left(x\right)\\
 & \le\epsilon U_{r}^{d}\left(\frac{9^{m}}{8}\right)^{d}c^{\otimes d}\left(x\right)\\
 & \le\frac{\varepsilon_{r,d}U_{r}^{d}}{8^{d}}c^{\otimes d}\left(x\right)\\
 & \le\frac{1}{8}c^{\otimes d}\left(x\right)\text{.}
\end{align*}
Hence 
\[
g_{\alpha}\left(x\right)\ge\frac{7}{8}c^{\otimes d}\left(x\right)\text{,}\qquad x\in\mathbb{R}^{d}\text{,}
\]
so each $g_{\alpha}$ is a probability density. Consequently, 
\[
f_{\alpha}\in{\cal M}_{\tau_{r,d},1}\qquad\text{for every }\alpha\in\left\{ 0,1\right\} ^{{\cal J}_{m}}\text{.}
\]

Using the identity $\tau_{r,d}*v_{r,\mathbf{j}}=\Gamma_{\mathbf{j}}$,
we may write 
\[
f_{\alpha}=f_{\ast}+\epsilon\sum_{\mathbf{j}\in{\cal J}_{m}}\alpha_{\mathbf{j}}\Gamma_{\mathbf{j}}\text{,}\qquad f_{\ast}=\tau_{r,d}*c^{\otimes d}\text{.}
\]
Since $g_{\alpha}\ge\frac{7}{8}c^{\otimes d}$, we also have 
\[
f_{\alpha}\ge\frac{7}{8}f_{\ast}\text{.}
\]

We next prove a pointwise lower bound on $f_{\ast}$. For $\left|u\right|\le1$,
\[
1+\left(x-u\right)^{2}\le1+x^{2}+2\left|x\right|+1\le3\left(1+x^{2}\right)\text{,}
\]
so 
\[
c\left(x-u\right)\ge\frac{1}{3}c\left(x\right)\ge\frac{1}{6}c_{2}\left(x\right)\text{,}\qquad x\in\mathbb{R}\text{.}
\]
Therefore, for every $x\in\mathbb{R}$, 
\[
\left(\tau_{r}*c\right)\left(x\right)\ge\int_{-1}^{1}\tau_{r}\left(u\right)c\left(x-u\right)\mathrm{d}u\ge\frac{p_{r}}{6}c_{2}\left(x\right)\text{.}
\]
Since 
\[
f_{\ast}\left(x\right)=\prod_{k=1}^{d}\left(\tau_{r}*c\right)\left(x_{k}\right)\text{,}\qquad x=\left(x_{1},\dots,x_{d}\right)\in\mathbb{R}^{d}\text{,}
\]
it follows that 
\[
f_{\ast}\left(x\right)\ge\left(\frac{p_{r}}{6}\right)^{d}c_{2}^{\otimes d}\left(x\right)\text{,}\qquad x\in\mathbb{R}^{d}\text{.}
\]

Let us now compute the ${\cal L}_{2}$ separation. By Lemma~\ref{lem:OddStudentTensorBlocks}(i),
for any $\alpha,\beta\in\left\{ 0,1\right\} ^{{\cal J}_{m}}$, 
\begin{align*}
\left\Vert f_{\alpha}-f_{\beta}\right\Vert _{2}^{2} & =\epsilon^{2}\left\Vert \sum_{\mathbf{j}\in{\cal J}_{m}}\left(\alpha_{\mathbf{j}}-\beta_{\mathbf{j}}\right)\Gamma_{\mathbf{j}}\right\Vert _{2}^{2}\\
 & =\epsilon^{2}\sum_{\mathbf{j}\in{\cal J}_{m}}\left(\alpha_{\mathbf{j}}-\beta_{\mathbf{j}}\right)^{2}\\
 & =\epsilon^{2}\Upsilon\left(\alpha,\beta\right)\text{.}
\end{align*}
Hence 
\[
\min_{\alpha\neq\beta}\frac{\left\Vert f_{\alpha}-f_{\beta}\right\Vert _{2}^{2}}{\Upsilon\left(\alpha,\beta\right)}=\epsilon^{2}\text{.}
\]

Now suppose that $\Upsilon\left(\alpha,\beta\right)=1$, and let $\mathbf{j}_{*}$
be the unique coordinate at which $\alpha$ and $\beta$ differ. Then
\[
f_{\alpha}-f_{\beta}=\pm\epsilon\Gamma_{\mathbf{j}_{*}}\text{.}
\]
Using the two preceding displays together with Lemma~\ref{lem:OddStudentTensorBlocks}(ii),
we obtain 
\begin{align*}
\chi^{2}\left(f_{\alpha}\|f_{\beta}\right) & =\int_{\mathbb{R}^{d}}\frac{\left(f_{\alpha}\left(x\right)-f_{\beta}\left(x\right)\right)^{2}}{f_{\beta}\left(x\right)}\mathrm{d}x\\
 & \le\frac{8}{7}\epsilon^{2}\int_{\mathbb{R}^{d}}\frac{\Gamma_{\mathbf{j}_{*}}\left(x\right)^{2}}{f_{\ast}\left(x\right)}\mathrm{d}x\\
 & \le\frac{8}{7}\epsilon^{2}\left(\frac{6}{p_{r}}\right)^{d}\int_{\mathbb{R}^{d}}\frac{\Gamma_{\mathbf{j}_{*}}\left(x\right)^{2}}{c_{2}^{\otimes d}\left(x\right)}\mathrm{d}x\\
 & \le\frac{8}{7}\epsilon^{2}\left(\frac{24\pi}{p_{r}}\right)^{d}\text{.}
\end{align*}
By the choice of $\varepsilon_{r,d}$, the right-hand side is at most
$1/\left(2n\right)$. Therefore 
\[
\max_{\Upsilon\left(\alpha,\beta\right)=1}\chi^{2}\left(f_{\alpha}\|f_{\beta}\right)\le\frac{1}{2n}\text{.}
\]

Assouad's lemma now yields 
\[
R_{n}\left({\cal M}_{\tau_{r,d},1}\right)\ge\frac{\left|{\cal J}_{m}\right|}{8}\,\epsilon^{2}\left(1-\sqrt{\frac{n}{2}\cdot\frac{1}{2n}}\right)=\frac{m^{d}\epsilon^{2}}{16}\text{.}
\]
Since $m^{d}\asymp\left(\log n\right)^{d}$ and $\epsilon^{2}=\varepsilon_{r,d}^{2}/n$,
this proves the required lower bound at scale $1$. Combining it with
the preceding rescaling display gives 
\[
R_{n}\left({\cal M}_{\tau_{r,d},\nu}\right)\ge\kappa_{r,d}\nu^{d}\frac{\left(\log n\right)^{d}}{n}
\]
for some constant $\kappa_{r,d}>0$ depending only on $r$ and $d$.

Finally, let $t_{2r+1}^{\otimes d}$ denote the product of $d$ copies
of the standard Student-$t$ density with $2r+1$ degrees of freedom.
Then 
\[
t_{2r+1}^{\otimes d}\left(x\right)=\left(2r+1\right)^{-d/2}\tau_{r,d}\left(\frac{x}{\sqrt{2r+1}}\right)=\left(\tau_{r,d}\right)_{1/\sqrt{2r+1}}\left(x\right)\text{.}
\]
Hence ${\cal M}_{t_{2r+1}^{\otimes d},\nu}={\cal M}_{\tau_{r,d},\nu/\sqrt{2r+1}}$,
and the preceding bound implies that 
\[
R_{n}\left({\cal M}_{t_{2r+1}^{\otimes d},\nu}\right)\ge\widetilde{\kappa}_{r,d}\nu^{d}\frac{\left(\log n\right)^{d}}{n}
\]
for some constant $\widetilde{\kappa}_{r,d}>0$ depending only on
$r$ and $d$. Thus the same lower-bound order holds in the product
Student-$t$ parametrization. 
\end{proof}
\begin{proof}[Proof of Proposition~\ref{prop:GaussianCoreSubmodelRelation}]
Let $\eta_{\nu}$ denote the image of $\eta$ under the map $x\mapsto x/\nu$.
If $f\in{\cal M}_{\varphi,\nu}$, then there exists a probability
measure $G$ on $\mathbb{R}^{d}$ such that 
\[
f=\varphi_{\nu}*G=\left(\gamma_{\nu}*\eta_{\nu}\right)*G=\gamma_{\nu}*\left(G*\eta_{\nu}\right)\text{.}
\]
Since $G*\eta_{\nu}$ is again a probability measure, it follows that
$f\in{\cal M}_{\gamma,\nu}$. This proves the inclusion.

If $\eta=\delta_{a}$ for some $a\in\mathbb{R}^{d}$, then $\eta_{\nu}=\delta_{a/\nu}$
and 
\[
\varphi_{\nu}=\gamma_{\nu}*\delta_{a/\nu}=\gamma_{\nu}\left(\cdot-a/\nu\right)\text{.}
\]
Thus convolution with $\eta_{\nu}$ merely translates the mixing law,
and therefore ${\cal M}_{\varphi,\nu}={\cal M}_{\gamma,\nu}$.

Conversely, suppose that 
\[
{\cal M}_{\varphi,\nu}={\cal M}_{\gamma,\nu}\text{.}
\]
Since $\gamma_{\nu}\in{\cal M}_{\gamma,\nu}$, there exists a probability
measure $G$ on $\mathbb{R}^{d}$ such that 
\[
\gamma_{\nu}=\varphi_{\nu}*G=\gamma_{\nu}*\left(G*\eta_{\nu}\right)\text{.}
\]
Taking Fourier transforms and using the identity 
\[
\widehat{g*\mu}=\left(2\pi\right)^{d/2}\widehat{g}\,\widehat{\mu}
\]
from the Appendix, we obtain 
\[
\widehat{\gamma_{\nu}}=\left(2\pi\right)^{d/2}\widehat{\gamma_{\nu}}\,\widehat{G*\eta_{\nu}}\text{.}
\]
Since $\widehat{\gamma_{\nu}}\left(\xi\right)\neq0$ for every $\xi\in\mathbb{R}^{d}$,
it follows that 
\[
\widehat{G*\eta_{\nu}}\left(\xi\right)=\left(2\pi\right)^{-d/2}\qquad\text{for every }\xi\in\mathbb{R}^{d}\text{.}
\]
The right-hand side is the Fourier transform of $\delta_{0}$. To
justify the corresponding uniqueness step directly, set 
\[
\mu=G*\eta_{\nu}-\delta_{0}\text{.}
\]
Then $\mu$ is a finite signed Borel measure on $\mathbb{R}^{d}$
and $\widehat{\mu}\left(\xi\right)=0$ for every $\xi\in\mathbb{R}^{d}$.
If $\psi\in{\cal C}_{c}^{\infty}\left(\mathbb{R}^{d}\right)$, then
$\widehat{\psi}\in{\cal L}_{1}$ and Fourier inversion gives 
\[
\psi\left(x\right)=\left(2\pi\right)^{-d/2}\int_{\mathbb{R}^{d}}e^{i\left\langle x,\xi\right\rangle }\widehat{\psi}\left(\xi\right)\mathrm{d}\xi\text{.}
\]
Therefore Fubini's theorem yields 
\begin{align*}
\int_{\mathbb{R}^{d}}\psi\left(x\right)\mu\left(\mathrm{d}x\right) & =\left(2\pi\right)^{-d/2}\int_{\mathbb{R}^{d}}\widehat{\psi}\left(\xi\right)\left[\int_{\mathbb{R}^{d}}e^{i\left\langle x,\xi\right\rangle }\mu\left(\mathrm{d}x\right)\right]\mathrm{d}\xi\\
 & =\int_{\mathbb{R}^{d}}\widehat{\psi}\left(\xi\right)\widehat{\mu}\left(-\xi\right)\mathrm{d}\xi\\
 & =\int_{\mathbb{R}^{d}}\widehat{\psi}\left(-\xi\right)\widehat{\mu}\left(\xi\right)\mathrm{d}\xi=0\text{.}
\end{align*}
Hence $\mu$ annihilates ${\cal C}_{c}^{\infty}\left(\mathbb{R}^{d}\right)$.
Since ${\cal C}_{c}^{\infty}\left(\mathbb{R}^{d}\right)$ is uniformly
dense in ${\cal C}_{c}\left(\mathbb{R}^{d}\right)$, it follows that
$\mu$ also annihilates ${\cal C}_{c}\left(\mathbb{R}^{d}\right)$.
Finite signed Borel measures on $\mathbb{R}^{d}$ are determined by
their integrals against ${\cal C}_{c}\left(\mathbb{R}^{d}\right)$,
and therefore $\mu=0$. Thus 
\[
G*\eta_{\nu}=\delta_{0}\text{.}
\]

Now let $X\sim G$ and $Y\sim\eta_{\nu}$ be independent. Then $X+Y=0$
almost surely. Suppose, toward a contradiction, that $G$ is not a
point mass. Then there exist disjoint open balls $B_{1},B_{2}\subset\mathbb{R}^{d}$
such that $G\left(B_{j}\right)>0$ for $j=1,2$. For each $j$, let
\[
\check{B}_{j}=\left\{ -x\text{ | }x\in B_{j}\right\} \text{.}
\]
Since $X+Y=0$ almost surely, we have 
\[
\left\{ X\in B_{j}\right\} \subseteq\left\{ Y\in\check{B}_{j}\right\} \qquad\text{a.s.}
\]
for each $j$. Therefore 
\[
G\left(B_{j}\right)={\rm P}\left(X\in B_{j},Y\in\check{B}_{j}\right)=G\left(B_{j}\right)\eta_{\nu}\left(\check{B}_{j}\right)
\]
by independence, and so $\eta_{\nu}\left(\check{B}_{j}\right)=1$
for $j=1,2$. This is impossible because $\check{B}_{1}$ and $\check{B}_{2}$
are disjoint. Hence $G$ is a point mass. Interchanging the roles
of $G$ and $\eta_{\nu}$ in the same argument shows that $\eta_{\nu}$,
and therefore $\eta$, is also a point mass. This proves the converse. 
\end{proof}

\subsection*{Sources of auxiliary technical results}
Lemma~\ref{lem:convex-approximation}
is Corollary~3.6 of \citet{donahue1997rates}. That corollary gives the same
powers of $m$ with the additional factors
$\left[1+(p-1)\log_{2}m/m\right]^{1/p}$ for $1<p\le2$ and
$\left[1+\log_{2}m/m\right]^{1/2}$ for $2\le p<\infty$; these factors are uniformly
bounded for $m\ge1$ and are absorbed into $C_{p}$ in Lemma~\ref{lem:convex-approximation}.
Lemma~\ref{lem:minkowski-integral-inequality} appears as Theorem~14.16
in \citet{Schilling2017Measures}. Lemma~\ref{lem:smoothness-W1p}
appears as Lemma~12.19 in \citet{Leoni2017SobolevSpaces}. Lemma~\ref{lem:smoothness-Wsp}
is the standard translation estimate stated in Lemma~6.14 of \citet{Leoni2023FractionalSobolev},
with the implicit constant written here as $C_{d,s,p}$. Lemma~\ref{lem:VCBoundMaximalIneq}
follows from the maximal inequalities for pointwise measurable VC-subgraph
classes in \citet[Sec.~2.14.4]{VanderVaartWellner2023}. Lemma~\ref{lem:Assouad}
is Assouad's lemma; see, for example, \citet[Thm.~2.12]{Tsybakov2009Introduction}.
Lemma~\ref{lem:LeastSquaresBasicIneq} is the standard least-squares
basic inequality (cf. \citealp[Lem.~15.1]{klemela2007density}).
Theorem~\ref{thm:MassartNedelecLeastSquares} is a direct specialization
of Theorem~2 of \citet{massart2006risk}.

\bibliographystyle{apalike2}
\bibliography{references}

\end{document}